\documentclass{article}
\font\eufm=eufm10
\def\frak#1{\hbox{\eufm#1}}
\newcommand{\bd}{
\begin{document}}
\newcommand{\ed}{\end{document}}
\newcommand{\be}{\begin{enumerate}}
\newcommand{\ee}{\end{enumerate}}
\newcommand{\bi}{\begin{itemize}}
\newcommand{\ei}{\end{itemize}}
\newcommand{\ba}{\begin{array}}
\newcommand{\ea}{\end{array}}
\newcommand{\sd}{\hspace{0.3ex}\tiny{\rhd\mbox{\hspace{-2ex}}<}\hspace{0.3ex}}
\newcommand{\g}{\frak g}
\newcommand{\ab}{\frak a}
\newcommand{\bb}{\frak b}
\newcommand{\h}{\frak h}
\newcommand{\al}{\frak a}
\newcommand{\n}{\frak n}
\newcommand{\m}{\frak m}
\newcommand{\got}{\frak t}
\newcommand{\ms}{\oplus}
\newcommand{\mt}{\otimes}
\newcommand{\dz}{\wedge}
\newcommand{\D}{\Delta}
\newcommand{\Om}{\Omega}
\newcommand{\om}{\omega}
\newcommand{\oml}{\Omega_L^{1/2}}
\newcommand{\omr}{\Omega_R^{1/2}}
\newcommand{\omh}{\Omega^{1/2}}
\newcommand{\lo}{\lambda_0}
\newcommand{\ro}{\rho_0}
\newcommand{\sA}{\mbox{$\cal A\,$}}
\newcommand{\lra}{\longrightarrow}
\newcommand{\lma}{\Lambda^{max}}
\newcommand{\la}[2]{\Lambda_{#1#2}}
\newcommand{\kad}{ad^{\#}}
\newcommand{\ove}{\overline}
\newtheorem{defi}{Definition}[section]
\newtheorem{tw}[defi]{Theorem}
\newtheorem{prop}[defi]{Proposition}
\newtheorem{lem}[defi]{Lemma}
\newtheorem{re}[defi]{Remark}
\newtheorem{col}[defi]{Corollary}
\newtheorem{ex}[defi]{Examples}
\newcommand{\dow}{\rule{1.6ex}{1.6ex}}
\newcommand{\dowl}{\rule{1ex}{1ex}}
\renewcommand{\arraystretch}{0.8}
\newcommand{\rel}{\mbox{$\,$\rule[0.5ex]{1.1em}{0.2pt}$\triangleright\,$}}
\newcommand{\timh}{\times_h}
\newcommand{\G}{\Gamma}
\newcommand{\Gd}{\G^{(2)}}
\newcommand{\wl}[1]{\vphantom{X}_{#1}{\G}}
\newcommand{\el}{e_L}
\newcommand{\er}{e_R}
\newcommand{\te}{\tilde{e}}
\renewcommand{\baselinestretch}{1.1}
\newcommand{\vs}{\vspace*{0.3\baselineskip}}
\newcommand{\vsm}{\vspace*{-0.3\baselineskip}}
\newcommand{\tran}{\frown\hspace{-2.2ex}|\hspace{1.9ex}}
\newcommand{\R}{\Bbb R}
\renewcommand{\thefigure}{\arabic{section}.\arabic{figure}}
\usepackage{graphicx}
\usepackage{psfrag}
\usepackage{placeins}
\usepackage{amssymbo}
\renewcommand{\textfraction}{0.1}
\renewcommand{\topfraction}{0.89}
\renewcommand{\bottomfraction}{0.8}
\renewcommand{\floatpagefraction}{0.65}

\newcommand{\p}{\psfrag}
\p{g}{ $\Gamma$}
\p{g1}{ $\Gamma'$}
\p{g2}{ $\Gamma''$}
\p{x}{$x$}
\p{x1}{$x_1$}
\p{y}{$y$}
\p{z}{$z$}
\p{z1}{$z_1$}
\p{y1}{$y_1$}
\p{fh}{$f_h$}
\p{fk}{$f_k$}
\p{hbl}{$h_b^L$}
\p{kal}{$k_a^L$}
\p{z2}{$z_2$}
\p{z3}{$z_3$}
\p{y2}{$y_2$}
\p{a}{$a$}
\p{b}{$b$}
\p{c}{$c$}
\p{yp}{$y'$}
\p{zp}{$z'$}
\p{zpp}{$z''$}

\bd
\title{Differential groupoids and $C^*$-algebras.}
\author{Piotr Stachura\\\\
Departement of Mathematical Methods in Physics,\\
Faculty of Physics, University of Warsaw,\\
ul. Hoza 74, 00-682 Warszawa\\
Poland\\
e-mail: stachura@fuw.edu.pl}
\maketitle
\begin{abstract}For a differential groupoid $\G$ we construct a $C^*$ algebra 
$C^*(\G)$ in a way that the correspondance $\G\lra C^*(\G)$ is a covariant
functor from the category of differential groupoids in a sense of S. Zakrzewski
to the category of $C^*$ algebras.\end{abstract}

\section{Introduction}

In this paper we construct a covariant functor from the category
of differential groupoids to the category of $C^*$ algebras in the sense of
\cite{slw2}. However our definition of morphism of differential groupoids is 
different from the standard one i.e. a mapping which satisfies obvious 
compatibility
condition with respect to groupoid structure. Let us argue that there is rather
no hope to construct such a functor with a  standard notion of morphism. 
The main problem can be shown in the discret case so let us
 assume that all sets have discrete topology.

Let $(\G,m,E,s)$ be a groupoid (see the
next section for the notation) and let $\sA(\G)$ denote linear space of
complex functions with compact support. (i.e.  for $f\in \sA(\G)$ we have 
$f(x)\neq0$ only for finite number of $x$'s). There is a natural notion of 
convolution and
star operation in $\sA(\G)$ which make it a *-algebra. Namely 
$$(f_1\,f_2)(x):=\sum_{yz=x}f_1(y)f_2(z)=\sum_{y\in F_l(x)}f_1(y)f_2(s(y)x)=
\sum_{z\in F_r(x)}f_1(xs(z))f_2(z)\,,\,\,f^*(x):=\overline{f(s(x))}$$ 
($F_l(x),F_r(x)$  denote left and right fiber
containing $x$). 
We expect that $C^*(\G)$ will
be a completion of $\sA(\G)$ with respect to some $C^*$-norm.

The ``extremal'' examples of groupoids are sets and groups. For sets the above
multiplication is equal to a pointwise multiplication and for groups it is 
usual convolution. The standard definition of morphism of groupoids reduces
to a mapping if groupoids are sets and to a group homomorphism if they are 
groups. If $h: \G\lra\G'$ is a group homomorphism, we can 
{\em push forward} convolution algebra by a formula 
$(\hat{h}f)(x'):=\sum_{x\in h^{-1}(x')}f(x)$, which defines mapping 
 $\hat{h}:\sA(\G)\lra\sA(\G')$.
But if $h: \G\lra\G'$ is a mapping of sets, functions with pointwise
multiplication can be {\em pulled back} by $(\hat{h}f')(x):=f'(h(x))$. In fact
$\hat{h} f'$ can have non compact support, but this is not a problem, since
we know that it should belong to a (some kind of) multiplier algebra
of $\sA(\G)$ so  $(\hat{h}f')f$ should be in $\sA(\G)$ for any $f\in \sA(\G)$
and certainly this is true. Disregarding the subtelity  in 
this case have $\hat{h}:\sA(\G')\lra \sA(\G)$. We expect that $C^*(h)$ will be 
some extension of $\hat{h}$. And here we are in trouble, since our 
``$C^*$-functor'' is covariant in the first case and contravariant in the 
second. So to achive our goal we {\em need a definition of morphism
between groupoids which reduces to a group homomorphism if groupoids
are groups and to a mapping in the reverse direction if groupoids are sets.}
In particular this suggest that morphisms should be rather {\em relations} 
instead of {\em mappings.}
Such a definition was given in \cite{SZ1} and extended to a differential
setting in \cite{SZ2}.
Let us briefly explain the main idea of the construction, still in the 
discrete  setting. We suggest to look at the next section before the following.

Let $(\G,m,E,s)$ and $(\G',m',E',s')$ be groupoids. A morphism
from $\G$ to $\G'$ is a relation, which satisfies some obvious compatibility
conditions. In particular it turns out that it defines a mapping
$f_h:E'\lra E$ and for each $b\in E'$ a mapping $h_b^R:F_r(f_h(b))\lra F_r(b)$.
For a morphism $h:\G\rel \G'$ and $f\in\sA(\G)$ we define $\hat{h}f$ --a linear
mapping on $\sA(\G')$ by the formula: 
$$((\hat{h} f)f')(x'):=\sum_{x\in F_l(b)} f(x)f'(s'(h_a^R(x))x'),$$ 
where  $a:=\el'(z)\,,\,b:=f_h(a)$. 
By the same formula we define $\pi_h(f) f'$ where 
we view  $f'$ as an element of $L^2(\G')$ -- the Hilbert space of square 
summable functions on $\G'$. Let us also define norms on $\sA(\G)$:
$||f||_l:=sup_{a\in E}\sum_{x\in F_l(a)} |f(x)|\,,\,
||f||_r:=sup_{a\in E}\sum_{x\in F_r(a)} |f(x)|\,,\,$ and 
$||f||:=max\{||f||_l,||f||_r\}.$ It is not difficult to prove the following:
\begin{prop}
a) $(\sA(\G),*,||.||)$ is a normed *-algebra.\\
b) $||\pi_h(f)||\leq ||f||$ and $\pi_h$ is a representation of 
a *-algebra $\sA(\G)$\\
c)$f_3^*(\hat{h}(f_1)f_2)=(\hat{h}(f_1^*)f_3)^*f_2$ for any
$f_1\in \sA(\G)\,,\,f_2,f_3\in\sA(\G')$.\\
d) If $k:\G'\rel\G''$ is a morphism of groupoids then: \\
$\pi_k(\hat{h}(f_1)f_2)f_3=\pi_{kh}(f_1)\pi_k(f_2)f_3$ for any
$f_1\in\sA(\G)\,,\,f_2\in\sA(\G')\,,\,f_3\in L^2(\G'')$--with compact 
support.\\
$\hat{k}(\hat{h}(f_1)f_2)f_3=\hat{kh}(f_1)(\hat{k}(f_2)f_3)
\,$ for any $f_1\in\sA(\G)\,,\,f_2\in\sA(\G')\,,\, f_3\in \sA(\G'')$.
\end{prop}
Using these facts  one can define the $C^*$ 
norm on $\sA(\G)$ by:
$||f||_{C^*}:=\sup ||\pi_h(f)||$ where the supremum is taken over all
morphism $h:\G\rel\G'$. (This is obviously $C^*$-seminorm, but one can show 
that there exists faithfull representation of $\sA(\G)$.) The completion of 
$\sA(\G)$ with respect to this norm
is $C^*-algebra$ of $\G$ and one can see that $\hat{h}$ extents to a 
$C^*(h)\in Mor(C^*(\G),C^*(\G'))$. The extension of this construction to
a differential setting is the main result of the paper.\\
 Of course in the above
case we can also proceed in the standard way: first one can complete 
$(\sA(\G),*,||.||)$ to get a Banach *-algebra and then take envelopping 
$C^*$ algebra. Howewer, in such construction the functoriality is lost and 
moreover it seems that there is no natural, geometric norm on
$\sA(\G)$ in the differential case, so we don't have Banach algebra.

Let us now say a few words about our motivations. One is to get ``geometric
models'' of quantum groups, especially non compact, from double Lie groups.
For a given double Lie group $(G;A,B)$ (see section 3) one can define 
diffeomorphism  $\Psi: G\times G\ni (x,y)\mapsto 
(x (a_L(y))^{-1}, b_R(x (a_L(y))^{-1}) y) \in G\times G$, 
which satisfies pentagonal equation: 
$\Psi_{23}\Psi_{12}=\Psi_{12}\Psi_{13}\Psi_{23}$, where  
$\Psi_{23}:G\times G\times G\ni(x,y,z)\mapsto (x,\Psi(y,z))
\in G\times G\times G$, etc.    
Since  $\Psi$ is a diffeomorphism, it defines (by a push-forward) unitary, 
multiplicative operator  \cite{Baj} on 
$L^2(G\times G)=L^2(G)\mt L^2(G)$ ($L^2(G)$ is a Hilbert space of square 
integrable half  densities on $G$).
This operator is manageable in the sense of Woronowicz \cite{slw},\cite{ps} 
so it  defines a quantum group.

From the other side if $(G;A,B)$ is a double Lie group 
 we can define two 
differential groupoids structures on $G$: $G_A:=(G,m_A,A,s_A)$ and 
$G_B:=(G,m_B,B,s_B)$. 
It turns out that $m_B^T$ is a morphism $G_A\lra G_A\times G_A$ which is 
coassociative: $(m_B^T\times id)m_B^T=(id\times m_B^T)m_B^T$. Applying our
$C^*$ functor we get a coassociative morphism 
$\D\in Mor (C^*(G_A), C^*(G_A\times G_A))$. We expect that $C^*(G_A\times G_A)$
be (a some sort of) $C^*(G_A)\mt C^*(G_A)$. In this way we get one of the
main ingredients of quantum group structure on $C^*(G_A)$. It seems that also
other ingredients as defined in \cite{W-M-N} have natural geometric 
interpretation in the groupoid setting. The details will be presented 
elsewhere. 

For  connections with symplectic geometry and quantisation see Appendix A.

Sections 2 and 3 are of introductory character and are more or less contained 
in \cite{SZ1},\cite{SZ2}. We hope that they make the paper as self contained 
as posiible, keeping its size limited. Section 2 is devoted to algebraic 
structure of groupoid in the language of relations. In section 3 we add 
differential structure. In section 4 - the central of the paper, we 
construct for any
morphism of diffrential groupoids $h$ the mappings $\hat{h},\pi_h$ with the 
properties given above. In section 5 we define a $C^*$ algebra of a 
differential groupoid and prove the functorial properties of the construction.
We also interpret bissections and functions on the set of identities as 
multipliers. In Appendices B and C we relate cocycles on $\G$ with one 
parameter groups and give a construction of weights on $C^*_{red}(\G)$.
In Appendix D we discuss some basic facts about subgroupoids and
Appendix E contains necessary results from  theory of  $C^*$-algebras.\vs

\noindent{\bf Aknewledgment}
The main ideas of this paper are due to \fbox{S. Zakrzewski} who suddenly
died in spring 1998. He knew what had to be done and was convinced that
it could be done.\vs  

\noindent
{\sf Supported by Polish KBN grant no 2 P03A 030 14.}
\section{Groupoids -- algebraic structure}

We begin by recalling some facts about relations.\\
{\em A relation $r$ from $X$ to $Y$} is a triple $r=(R;Y,X)$, where $X$ and
$Y$ are sets and $R$ is a subset of $Y\times X$. $R$ is a {\em graph} of $r$
and we denote it by $Gr(r)$. A relation $r$ from $X$ to $Y$ will be denoted
by $r:X\rel Y$ (note the special type of arrow). 
Relations can be composed, if $s:X\rel Y$ and
$r:Y\rel Z\,$, then a composition $rs$ is a relation from $X$ to $Z$
defined by $Gr(rs):=\{(z,x)\in Z\times X : \exists \, y\in Y\,\,[(z,y)\in
Gr(r)\,\, {\rm and}\,\,(y,x)\in Gr(s)]\}$. \\
We say that the composition $r s$ is {\em simple} iff for any $(z,x)\in Gr(r
s)$ there exists unique $y\in Y$ such that $(y,x)\in Gr(s)$ and $(z,y)\in
Gr(r)$. \\
If $r:X\rel Y$ then its {\em
transposition} is a relation $r^T:Y\rel X$ with
$Gr(r^T):=\{(x,y)\in X\times Y : (y,x)\in Gr(r)\}$.
The cartesian product of relation is also naturaly defined: if $r:X\rel
Y$ and $s:Z \rel T$ then $r\times s:X\times Z\rel Y\times T$ is a
relation with graph $Gr(r\times s):=\{ (y,t,x,z)\in Y\times T\times X\times
Z\,:\,(y,x)\in Gr(r)\,\,and\,\,(t,z)\in Gr(s)\}$.\\
If $r:X\rel Y$ and $A\subset X$ we denote be $r(A)$ the image of $A$ by
$r$: $r(A):=\{y\in Y : \exists\, x\in A \,\,(y,x)\in Gr(r)\}$. Let $\{1\}$ 
denote one point set. 
Now we can formulate the basic definition:
\begin{defi}{\em \cite{SZ1}} Groupoid {\em is a quadrouple $(\G,m,e,s)$ 
where $\G$ is a set, $m:\G\times \G \rel \G$ and $e:\{1\}\rel \G$ are 
relations, $s:\G\lra \G$ is an involution which satisfy:\vsm
\begin{enumerate}
\item[]{\em assiciativity } $m (m\times id)=m (id\times m)$\vsm\vsm,
\item[]{\em identity } $m (e\times id)=m(id \times e)=id$\vsm\vsm,
\item[]{\em inverse } $s m =m (s\times s) \sim$ where 
$\sim:\G\times \G\ni (x,y) \mapsto (y,x)\in \G\times \G$\vsm\vsm,
\item[]{\em strong positivity } for any $x\in \G$ $\emptyset\neq
m(s(x),x)\subset e(\{1\})$
\end{enumerate}\vsm\vsm
Relation $m$ is called {\em multiplication}, $s$ - {\em inverse} and 
$E:=e(\{1\})$ - the set of {\em identities}}.
\end{defi}
\begin{re}
{\em Notice that first three conditions are formally the same as in a group
case but instead of mappings we use relations. The above definition is
equivalent (cf. proposition below) to the "ordinary" definition of groupoid:}
A groupoid is a small category in which every morphism is an isomorphism.
{\em But if we think of groupoid as of category, the natural candidates for
morphism are functors--this is not our point of view, so we prefer the
definition based on relations.} 
\end{re}
\begin{prop} {\em \cite{SZ1}} Let $(\G,m,e,s)$ be a groupoid. Then:\\
a) If $a,b\in E$ then $m(a,b)\neq\emptyset$ iff $a=b$ and in this case
$m(a,a)=a$. \\
b) There exist unique mappings $\el,\er:\G\lra E$ such that
$m(\el(x),x)=x=m(x,\er(x))$ for any $x\in\G$ and $\el(a)=\er(a)=a$ for any
$a\in E$. \\
c)  $m(s(x),x)=\er(x)\,,\,m(x,s(x))=\el(x)$\\
d)  $m(x,y)\neq\emptyset$ iff $\er(x)=\el(y)$\\
e)  $m(x,y)\cap E\neq\emptyset$ implies $y=s(x)$.\\
f)  $m(x,y)$ consists of at most one point.\vs \\
Proof: {\em We prove these statements  here as an exercise in dealing with 
relations.\\
a) From the {\em identity} axiom: if $(x;a,y)\in m$ for some $a\in E$ then
$x=y$ and if $(x;y,b)\in m$ for some $b\in E$ then $x=y$. So for $a,b \in E$,
if $(x;a,b)\in m$ then $x=a=b$. Also for any $x\in\G$ there exists some 
$b\in E$ such that $(x;b,x)\in m$, so for any $b\in E\,\, (b;b,b)\in m$. 
Moreover
from the {\em strong positivity} we have $s(a)=a$ for $a\in E$.\\
b) Suppose that $(x;b_1,x)\in m$ and $(x;b_2,x)\in m$ for $b_1,b_2\in E$. Then
$(x;b_1,b_2,x)\in m(id\times m)=m(m\times id)$ ({\em associativity}). This 
means that $m(b_1,b_2)$ is not empty
and then $b_1=b_2$. This proves the existence and uniqueness of $\el$. In 
the same way one deals with $\er$.\\
c) From {\em strong positivity} we have: $(a;s(x),x)\in m$ for $a\in E$ so 
$(a;s(x),x,\er(x))\in m(id\times m)=m(m\times id)$. From this it follows that 
 exist $(x_1,x_2)$ such
that $(a;x_1,x_2)\in m$ and $(x_1,x_2;s(x),x,\er(x))\in (m\times id)$.
 Then $x_2=\er(x)$ and $a=x_1=\er(x)$. So $m(s(x),x)$ is one point: $\er(x)$.
If $(\er(x);s(x),x)\in m$ then $(\er(x);s(x),\el(x),x)\in m(id\times m)=
m(m\times id)$ so $(\er(x);x_1,x_2)\in m\,$ and  $(x_1,x_2;s(x),\el(x),x)\in 
(m\times id)$. From this we infer that $x_2=x$ and 
$(x_1;s(x),\el(x))\in m$ so 
$\er(s(x))=\el(x)$ and $\el(x)=xs(x)$ since $s$ is an involution.\\
d) If $(z;x,y)\in m$ then $(z;x,\er(x),y)\in m(m\times id)=m(id\times m)$ so 
$m(\er(x),y)$ is not empty and $\er(x)=\el(y)$. Conversely, if $\er(x)=\el(y)$ 
then $(x;x,\er(x)=\el(y))\in m$, so $(x;x,y,s(y))\in m(id\times m)=
m(m\times id)$ and $m(x,y)$ is not empty.\\
e) If $(a;x,y)\in m$ for some $a\in E$ then $(a;\el(x),x,y)\in m(m\times id)=
m(id\times m)$ so $(a;\el(x),x_2)\in m$ and $(x_2;x,y)\in m$ so $a=\el(x)$.
In the same way $a=\er(y)=\el(s(y))$. From this fact: $(s(y);a,s(y))\in m$ and
$(s(y);x,y,s(y))\in m(m\times id)=m(id\times m)$. We have: 
$(s(y);x_1,x_2)\in m$ and $(x_1,x_2;x,y,s(y))\in (id\times m)$ so 
$(s(y);x,\el(y)=\er(x))\in m$ so $s(y)=x$.\\
f) Let $(x_1;x,y)\in m$ and $(x_2;x,y)\in m$. In this situation 
$\el(x_1)=\er(s(x_2))$ so $(z;x_1,s(x_2))\in m$ for some $z$. So 
$(z;x,y,s(x_2))\in m(m\times id)$ and $(z;x,x_4)\in m$ and 
$(x_4;y,s(x_2)\in m$. But then $(x_4;y,s(y),s(x))\in m(id\times m)$ and 
$x_4=s(x)$. So $z\in E$ and $x_1=x_2$.}\\ 
\dow
\end{prop}
Now we explain our notation.\\
The set of composable pairs will be denoted by $\Gd
:=m^T(\G)=\{(x,y)\in\G\times\G : \er(x)=\el(y)\}$.  From the statements
c) and f) of the above proposition follows that $m$ restricted to $\Gd$ is a
surjective mapping on $\G$. 
The set of identities $E$ will also be denoted by $\G^0$. If it doesn't lead
to any confusion we write $x=x_1 x_2$ instead of $(x;x_1,x_2)\in Gr(m)$.
For $A,B\subset \G$ let $A B:=\{a b\,:\,a\in A\,,\,b\in B\}$.
If $r:X\rel Y$ we also write $(y,x)\in r$ instead of $(y,x)\in
Gr(r)$. \\
For $x\in\G\,$,  by $F_l(x)$ and $F_r(x)$ we denote left and right fibers
containing $x$ i.e. $F_l(x):=\el^{-1}(\el(x))$ and $F_r(x):=\er^{-1}(\er(x))$.
If $a\in E$ we also write $\vphantom{\G}_a\G:=F_l(a)$ and $\G_a:=F_r(a)$.
Clearly $\vphantom{\G}_a\G\cap \G_a$ is a group. \\
For $x\in\G$ let $L_x:F_l(\er(x))\ni z\mapsto x z\in F_l(x)$ and
$R_x:F_r(\el(x))\ni z\mapsto z x\in F_r(x)$ denote left and right translation
by $x$. They are bijections.
\begin{ex}\label{ex1}
a) Sets. {\em If $X$ is a set then $(X,d^T,X,id)$, where $d: X\lra
X\times X$ is a diagonal mapping, is a groupoid. Conversely, any groupoid such
that $m^T$ is a mapping is of this type.}\vs\\
b) Groups. {\em  If $G$ is a group then $(G,m_G,\{e\},s)$, where $m_G,s$
are group multiplication and group inverse, is a groupoid and any groupoid
for which $m$ is a mapping is a group.}\vs\\
c) Pair groupoids. {\em  Examples a) and b) are ``extremal''examples of 
groupoids. The ``middle'' and the simplest are pair groupoids. Let $X$ be a 
set. We put: $\G:=X\times X\,,\,
Gr(m):=\{((x,y);(x,z),(z,y))\,:\,x,y,z\in X\}\,,\\
s(x,y):=(y,x)$ and $\G^0:=\{(x,x)\,,\,x\in X\}$. 
Then $(\G,m,\G^0,s)$ is a groupoid. }\vs\\
d) Equivalence relations. {\em If $R\subset X\times X$ is an equivalence
relation on $X$, then $(R,m,\G^0,s)$ where $m,\G^0,s$ are as above is a
groupoid.}\vs\\ 
e) Transformation groupoids. {\em Let a group $G$ acts on a set $X$. We
denote the action by $G\times X\ni(g,x)\mapsto g x\in X$. Define $\G:=G\times
X$, $s(g,x):=(g^{-1},g x)$, $E:=\{e\}\times X$ and $m$ by
 $$Gr(m):=\{((g_1 g_2,x);(g_1,g_2 x),(g_2,x))\,:\,g_1,g_2\in G\,,\,x\in
X\}\subset\G\times\G\times\G.$$ 
Then $(\G,m,E,s)$ is a groupoid.}\vs\\
f) Double groups. {\em \cite{SZ1} Let $(G;A,B)$ be  double group i.e.
$A,B\subset G$ are subgroups, $A\cap B=\{e\}$ and $G=A B$.
In this situation each element of $G$ can be written uniquely as: 
$g=a_L(g)b_R(g)=b_L(g)a_R(g)$. This decomposition defines four mappings:  
$a_R,a_L:G\lra A$  and $b_R,b_L:G\lra B$. 
Let $m_A:G\times G\rel G$ be a relation defined by
$$m_A(g_1,g_2):=\left\{\ba{ll} g_1b_R(g_2)=b_L(g_1)g_2 &
if\;a_R(g_1)=a_L(g_2)\\ 
\emptyset & otherwise \ea \right.$$
Let  $s_A:G\ni g\mapsto (b_L(g))^{-1} a_L(g)\in G.$ 
Then $G_A:=(G,m_A,A,s_A)$ is a groupoid. The same holds for
$G_B:=(G,m_B,B,s_B)$. \vs\\
g) And many, many more. See e.g. \cite{McK}}
\end{ex}
The cartesian product of groupoids is defined in a natural way. For 
groupoids  $(\G_i,m_i,E_i,s_i)\,,\,i=1,2\,$, their
cartesian product, which we denote simply by 
$\G_1\times \G_2$,  is a groupoid \\
$(\G_1\times\G_2,(m_1\times m_2)(id\times\sim\times id),
E_1\times E_2,s_1\times s_2)$.\vs

\noindent
{\bf Morphisms of groupoids}.
\begin{defi} {\em \cite{SZ1} Let $(\G,m,e,s)$ and $(\G',m',e',s')$ be 
groupoids.} A
morphism from $\G$ to $\G'$ {\em is a relation $h:\G\rel\G'$ such that:\vs\\
1. $hm=m'(h\times h)$\hspace{0.1\textwidth}2. $hs=s'h$\hspace{0.1\textwidth}
3. $he=e'$.}
\end{defi}
Note the following proposition \cite{SZ1}:
\begin{prop}\label{mor}
Let $h:\G\rel \G'$ be a morphism of groupoids. Then:\\
a) The compositions in the definition above are simple.\vs\\
b) Let the relation $h_0:E\rel E'$ be defined by $Gr(h_0):=Gr(h)\cap(E'\times
E)$. \\
Then $(\er'\times \er)Gr(h)=(\el'\times \el)Gr(h)=Gr(h_0)$ and
$f_h:=h_0^T$ is a mapping.\vs\\
c) Let $b\in E'$ and $a:=f_h(b)$.  Let us define  two relations  $h_b^R:
F_r(a)\rel F_r(b)$ and $h_b^L: F_l(a)\rel F_l(b)$  by 
$Gr(h_b^R):=Gr(h)\cap(F_r(b)\times F_r(a))\,$ and $\,
Gr(h_b^L):=Gr(h)\cap(F_l(b)\times F_l(a))$. Then $h_b^R\,,\,h_b^L$
are mappings.\\
\dow
\end{prop}
Morphisms can also be characterised in terms of mapping. This is the contents of the following: 
\begin{prop} \label{mormap} 
Any morphism $h:\G\rel\G'$
determines and is uniquely determined 
by mappings \\
$f:E'\lra E$ and $g:\G\times_f E'\lra \G'$, where
$\G\times_fE':=\{(x,e')\in\G\times E'\,:\,\er(x)=f(e')\}$ which satisfy
conditions: \\
a) $\el\er^{-1}(f(E'))=f(E')$ (then also $\er\el^{-1}(f(E'))=f(E')$)\\
b) $\er'g(x,e')=e'$\\
c) $s'g(x,e')=g(s(x),\el'g(x,e'))$\\ 
d) $\forall \,x_1,x\in\G\,\, [(x,e')\in \G\times_fE'$ and
$\er(x_1)=\el(x)]\Rightarrow \,g(x_1 x,e')=g(x_1,\el'g(x,e'))\, g(x,e')$ \vs\\
Proof: {\em First, we show that such two mappings define morphism of
groupoids. Let  a relation $h$ be given by the graph:
$Gr(h):=\{(g(x,e'),x)\,:\,(x,e')\in \G\times_fE'\}$. Notice that
$$e'=\er'g(f(e'),e')= s'(g(f(e'),e'))\,
g(f(e'),e')=g(f(e'),\el'g(f(e'),e'))\,g(f(e'),e')= g(f(e'),e').$$ First
eqality follows from statement b), third from c) and the last one from d). 
This shows that relation $h$ satisfies: $h E=E'$.\\
We have the following sequence of equivalences: 
$$(y,x)\in h s\iff (y,s(x))\in h\iff y=g(s(x),e')\iff
s'(y)=g(x,\el'g(s(x),e'))\iff$$ 
$$(s'(y),x)\in h\iff (y,x)\in s'h.$$ 
In this way $h s=s' h$.

Let $(y;x_1,x_2)\in h m$. Now we have:
$$(y;x_1,x_2)\in h m\iff [\er(x_1)=\el(x_2)\,\,{\rm and}\,\,
(y,x_1 x_2)\in h ]
\iff [\er(x_1)=\el(x_2)\,\,{\rm and}\,\,y=g(x_1 x_2,e')].$$ 
It follows that $\er(x_2)=f(e')$, so for $y_2:=g(x_2,e')$ and 
$y_1:=g(x_1,\el'(y_2))$ we
have: $(y_2,x_2)\in h\,,\,(y_1,x_1)\in h$ and $y=y_1 y_2$. From this:
$(y;x_1,x_2)\in m'(h\times h)$.\\
Conversely, for $(y;x_1,x_2)\in m'(h\times h)$ we have $y=y_1 y_2$ with
$\er'(y_1)=\el'(y_2)$ and $y_1=g(x_1,e_1')\,,\,y_2=g(x_2,e_2')$. Now
$s'(y_2)=g(s(x_2),\el'(y_2))=g(s(x_2),\er'(y_1))$, so $\el(x_2)=\er(x_1)$ and
$x_1,x_2$ are composable. So we get  that $m'(h\times h)\subset h m$.

Now the ``determines'' part. If $h: \G\rel\G'$ is a morphism, 
define $f:=f_h$ and $g(x,e'):=h_{e'}^R(x)$.
Then a) and b) follows directly from prop. \ref{mor}. 
Let us show c). 
$$y=s'g(x,e')\iff s'(y)=g(x,e')\iff [(s'(y),x)\in
h\,\,{\rm and}\,\,\er's'(y)=\el'(y)=e'] \iff$$ 
$$ [(y,s(x))\in h\,\,{\rm and}\,\,\el'(y)=e'].$$ 
But $\er'(y)=\el' g(x,e')$, so 
$y=g(s(x),\el'g(x,e'))$.\\ 
d) Let $\er(x)=f(e')\,,\,\er(x_1)=\el(x)\,,\,y:=g(x,e')\,\,{\rm and}
\,\, z:=g(x_1 x,e')$. This means that  $y$ is given by the conditions:
$(y,x)\in h\,,\,\er'(y)=e'$ and $z$ by $(z,x_1 x)\in h\,,\,\er'(z)=e'$. 
So $z, s'(y)$ are composable and\\
$(z s'(y),x_1 x, s(x))\in m'(h\times h)=h m$. From this we have 
$(z s'(y),x_1)\in h\,,\,\er'(z s'(y))=\el'(y)\, ,\, 
z s'(y)=g(x_1,\el'g(x,e'))$ and $z=g(x_1,\el'g(x,e'))\, y$. \\
\dow
}\end{prop}
{\newcommand{\tg}{\mbox{$\tilde{\G}$}}
\newcommand{\ts}{\mbox{$\tilde{s}$}}
\renewcommand{\te}{\mbox{$\tilde{E}$}}
\newcommand{\tm}{\mbox{$\tilde{m}$}}
Let $h: \G\rel \G'$ be a morphism and $f,g$ be as above. Denote 
$\tg:=\G\times_fE'$ and define 
$$\ts:\tg\ni(x,b)\mapsto 
(s(x),\el'g(x,b))\in\tg\,,\,\, \te:=E\times_f E'$$ and a relation
$\tm:\tg\times\tg\rel\tg$ by: 
$$Gr(\tm):=\{(x_1 x_2,b_2;x_1,\el'g(x_2,b_2),x_2,b_2):
\er(x_1)=\el(x_2)\,,\,(x_2,b_2)\in \tg\}$$.
\begin{lem}
$(\tg,\tm,\te,\ts)$ is a groupoid.\vs\\
Proof: {\em Let us first check that $\ts$ is an involution: 
$$\ts\ts(x,b)=\ts(s(x),\el'g(x,b))=(x,\el'g(s(x),\el'g(x,b)))=
(x,\er's'g(s(x),\el'g(x,b)))=(x,\er'g(x,b))=(x,b).$$
1. $\tm(\tm\times id)=\tm(id\times\tm)$. \\
Compute the left hand side:
$$(x_1,b_1;x_2,b_2,x_3,b_3,x_4)\in \tm(\tm\times id)\iff $$
$$\exists\,
(x_5,b_5),(x_6,b_6):[(x_1,b_1;x_5,b_5,x_6,b_6)\in\tm\,\,{\rm and}\,\,
(x_5,b_5,x_6,b_6;x_2,b_2,x_3,b_3,x_4,b_4)\in (\tm\times id)]\iff$$
$$ \exists \,
(x_5,b_5):[(x_1,b_1;x_5,b_5,x_4,b_4)\in\tm\,\,{\rm and}\,\,
(x_5,b_5;x_2,b_2,x_3,b_3)\in\tm]\iff$$
$$ [\er(x_2)=\el(x_3)\,{\rm and}\,b_2=\el'g(x_3,b_3)\,
{\rm and}\,(x_1,b_1;x_2 x_3,b_3,x_4,b_4)\in\tm].$$
So $Gr(\tm(\tm\times id))=$
$$=\{(x_2 x_3 x_4,b_4;x_2,\el'g(x_3,\el'g(x_4,b_4)),
x_3,\el'g(x_4,b_4),x_4,b_4):\er(x_2)=\el(x_3)\,,\,\er(x_3)=\el(x_4)\,,\,
(x_4,b_4)\in\tg\}.$$
And the right hand side:
$$(x_1,b_1;x_2,b_2,x_3,b_3,x_4)\in \tm(id\times \tm)\iff$$
$$ \exists\,(x_5,b_5),(x_6,b_6):[(x_1,b_1;x_5,b_5,x_6,b_6)\in\tm\,\,{\rm and}
\,\,(x_5,b_5,x_6,b_6;x_2,b_2,x_3,b_3,x_4,b_4)\in (id\times \tm)]\iff $$
$$[\er(x_3)=\el(x_4)\,{\rm and}\,
b_3=\el'g(x_4,b_4)\,{\rm and}\,(x_1,b_1;x_2,b_2,x_3 x_4,b_4)\in\tm].$$
It follows that $Gr(\tm(id\times\tm ))=$
$$=\{(x_2 x_3 x_4,b_4;x_2,\el'g(x_3 x_4,b_4),
x_3,\el'g(x_4,b_4),x_4,b_4)\,:\,\er(x_2)=\el(x_3)\,,\,\er(x_3)=\el(x_4)\,,\,
(x_4,b_4)\in\tg\}.$$ 
But $\el'g(x_3 x_4,b_4)=
\el'(g(x_3,\el'g(x_4,b_4))\,g(x_4,b_4))=\el' g(x_3,\el'g(x_4,b_4))$. 
In this way
$\tm(\tm\times id)=\tm(id\times \tm).$\vs\\
2. $\tm(\te\times id)=\tm(id\times \te)=id$.\\ 
If $(x,b;a_1,b_1,x_2,b_2)\in \tm$
for some $(a_1,b_1)\in \te$ then $x=a_1 x_2$ and $b=b_2$ so $(x,b)=(x_2,b_2)$.
Conversely for any $(x,b)\in \tg$ we have $(x,b;\el(x),\el'g(x,b),x,b)\in\tm$. 
So $\tm(\te\times id)=id.$ In the same way one shows that 
$\tm(id\times \te)=id$.\vs\\
3. $\ts\tm=\tm(\ts\times\ts)\sim$.\\
The left hand side: 
$$(x,b;x_1,b_1,x_2,b_2)\in \ts\tm\iff
(s(x),\el'g(x,b);x_1,b_1,x_2,b_2)\in\tm\iff $$
$$\er(x_1)=\el(x_2)\,,\,
x=s(x_2)s(x_1)\,,\,\el'g(x,b)=b_2\,,\,b_1=\el'g(x_2,b_2)$$ 
but $b=\er'g(x,b)=\el'g(s(x),\el'g(x,b))=\el'g(x_1 x_2,b_2)=
\el'g(x_1,\el'g(x_2,b_2)).$ \\
In this way 
$$Gr(\ts\tm)=\{(s(x_2)s(x_1),\el'g(x_1,\el'g(x_2,b_2));x_1,
\el'g(x_2,b_2),x_2,b_2)\,:\,\er(x_1)=\el(x_2)\,,\,(x_2,b_2)\in\tg\}.$$
And the right hand side:
$$(x,b;x_1,b_1,x_2,b_2)\in \tm(\ts\times\ts)\sim\iff
(x,b;s(x_2),\el'g(x_2,b_2),s(x_1),\el'g(x_1,b_1))\in\tm\iff$$
$$x=s(x_2)s(x_1)\,,\,b=\el'g(x_1,b_1)\,,\, \er(x_1)=\el(x_2)\,,\,
\el'g(x_2,b_2)=\el'g(s(x_1),\el'g(x_1,b_1))=\el's'g(x_1,b_1)=b_1.$$ 
So we have:
$Gr(\ts\tm)=Gr(\tm(\ts\times\ts)\sim)$.\vs\\
4. $(s(x) x,b;s(x),\el'g(x,b),x,b)\in\tm$ for any $(x,b)\in \tg$ and 
$(s(x) x,b)\in \te$.}\\
\dowl
\end{lem}
Consider relations $h_1:\G\rel\tg$ and $h_2:\tg\rel\G'$ defined by:
$Gr(h_1):=\{(x,b;x)\,:\,(x,b)\in \tg\}\,,\,Gr(h_2):=\{(g(x,b);x,b)\,:\,
(x,b)\in\tg\}$. Clearly we have $h=h_2h_1$, moreover $h_1$ is a morphism
from $\G$ to $\tg$ and $h_2$ is a morphism from $\tg$ to $\G'$. For $h_1$ the 
mappings between fibers are bijective, and  $f_{h_2}$ is
bijective mapping. In this way we have the following structure:
\begin{prop} If $h:\G\rel\G'$ is a morphism of groupoids, then exists groupoid
$\tg$, morphisms $k:\G\rel\tg$ and $l:\tg\rel\G'$ such that $h=lk$ and:\\
a) For each $a\in\te$ the mappings $k_a^R$ and $k_a^R$ are bijections.\\
b) $l$ is a mapping from $\tg\lra \G'$ which is bijective when resticted to  
$\te$.\\
\dow
\end{prop}} 
Groupoids together with just defined morphisms form a category as the 
following lemma states.
\begin{lem}{\em \cite{SZ1}} Let $h: \G\rel \G'$ and $k:\G'\rel \G''$ be
morphisms of groupoids. Then $h$ and $k$ have simple composition and 
$k h$ is a morphism from $\G$ to $\G''$.\\
\dowl
\end{lem}
\begin{ex}\label{exmor}
{\em 
a) If $X$ is a set and $(\G,m,e,s)$ is a groupoid then any morphism
$h:X\rel \G$ is equal $f^T$ for some mapping $f: E\lra X$.\vs\\
b) If $G,H$ are groups then morphisms from $G$ to $H$ are just group
homomorphisms. \vs \\
c) If $X$ is a set and $G$ is a group then morphisms $X\rel G$ are points
of $X$. \vs\\
d) If $(G;A,B)$ is a double group then $m_B^T: G_A\rel G_A\times G_A$ and
$m_A^T: G_B\rel G_B\times G_B$ are morphisms of groupoids. \cite{SZ1}\vs\\
e) For any groupoid $\G$ the mapping $\G\ni x\mapsto (\el(x),\er(x))\in
E\times E$ is a morphism from $\G$ to the pair groupoid $E\times E$. We
denote this relation by $\te$.\vs\\
f)  The relation $l:\G\rel\G\times \G$ given by: $(x,y;z)\in
Gr(l)\Leftrightarrow (x;z,y)\in Gr(m)$ is a morphism from $\G$ to the pair
groupoid $\G\times\G$. It is called {\em left regular representation}.
\cite{SZ1}}
\end{ex}
\begin{re} {\em
The above defined morphisms differ from the standard one, but
later on we will see, that this definition is proper for defining the algebra
of groupoid and the functorial properties of the construction. Also we want 
to point out that our definition {\em is not a generalisation} of a usual
definition. Below we show
that our  morphisms are related to actions of groupoids on sets.}
\end{re}
\begin{defi} {\em \cite{Wein}  Let $(\G, m, e, s)$ be a groupoid, $Y$ be a set
and $\mu:\, Y\lra \G^0$ be a mapping. Denote $\G\times_{\mu} Y:=\{(x,y)\in
\G\times Y\,:\,\er(x)=\mu(y)\}$.} \\
The (left) action of $\G$ on $Y$ {\em is a
mapping $\phi:\G \times_{\mu} Y\ni(x,y)\mapsto \phi(x,y)\in Y$ which
satisfy conditions:\\ 
a) $\mu\phi(x,y)=\el(x)$\\
b) $\phi(x_1 x_2,y)=\phi(x_1,\phi(x_2,y))$ (i.e. if one side of the equality is
defined the other also and are equal)\\
c) $\phi(\mu(y),y)=y$}
\end{defi}
Now let $\G$ acts on $Y$. Put $f:=\mu$ and $g:\G\times_{f} Y\ni(x,y)\mapsto
(\phi(x,y),y)\in Y\times Y$. Then it is easy to see that these mappings
satisfy the conditions given in Prop. \ref{mormap}, so they determine
morphism from $\G$ to the pair groupoid $Y\times Y$. Conversely, if $h:\G\rel
Y\times Y$ is a morphism then putting: $\mu:=f_h$ and
$\phi(x,y):=\el'h_y^R(x) $ we get the action of $\G$ on $Y$. 
Also for any morphism $h:\G\rel\G'$ the mappings $\mu:=f_h \el':\G'\lra E$
and $\phi(x,x'):=h_{a'}^R(x) x'$ where $a':=\el'(x')$ define action of $\G$
on $\G'$. \vs\vs

\noindent
{\bf Bissections}
\begin{defi} A bissection $B$ is a subset of $\G$ such that:
$\el\mid_B:B\lra \G^0\;and \;\er\mid_B:B\lra \G^0$ are
bijections. \end{defi}
The set of bissections of $\G$ will be denoted by ${\cal B}(\G)$.
Bissections can also be characterized as follows:
\begin{lem}
A subset $B\subset\G$ is a bissection $\Leftrightarrow\,$ $B s(B) = s(B) B =
\G^0$.\vs \\
Proof: {\em $\Rightarrow$ Let $x=b_1 s(b_2)$ for some $b_1,b_2\in B$. Then
$\er(b_2)=\er(b_1)$ and, since $B$ is a bissection, $b_2=b_1$ and $x\in 
\G^0$, so $B s(B)\subset \G^0$. Moreover for any $x\in\G^0$ we can find $b\in
B$ with $x=\el(b)=b s(b)$ so $\G^0\subset B s(B)$ and $B s(B)=\G^0$. In the
same way we have $s(B) B=\G^0$.\\ 
$\Leftarrow$ Suppose that for $b_1,b_2\in B$ we have $\er(b_1)=\er(b_2). $
Then $b_1$ and $s(b_2)$ are composable, so $b_1 s(b_2)\in\G^0$ and $b_1=b_2$.
The same holds for the left projection. So $\el|_B$ and $\er|_B$ are
injective. But for any $x\in \G^0$ we can find $b_1,b_2\in B$ with
$x=\el(b_1)=\er(b_2)$. }\\
\dowl
\end{lem}
For a bissection $B$ let $Bx:=B\{x\}$ and $xB:=\{x\} B$.
The following proposition is easy to prove:
\begin{prop} Let $B$ be a bissection.\\
a) $L_B:\G\ni x\mapsto B x\in\G$ and $R_B:\G\ni x\mapsto x B\in\G$ are
bijections.\\ 
b) $L_B\,(R_B)$ preserves right (left) fibers.\\
c) $L_B(F_l(x))=F_l(B x)\,,\,R_B(F_r(x))=F_r(x B).$\\
d) $B (x y)=(B x) y\,,\,(x y) B=x (y B)$.\\
e) If $B,C\in{\cal B}(\G)$ then $B C \in {\cal B}(\G)$ and ${\cal B}(\G)$ is
a group. \\
f) $L_B\,(R_B)$ is a left (right) action of ${\cal B}(\G)$ on $\G$.\\
\dow
\end{prop}
\begin{ex} {\em a) For any groupoid the set of identities is a bissection.\\
b) If $\G$ is a group then bissections are just group
elements. \vs\\
c) If $\G:= X\times X$ is a pair groupoid then any bissection is of the form
$B:=\{ (f(x),x)\,:\,x\in X\}$ for some bijection $f:X\lra X$.\vs\\
d) If $(G;A,B)$ is a double group then for any $b\in B$ the sets $bA,Ab$ are
bissections of $G_A$, the mapping $B\ni b\mapsto bA\in {\cal B}(\G)$ is a 
group homomorphism.} 
\end{ex}
\begin{prop}
Let $h:\G\rel\G'$ be a morphism and $B$ a bissection of $\G$ then the set
$h(B)$ is a bissection of $\G'$.\vs\\
Proof: {\em Let us take any $a'\in E'$ and let $E\ni a:=f_h(a')$. Then
there exist unique points $x,y\in B$ with $\el(x)=a=\er(y)$, so there exist
unique $x',y'\in\G'$ such that $(x',x)\in Gr(h)\,,\,\el'(x')=a'$ and
$(y',y)\in Gr(h)\,,\,\er'(y')=a'$.} \\
\dow
\end{prop}
\begin{re} {\em 
Since bissections defines bijections of $\G$ they acts, in this purely 
algebraic context, on $\sA(\G)$. This action commutes with right 
multiplication. Since morphisms acts also on bissections we can expect that
they are unitary multipliers on  $C^*(\G)$. And this is true. Later on we will
see that also in the differential setting bissections can be interpreted as 
multipliers.}
\end{re}
\begin{re} {\em
One can think of groupoids as of some generalisation of groups and treat 
groupoid elements in the same way as group elements. But this analogy can be 
misleading since for groups bissections are just elements. So group elements
have some properties ``because they are groupoid elements'' and others 
``because they are bissections''.}
\end{re}

\section{Differential groupoids}
From now on, when we use the word {\em manifold} without any comments, we
mean Hausdorff, finite dimensional, smooth manifold with a countable basis of
neighbourhoods. {\em Submanifold} is a nonempty, embedded submanifold ( with
the relative topology).

{\em A differentiable relation} $r:X\rel Y$ is a triple $r=(R;Y,X)$ such that
$X,Y$ are manifolds and $R$ is a submanifold in $Y\times X$.

If $r=(R;Y,X)$ is a differentiable relation then its {\em tangent lift} is a
relation $Tr:TX\rel TY$ with a graph $Gr(Tr):=TGr(r)$. A {\em phase lift}
of $r$ is a relation $Pr:T^*X\rel T^*Y$ such that: 
$$(\alpha,\beta)\in
Gr(Pr) \iff  <\alpha,u>=<\beta,v>\,\,{\rm  for\,\, any}\,\, 
(u,v)\in T_{(y,x)}Gr(r),$$ 
where $y:=\pi_Y(\alpha) \,,\,x:=\pi_X(\beta)$ and 
$\pi_X,\pi_Y$-are the cannonical projections in the cotangent boundles.

We say that relations $r:X\rel Y$ and $s:Y\rel Z$ are {\em transverse} iff
$Tr,Ts$ and $Pr,Ps$ have simple composition, and $sr$ is a differentiable
relation. Such situation will be donoted by $r\tran s$.

Let us also recall that {\em a differentiable reduction} is a differentiable
relation $r:X\rel Y$ of the form $r=fi^T$, where $i:C\lra X$ is an inclusion
map of the submanifold $C\subset X$ and $f:C\lra Y$ is a surjective
submersion. 

\begin{defi}{\em \cite{SZ2}} A differential groupoid {\em $(\G,m,e,s)$ is
a groupoid such that $\G$ is a manifold, $m$ is a differentiable reduction,
$e$ is a differentiable relation, $s$ is a diffeomorphism and the following
transversality conditions  hold: 
$m\tran (m\times id)\,,\,m\tran (id\times m)\,,\,m\tran (e\times
id)\,,\,m\tran (id\times e)$.}
\end{defi}
It can be shown \cite{SZ2} that in this situation $\el,\er$ are
submersions. 
\begin{ex} {\em a) Examples \ref{ex1} a)-e) after obvious smoothness
conditions are differential groupoids.\vs\\
b)}  Double Lie groups. {\em We say that a double group $(G;A,B)$ is a double
Lie group iff $G$ is a Lie group and $A,B$ are closed subgroups of $G$. Then
$G_A,G_B$ are differential groupoids.\vs\\
c)} Tangent and cotangent bundles. {\em If $X$ is a manifold then
$(TX,+,X,-)$ and $(T^*X,+,X,-)$ are differential groupoids. More generally,
if $(P,X)$ is a vector bundle then it is a groupoid in a natural way 
$(P,+,X,-)$.\vs\\
d)} Tangent and phase lifts of differential groupoids {\em \cite{SZ2}.
If $(\G,m,e,s)$ is a differential groupoid then $(T\G,Tm,Te,Ts)$ and
$(T^*\G,Pm,Pe,-Ps)$ are differential groupoids. If $\G:=(X,d^T,X,id)$ is
a manifold groupoid then its tangent lift $T\G=(TX,d_{TX}^T,TX,id)$ is 
again manifold groupoid but its cotangent lift $P\G=(T^*X,+,X,-)$ is a 
cotangent bundle with usual groupoid structure.\vs\\
e) If $\G=(G,m,e,s)$ is a Lie group, then its tangent lift is a Lie group
$TG$.  But the phase lift is $T^*G$ as a transformation groupoid:
$T^*G=G\times \g^*$ with a coadjoint action.}
\end{ex}

\noindent
{\bf Morphisms of differential groupoids.}
\begin{defi}{\em  \cite{SZ2} Let $\G\,,\,\G'$ be  differential
groupoids and $h:\G\rel \G'$ a differentiable relation which is a
(algebraic) morphism of groupoids. Then $h$ is a morphism of differential
groupoids iff $m'\tran (h\times h)$ and $h\tran e$.}
\end{defi}
\begin{prop} {\em \cite{SZ2}\label{mordif}
}
If $h:\G\rel \G'$ is a morphism of differential groupoids and $f_h:=h_0^T$
then: \\
a) $f_h:E'\lra E$ is a smooth mapping.\\
b) $\G *_h E':=\{(x,b)\in \G\times E'\,:\,\el(x)=f_h(b)\}$ is a
submanifold of $\G\times E'$ (and of $\G\times\G'$). 
\\
c) The mapping: $Gr(h)\ni(y,x)\mapsto (x,\el'(y))\in \G *_h E'$ is a
diffeomorphism.\\
\dow
\end{prop}
In the next lemma we collect the properties of various sets and mappings
associated with morphisms of differential groupoids, which will be used later
on.
\begin{lem}\label{rfakty} Let $h:\,\G\rel \G'$ be a morphism of
differential groupoids. 
\begin{itemize}
\item[a)] $\G *_h\G':=\{(x,y)\in \G\times \G'\,:\,\el(x)=f_h(\el'(y))\}$ is a
submanifold of $\G\times\G'$. 
\item[b)] $\G\timh E':=\{(x,b)\in \G\times E'\,:\,\er(x)=f_h(b)\}$ is a
submanifold of $\G\times E'$. 
\item[c)] $\G\timh \G':=\{(x,y)\in \G\times \G'\,:\,\er(x)=f_h(\el'(b))\}$ is a
submanifold of $\G\times \G'$. 
\item[d)] Let $(x,y)\in \G\timh\G'$ and $b:=\el'(y)$, then the mappings: \\
$m_h: \G\timh\G'\ni
(x,y)\mapsto m'(h_b^R(x),y)\in\G'$ and $\pi_2 :\G\timh\G'\ni (x,y)\mapsto y\in
\G'$\\
 are  surjective submersions.
\item[e)] The mapping $\got_h: \G\timh\G'\ni (x,y)\mapsto (x,m_h(x,y))\in \G*_h\G'$ is a
diffeomorphism. 
\item[f)] For $b\in E'$, $h_b^R:\G\supset F_r(f(b)) \lra F_r(b)\subset\G'$ 
and $h_b^L:\G\supset F_l(f(b))\lra F_l(b)\subset\G'$ are smooth mappings.
\item[g)] $\tilde{\pi_2} :\G*_h\G'\ni (x,y)\mapsto y\in \G'$ is surjective
submersion.  
\item[h)] The sets $\G\timh\G_a'\,,\,\G*_h\G_a'$ are submanifolds of
$\G\timh\G'$ and $\G*_h\G'$. Also statements d),e),g) remains true after
suitable restriction of the corresponding mappings.
\end{itemize}
Proof: {\em 
\begin{itemize}
\item[a)] Consider $id\times \el':\G\times\G'\lra \G\times E'$ - this is
smooth submersion and $\G*_h\G'=(id\times\el')^{-1}(\G*_h E')$. The assertion
follows from second item of the previous proposition.
\item[b)] $\G\timh E'=(s\times id)(\G*_h E')$ but $(s\times id)$ is a
diffeomorphism so again we use second item of the previous proposition.
\item[c)] Write $\G\timh\G'=(id\times\el')^{-1}(\G\timh
E')$ and use b).
\item[d)] From the third statement of the previous proposition we know that
the mapping $h^L:\G*_h E'\ni(x,b)\mapsto h_b^L(x)\in\G'$ is smooth. Now the
mapping $h^R:\G\timh E'\ni(x,b)\mapsto h_b^R(x)\in\G'$ is the composition of :
$s'h^L(s\times id): \G\timh E'\lra \G'$ so it is smooth. Now $m_h$ is the
composition:\\
$m_h:\,\G\timh\G'\stackrel{id\times\el'\times id}{\lra}\G\times
E'\times\G' \stackrel{h^R\times
id}{\lra}\G'^{(2)}\stackrel{m'}{\lra}\G'\\
\hspace*{3em}(x,y)\mapsto(x,\el'(y),y)\mapsto(h_{\el'(y)}^R(x),y)\mapsto
m'(h_{\el'(y)}^R(x),y)$\\ 
and is a smooth mapping. Since for any $y\in \G'$ $m_h(f_h(\el'(y)),y)=y$ it
is clear that $m_h$ is surjective.\\
The mapping $m_h$ is illustrated on the picture below.\\
\begin{figure}[tbhp]\label{mh}
\psfrag{hbr}{$h_b^R$}
\p{mh}{$m_h(x,y)$}
\centering
\fbox{\includegraphics[height=0.25\textheight,width=0.7\textwidth]
{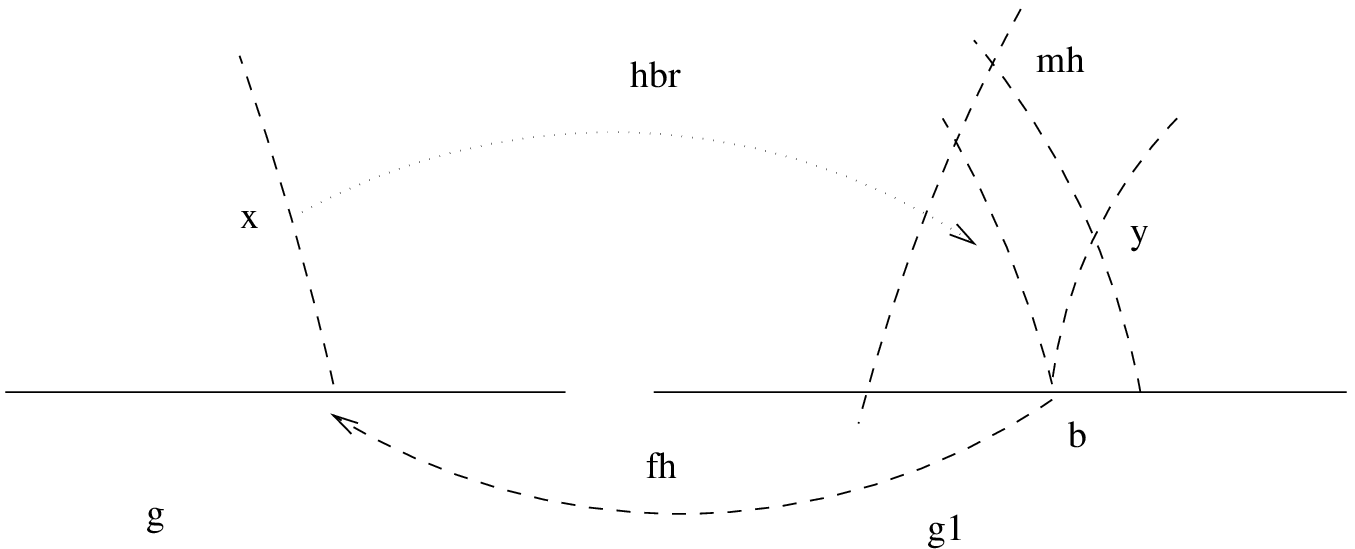}}
\caption{}
\end{figure}
\FloatBarrier
Now let $m_h(x_0,y_0)=z_0$ and let $z(t)\,,\,z(0)=z_0$
be a curve through $z_0$. Then $f_h(\el'(z(t)))$ is a curve in $E$ through
$\el(x_0)$. It can be lifted to $x(t)$ - a curve through $x_0$ with
$\el(x(t))=f_h(\el'(z(t)))$. Since $h^L$ is smooth
$w(t):=h^L(x(t),f_h(\el'(z(t))))$ is a curve in $\G'$ with
$\el'(w(t))=\el'(z(t))$ and $y(t):=m'(s'(w(t)),z(t))$ is a curve through
$y_0$. One can check that $m_h(x(t), y(t))=z(t)$. The statement about $\pi_2$
is obvious.
\item[e)] From the previous statement $\got_h$ is smooth and open. It is clear
that it is an immersion. So it is enaugh to show that it is surjective.
Let $(x,z)\in \G*_h\G'$ with $b:=\el'(z)$. Then it is easy to see that
$(x,z)=\got_h(x,s'(h_b^L(x)) z)$.
\item [f)] It follows from the fact that $h^L,h^R$ are smooth.
\item [g)] It is obvious.
\item[h)] The proofs are  simmilar to the proofs of points d), e) and g).
\end{itemize} }\dowl
\end{lem}
\begin{ex}\label{exmor2}
{\em
a) If $f:X\lra Y$ is a smooth mapping then $T^*f$ considered as a relation:
$ T^*Y\rel T^*X$ is a morphism of differential groupoids. The same is true for 
$(Tf)^T:TY\rel TX$. Note that here $TX,TY$ are considered as manifold 
groupoids, not vector bundle groupoids. Unless $f$ is a local diffeomorphism,
$(Tf)^T$ is not a morphism of $(TY,+,Y,-)$ and $(TX,+,X,-)$. \\
b) Let a Lie group $G$ acts on a manifold $X$. We form a transformation
groupoid $\G:=G\times X$. Let $Y$ be a manifold and $h:\G\rel Y\times Y$ be  
a morphism to the pair groupoid of $Y$. Consider the smooth mapping 
$\Phi: G\times Y\ni (g,y)\mapsto \el'h^R(g,f_h(y);y)\in Y$. Then $\Phi$
defines an action of $G$ on $Y$. Moreover $f_h$ is equivariant i.e. $g
f_h(y)=f_h \Phi(g,y)$. Conversely, if $G$ acts on $X$ and $Y$ with
equivariant mapping $f:Y\lra X$ then $h$ defined by: $Gr(h):=\{
(gy,y;g,f(y)))\,:\, y\in Y\,,\,g\in G\}$ is a morphism from $\G$ to $Y\times
Y$. \vs\\
c) Let $X,Y$ be manifolds and $\G:=X\times X$, $\G':=Y\times Y$ be
corresponding pair groupoids. Then using Props. \ref{mormap} and \ref{mordif}
one can see that any morphism $h: \G\rel \G'$ is determined by a smooth
surjection $f: Y\lra X$ and a smooth mapping $g: X\times Y\lra Y$ which
satisfy for any $x,x_1\in X\,,\,y\in Y\,$ conditions:\\
a) $f g(x,y)=x\,$, b) $g(f(y),g(x,y))=y\,$, c) $g(x,y)=g(x,g(x_1,y))$. \\
Then $Gr(h):=\{(g(x,y),y;x,f(y))\,:\,x\in X\,,\,y\in Y\}.$ 
From b) it follows that
$g$ is a surjection. From c) and b): $g(f(y),y)=g(f(y),g(x,y))=y$. If
$x_0=f(y_0)$ and $x(t)$ is a curve through $x_0$ then $g(x(t),y_0)$ is a
curve through $y_0$. Since $x(t)=f g(x(t),y_0)$ -- $f$ is a submersion.
Choose some $x_0\in X$ and let $Z:=f^{-1}(x_0)$-- this is a submanifold of
$Y$. We claim that the mapping: $\phi: X\times Z\ni(x,y)\mapsto g(x,y)\in Y$
is a diffeomorphism. It is clear that this mapping is smooth.  If
$g(x_1,y_1)=g(x_2,y_2)$ then from a) $x_1=x_2$. So
$y_1=g(f(y_1),g(x,y_1))=g(f(y_2),g(x,y_2))=y_2$ and $\phi$ is an injection.
For $y\in Y$ we have: $y=g(f(y),g(x_0,y))$ but $g(x_0,y)\in Z$ so it 
is a surjection. From a): $g(\dot{x},\dot{y})=0\Rightarrow \dot{x}=0$. If
$y(t)$ is a curve in $Z$ through $y$ then
$y(t)=g(x_0,y(t))=g(x_0,g(x,y(t)))$. From this follows that $\phi$ is an
immersion. So it is a diffeomorphism. In this way for any morphism $h: \G\rel
\G'$ there exists diffeomorphism $\phi: Y\lra X\times Z$ and
 $Gr((\phi\times\phi) h)=\{(x,z,x_1,z;x,x_1)\,:\,x,x_1\in X\,,\,z\in Z\}$.}
\end{ex}
The next proposition shows that differential groupoids with the above defined
morphisms form a category.
\begin{prop} {\em \cite{SZ2}}
Let $\G,\G',\G''$ be differential groupoids and let $h:\G\rel\G'$ and
$k:\G'\rel \G''$ be morphisms. Then $h\tran k$ and $k h:\G\rel\G''$ is
a morphism. \\
\dow
\end{prop}
A submanifold $B\subset\G$ is a {\em bissection} iff $\el|_B$ and $\er|_B$ are
diffeomorphisms. If $h:\G\rel \G'$ is a morphism of differential groupoids
and $B$ is a bisection then $h(B)$ is the image of $E'$ by the mapping:
$$E'\ni a'\mapsto (f_h(a'),a')\mapsto ((\er|_B)^{-1}f_h(a'),a')\mapsto
h^R((\er|_B)^{-1}f_h(a'),a')\in\G'$$ also $h(B)$ is the image of $E'$ by 
the mapping:
$$E'\ni a'\mapsto (f_h(a'),a')\mapsto ((\el|_B)^{-1}f_h(a'),a')\mapsto
h^L((\el|_B)^{-1}f_h(a'),a')\in \G'.$$ 
These mappings are smooth sections of
the projections $\er'$ and $\el'$ respectively. So $h(B)$ is a bissection of
$\G'$.

\section{A *-algebra of a differential groupoid.}

In this section we introduce a *-algebra of a differential groupoid. The
way we do it is rather non-standard and at the first look may be regarded
as too complicated, nevertheless it is very convenient for the further
development. 

Let $\G$ be a differential groupoid and let $\omh(\el),(\omh(\er))$ be the
smooth bundle of complex, half densities along the left (right) fibers of
$\G$. Following A. Connes  \cite{Con} our basic object is  the linear space 
of compactly supported smooth sections 
 of the bundle $\omh(\el)\mt\omh(\er)$. We denote this space by $\sA(\G)$. 
Its elements
will be called {\em bidensities} and usually denoted by $\om$. So
$\om(x)=\lambda(x)\mt\rho(x)\in \omh T^l_x\G\mt\omh 
T^r_x\G$, where we used notation:
$T^l_x\G:=T_x(F_l(x))\,,\,T^r_x\G:=T_x(F_r(x))$. In the following we also
write $\oml(x):=\Om^{1/2} T^l_x\G$ and $\omr(x):=\Om^{1/2} T^r_x\G$. \\
{\em We also use the following notation: if $M,N$ are manifolds,  
$F:M\lra N$ and $\Psi$ is some geometric object on $M$ which can be 
pushed-forward by $F$, then we denote the push-forward of $\psi$ simply 
by $F\psi$. What it really means will be clear from the context.}

The groupoid inverse induces the star operation on $\sA$ as follows 
$$\om^*(x)(v\mt w):=\overline{\om(s(x))(s(w)\mt s(v))}\,,\,v\in
\Lambda^{max}T^l_x\G\,,\,w\in \Lambda^{max}T^r_x\G$$ 
(for any vector space
$V$ by $\Lambda^{max} V$ we denote the maximal non zero 
exterior power of $V$). This
is well defined antylinear involution. (since $s$ is an involutive
diffeomorphism which interchanges left and right fibers).

We are going to equip $\sA$ with multiplication, which gives it a * - algebra
structure. First we will show that with any morphism
$h:\, \G\rel\G'$ is associated a mapping $\hat{h}: \sA(\G)\lra L\sA(\G')$
($L\sA(\G')$ denote linear endomorphisms of $\sA(\G')$), which ``well behave''
with respect to composition of morphisms and *-operation. Then putting $h=id$
we get algebra structure on $\sA(\G)$ and $\hat{h}$ will became *-algebra
homomorphism from $\sA(\G)$ to the left algebraic multipliers of $\sA(\G')$.
Before this we  define some special sections of $\omh(\el)\mt\omh(\er)$
which are very convenient for computations. \vs

\noindent
{\bf *-invariant bidensities.}\vs
 
Since left(right) translations are diffeomorphisms of left(right) fibers, we
can define left(right) invariant sections of $\omh(\el)(\omh(\er))$,
namely a section $\lambda$ is left invariant iff for any
$(x,y)\in\Gd\,\,\lambda(x y)(x v)=\lambda(y)(v)\,\,v\in\Lambda^{max}T^l_y\G$.
In the same way are defined right invariant half densities. Any left
invariant half density is determined by its value on $\G^0$ and conversely
any section of $\omh(\el)|_{\G^0}$ can be uniquely extended to left invariant
half density on $\G$.

So let $\tilde{\lambda}$ 
be a non-vanishing, real, half density on $\G^0$ along the left fibers. (one
constructs such a density by covering $\G^0$ with maps submitted to submersion
$\el$ and using appropriate partition of unity to glue them together) We
define: 
$$\lambda_0(x)(v):=\tilde{\lambda}(\er(x))(s(x) v)\,,\, v\in\lma T^l_x\G,$$
then $\lambda_0$ is a left invariant, non vanishing section of $\omh(\el)$.  
Now $\tilde{\rho}:=\tilde{\lambda} s $ is non vanishing, real, half
density on $\G^0$ along the right fibers, and $\ro$ defined by:
$\ro(x)(v):=\tilde{\rho}(\el(x)(v s(x))\,,\,v\in\lma T^r_x\G$ is a right
invariant, non vanishing, real half density along the right fibers. 
Let $\omega_0:=\lo\mt \ro$, then this is real, non vanishing
bidensity. \\
{\em From now on the symbol $\om_0$ will always mean bidensity
constructed in this way.} 
When $\om_0$ is choosen any element $\om\in\sA(\G)$ can be written uniquely
as $\omega=f\,\omega_0$ for some smooth, complex function $f$ with compact
support. Note the following:
\begin{lem} \label{f-gw} 
If $\om=f \om_0$ then $\om^*=f^*
\om_0$ where $f^*(x):=\overline{f(s(x))}.$\vs\\
Proof: {\em  $\om^*(x)(v\mt w):=\overline{\om(s(x))(s(w)\mt s(v))}\,,\,v\in
\lma T^l_x\G\,,\,w\in \lma T^r_x\G$.
$$\overline{\om(s(x))(s(w)\mt s(v))}=
\overline{f(s(x))}(\lo(s(x))\mt\ro(s(x)))(s(w)\mt s(v))=$$
$$=f^*(x) \lo(\el(x))(x s(w))\ro(\er(x))(s(v) x)=
f^*(x) \ro(\el(x))(w s(x))\lo(\er(x))(s(x) v)=$$
$$=f^*(x) \lo(x)(v)\ro(x)(w)=f^*(x)\om_0(x)(v\mt w).$$}
\dowl
\end{lem} 
\begin{re} {\em Choosing $\lo$ in fact we choose some left Haar system in the 
sense of \cite{Ren} on our groupoid. But all our constructions and in 
particular our $C^*$ algebra are independent of this choice. }
\end{re}
\noindent
{\bf Action of groupoid morphisms on bidensities}.\vs

Now we are going to construct for a morphism of diffrential groupoids $h$ 
the mapping $\hat{h}$.\\
Let  $h:\G\rel\G'$ be a morphism of differential groupoids. Then from
lemma \ref{rfakty}  we know that:\\
1. The set $\G\timh\G_a':=\{(x,y)\in
\G\times\G'\,:\,\er(x)=f_h(\el'(y)),\er('y)=a\} $ 
is a submanifold of $\G\timh\G'$. \\
2. The mapping: $\pi_2: \G\timh\G_a'\ni(x,y)\mapsto y\in \G_a'$ 
is a surjective submersion and $\pi_2^{-1}(y)$ is diffeomorphic to
$F_r(f_h(\el'(y)))$. \\
3. The mapping $$\got_h:\G\timh\G_a'\ni (x,y)\mapsto (x,m_h(x,y))\in
\G*_h\G_a':=\{(x,y)\in \G\times \G_a'\,:\, \el(x)=f_h(\el'(y))\}$$ 
is a diffeomorphism.\\ 
4.  $\tilde{\pi}_2:\G*_h\G_a'\ni(x,y)\mapsto y\in\G_a'$ 
is a surjective submersion and $\tilde{\pi}_2^{-1}(y)$ is diffeomorphic to
$F_l(f_h(\el'(y)))$. \vs

Before we go further, let us recall some facts about densities. Let $V$ be a 
finite dimensional vector space. For $p\geq 0$ we denote the linear space of 
complex $p$-densities on $V$ by $\Om^p(V)$. If $V=V_1\ms V_2$ and $\nu_1,\nu_2$
are $p$-densities on $V_1,V_2$ then the formula $(\nu_1\mt\nu_2)
(v_1\dz v_2):=\nu_1(v_1)\nu_2(v_2)$ for $v_1\in\lma V_1\,,\,
v_2\in \lma V_2$ defines isomorphism $\Om^p(V)=\Om^p(V_1)\mt\Om^p(V_2)$.
Also we have $\Om^p(V)=\Om^p(V_1)\mt \Om^p(V/V_1)$ defined by choosing some 
$V_2\subset V$ complementary to $V_1.$ The isomorphism does not depend from
the choice made. In this way if $F:V\lra W$ is a linear surjection, we have
cannonical isomorphism $\Om^p(V)=\Om^p(ker\, F)\mt\Om^p(W)$. This fact is 
constantly used in the following. Now we go back to groupoid morphisms.    

Let $(x,y)\in\G\timh\G_a'$ and $\got_h(x,y)=:(x,z)\,,\,b:=\el'(z)$. Due to
the point 2. we have
isomorphism: $i_1:\omr(x)\mt\omr(y)\lra \omh T_{(x,y)}(\G\timh\G_a')\,$. 
From point 3) $\got_h:\omh T_{(x,y)}(\G\timh\G_a')\lra \omh
T_{(x,z)}(\G*_h\G_a')$ is an isomorphism and from 4) 
$i_2:\oml(x)\mt\omr(z)\lra \omh T_{(x,z)}(\G*_h\G_a')$ is an isomorphism.
In this way we get equality: 
$(i_2)^{-1}\got_h i_1(\rho_x\mt\rho_y)=:\lambda_x\mt\rho_z$ for some
$\lambda_x\mt\rho_z\in\oml(x)\mt\omr(z)$. Moreover
the mapping $F_l(y)\ni u \mapsto h_{\el'(y)}^R(x) u \in F_l(z)$ is a
diffeomorphism, so 
for $\lambda_y\in\oml(y)$, we have $h_{\el'(y)}^R(x) \lambda_y\in\oml(z)$.\\
Now let $\om=\lambda\mt\rho\in\sA(\G)\,$ 
$\om'=\lambda'\mt\rho'\in\sA(\G')$. Then:
$(i_2)^{-1} \got_h i_1(\rho(x)\mt\rho'(y))=:
\tilde{\lambda_x}\mt\tilde{\rho_z}$ 
and $h_{\el'(y)}^R(x)\lambda'(y)=:\tilde{\lambda_z'}$. So the  expression:
$[\lambda(x)\tilde{\lambda_x}]\mt\tilde{\lambda_z'} \mt\tilde{\rho_z} $
defines  one-density on $F_l(f_h(b))$ 
with values in one dimensional vector space $\oml(z)\mt\omr(z)$. \\
Let us define  
$$(\hat{h}(\om) \om')(z):=
\int_{F_l(f_h(b))}[\lambda\tilde{\lambda}]\mt\tilde{\lambda_z'}
\mt\tilde{\rho_z}. $$
Choose:  $\om_0=\lambda_0\mt\rho_0$, $\om_0'=\lambda_0'\mt
\rho_0'$. Then $\om=f_1\,\om_0\,,\,\om'=f_2\,\om_0'$ and
$(i_2)^{-1} \got_h i_1 (\rho_0(x)\mt\rho_0'(y))=:
t_h(x,y)\lambda_0(x)\mt\rho_0'(z)$ for some smooth, nonvanishing function
$t_h:\G\timh\G'\lra R$ and $h_{\el'(y)}^R(x) \lambda_0'(y)=\lambda_0'(z)$. \\
We get the explicit expression:
$$(\hat{h}(\om)\om')(z):=\left[\int_{F_l(f_h(b))}\lo^2(x) f_1(x)
t_h(x,y)f_2(y)\right]\,\om_0'(z)=:(f_1*_hf_2)(z)\,\om_0'(z),$$ 
where $y$ is defined
by $\got_h(x,y)=(x,z)$, i.e. $y=s'(h_b^L(x))z$.

\begin{figure}[tbhp]
\centering
\fbox{\includegraphics[height=0.25\textheight,width=0.9\textwidth]
{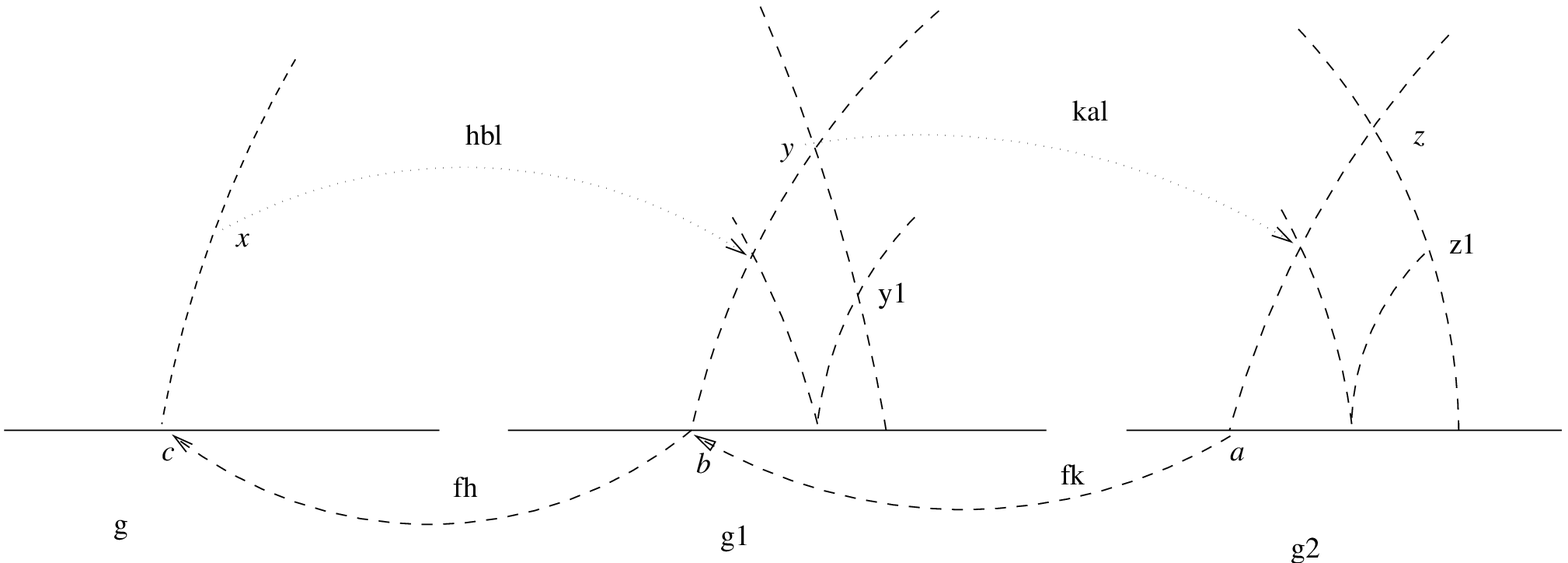}}
\caption{\label{funkt}}
\end{figure}

The next proposition is crucial for the construction, it describes how the
mapping $\hat{h}$ behaves with respect to a composition of morphisms. We cannot
simply write:
$\hat{k}(\hat{h}\om)=\hat{k h}\om$ since the left hand side is not defined. 
Instead of this equality we prove the other one, which is formally the 
same as for morphism of $C^*$-algebras. (see Appendix E)
\begin{prop}\label{funktor}
Let $h:\G\rel\G'\,,\,k:\G'\rel \G''$ be morphisms of differential groupoids.
Then  
$$\hat{k}(\hat{h}(\om_1)\om_2)\om_3=\hat{kh}(\om_1)(\hat{k}(\om_2)\om_3)
\,\,{\rm  for\,\, any\,\,} 
\om_1\in\sA(\G)\,,\,\om_2\in\sA(\G')\,,\, \om_3\in \sA(\G'').$$
\end{prop}
\noindent
{\em Proof:} Choose $\om_0:=\lo\mt\ro\,,\,\om_0':=\lo'\mt\ro'\,,\,
\om_0'':=\lo''\mt\ro''$ and write 
$\om_1=f_1\,\lambda_0\mt\rho_0\,,\,\om_2=f_2\,\lambda_0'\mt\rho_0'\,,\,
\om_2=f_3\,\lambda_0''\mt\rho_0''$. \\
Let
$z\in\G''\,,\,a:=\el''(z)\,,\,b:=f_k(a)\,,\,c:=f_h(b)=f_hf_k(a)=f_{kh}(a).$
\\ 
Compute the left hand side of the equality:\\
$(\hat{k}(\hat{h}(\om_1)\om_2)\om_3)(z)= ((f_1*_h f_2)*_k f_3)(z)\,\om''(z)\,$
and
$$((f_1*_h f_2)*_k f_3)(z)=\int_{F_l(b)}\lo'^2(y)\,(f_1*_hf_2)(y)
\,t_k(y,z_1)\,f_3(z_1)=$$
$$=\int_{F_l(b)} \lo'^2(y)\,\left[\int_{F_l(c)} \lo^2(x)\, f_1(x)\,t_h(x,y_1)\,
f_2(y_1)\right]\, t_k(y,z_1)\,f_3(z_1)=$$
$$=\int_{F_l(b)\times F_l(c)} (\lo'^2(y)\mt\lo^2(x))\,[f_1(x)\,
f_2(y_1)\,f_3(z_1)\,t_h(x,y_1)\,t_k(y,z_1)],$$ 
where $y_1,z_1$ are given by:
$\got_h(x,y_1)=(x,y)\,,\,\got_k(y,z_1)=(y,z)$.\\ 
The situation is illustrated on the figure \ref{funkt}.

\begin{figure}[tbhp]
\psfrag{khal}{$(kh)_a^L$}
\psfrag{ka2l}{$k_{a_2}^L$}\psfrag{a2}{$a_2$}\psfrag{b2}{$b_2$}
\centering
\fbox{\includegraphics[height=0.25\textheight,width=0.9\textwidth]
{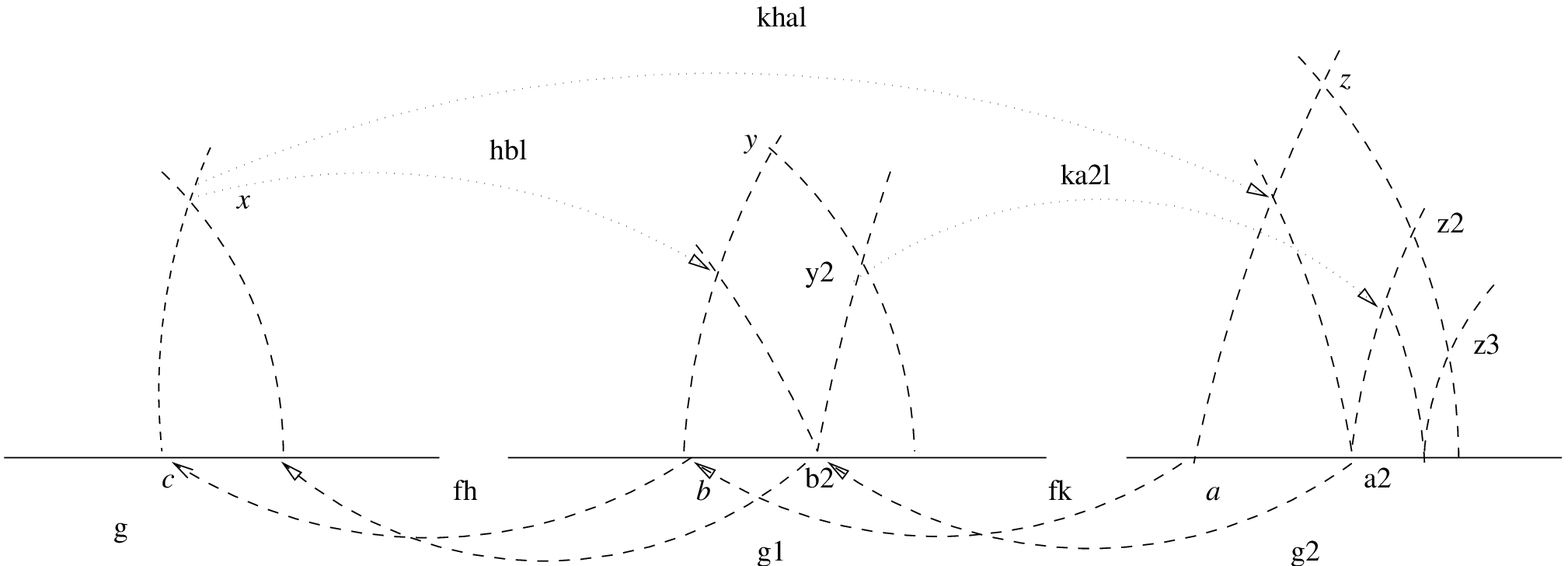}}
\caption{\label{funktor2}}
\end{figure}

\noindent
The right hand side: $(\hat{kh}(\om_1)(\hat{k}(\om_2)\om_3))(z)=
(f_1*_{k h}(f_2*_kf_3))(z)\, \om''(z)\,$ and 
$$(f_1*_{k h}(f_2*_kf_3))(z)=\int_{F_l(c)}\lo^2(x)\,
f_1(x)\, t_{kh}(x,z_2)\,(f_2*_kf_3)(z_2)=$$
$$=\int_{F_l(c)}\lo^2(x)\,f_1(x)\,t_{kh}(x,z_2)\,
\int_{F_l(b_2)}\lo'^2(y_2)\,f_2(y_2)\,
t_k(y_2,z_3)\,f_3(z_3),$$ 
where $y_2,z_2,z_3$ are given by: $\got_{kh}(x,z_2)=(x,z)
\,,\,\got_k(y_2,z_3)=(y_2,z_2) $ and $b_2:=f_k(\el''(z_2))$.\\ 
For fixed  $x\in F_l(c)$ the mapping: $F_l(b_2)\ni
y_2\mapsto m_h(x,y_2)\in F_l(b)$ is a diffeomorphism.
Using this fact, above integral is equal: 
$$\int_{F_l(b)\times F_l(c)}(\lambda_0'^2(y)\mt\lambda^2_0(x))
\,[f_1(x)\,f_2(y_2)\,f_3(z_3)\,
t_{kh}(x,z_2)\,t_k(y_2,z_3)],$$
where $y_2,z_2,z_3$ are defined by:
$\got_{kh}(x,z_2)=(x,z)\,,\,\got_h(x,y_2)=(x,y)
\,,\,\got_k(y_2,z_3)=(y_2,z_2).$\\ 
The situation is illustrated on the figure \ref{funktor2}.
\FloatBarrier
Comparing this expression with the left hand side, we see that $y_2=y_1$. \\
We prove that $z_3=z_1$. We have: $z_3=s''(k_{a_2}^L(y_2))
z_2\,,\,y_2=s'(h_b^L(x)) y$ and $z_2=s''((kh)_a^L(x)) z$. So
$z_3=s''(\tilde{z}) z$ for 
$\tilde{z}:=(kh)_a^L(x)\,k_{a_2}^L(y_2)=(kh)_a^L(x)\,k_{a_2}^L[s'(h_b^L(x))
y]$ and 
$z_1=s''(k_a^L(y)) z$. We have $\el''(\tilde{z})=\el''(z)=a$ and 
$\el''(k_a^L(y))=a$ so it is enough to show that
$(\tilde{z},y)\in Gr(k)$.\\
Since $((kh)_a^L(x),x)\in Gr(k h)$ there exists
unique (composition of morphisms is simple!) $\tilde{y}$ satisfying
$((kh)_a^L(x),\tilde{y})\in Gr(k)$ and $(\tilde{y},x)\in Gr(h)$ but then
$\el'(\tilde{y})=f_k(a)=b$ so 
$\tilde{y}=h_b^L(x)$. Then $(\tilde{z};h_b^L(x),y_2)\in
 Gr(m''(k\times k)=Gr(k m')$ and $(\tilde{z}, h_b^L(x) y_2)\in Gr(k)$. 
But $h_b^L(x) y_2=h_b^L(x) s'(h_b^L(x)) y=y$ and $(\tilde{z},y)\in Gr(k)$.\\
So to prove equality it remains to show that:
$t_h(x,y_1)\,t_k(y,z_1) = t_{kh}(x,z_2)\,t_k(y_1,z_1)$ . This is done in the
following: 
\begin{lem} \label{app} 
Let $(x,y,z)\in\G\times\G'\times\G''$ 
satisfy: $\er(x)=f_h(\el'(y))\,,\, \er'(y)=f_k(\el''(z))$ and let $y',z'$ be
defined by: $\got_h(x,y)=(x,y')$ and $\got_k(y,z)=(y,z')$. Then
$$t_h(x,y) t_k(y',z)=t_{k h}(x,z') t_k(y,z).$$
Proof: {\em The proof is given in the end of this section.}
\end{lem}
\noindent \dow
\begin{col}\label{algebra} 
The algebra structure on
$\sA(\G)$. {\em Take 
$h=id:\G\rel\G$, then $\G\timh\G=\Gd$ and put: $\om_1
\om_2:=\hat{id}(\om_1)\om_2$. Due to the above proposition, 
this product is associative. Chosen $\om_0$ we can write
$\om_1=f_1\om_0\,,\,\om_2=f_2\om_0$. Let us show that in this situation
$t_{id}\equiv 1$.\\
We have
$\G\times_{id}\G_c=\{(x,y)\in\G\times\G\,:\,\er(x)=\el(y)\,,\,\er(y)=c\}\,$,
$\G*_{id}\G_c=\{(x,y)\in\G\times\G\,:\,\el(x)=\el(y)\,,\,\er(y)=c\}$ and 
$\got_{id}(x,y)=(x,x y)$. Let $(x,y)\in\G\times_{id}\G_c\,,\,z:=x y
\,,\,\er(x)=:b\,,\,\el(x)=\er(y)=:a\,$,
$V:=T_{(x,y)}\G\times_{id}\G_c\,,\,W:=T_{(x,z)}\G*_{id}\G_c$. It is easy to
see that the following decompositions hold: $V=V_1\oplus V_2\,,\,W=W_1\oplus
W_2$ where $V_1:=\{(\dot{b} x,0)\,:\,\dot{b}\in T_b^r\G\}\,,\,V_2:=\{(x
s(\dot{a}),\dot{a} y)\,:\,\dot{a}\in T_a^r\G\}\,,\,W_1:=\{(\dot{b} x,\dot{b}
z)\}\,,\, W_2:=\{(x s(\dot{a}),0)\}$. Moreover $\got_{id}(\dot{b}
x,0)=(\dot{b} x,\dot{b} z)$ and $\got_{id}(x s(\dot{a}),\dot{a} y)=(x
s(\dot{a}),0)$. Using the definition of $t_{id}$ and $\om_0$ we get the
desired result. So the explicit formula for the product is:
$\om_1\om_2=:(f_1*f_2) \om_0$ and  $$(f_1* f_2)(x):=\int_{F_l(x)}
\lo^2(y) f_1(y) f_2(s(y) x)=\int_{F_r(x)} \ro^2(y) f_1(x s(y)) f_2(y).$$
The second equality follows from the fact that right and left fibers are
diffeomorphic by $s$. It is clear then it is smooth bidensity with support
contained in $m(supp\,\om_1,\,supp\,\om_2)$ and this is compact set.}
\end{col}
\begin{col} {\em We define the set $LM\sA(\G)$ of} left algebraic multipliers
of $\sA(\G)$ {\em as those linear mappings from $L\sA(\G)$ which commutes
with right multiplication. This is evidentely algebra. Then from the above
proposition follows that $\hat{h}$ is a homomorphism from $\sA(\G)\lra
LM\sA(\G')$.} 
\end{col}

That the multiplication is compatible with the star operation (i.e $\sA$ is
in fact * - algebra) can be shown directly, but again it will follow from
more general facts which we prove later on.\\
The above defined multiplication is non degenerate as is easy seen from the
following: 
\begin{lem} \label{nondegmul}
For any $\om\in\sA(\G)$: $\om^*\om=0\iff\om=0$.\vs \\
Proof: {\em Choose $\om_0$ and write $\om=f \om_0$. Then for $a\in \G^0$ we
have: $$f^**f(a)=\int_{\G_a}\ro^2(x) f^*(s(x)) f(x)=\int_{\G_a}\ro^2(x)
\overline{f(x)} f(x).$$ 
From this formula the statement is clear.}\\
\dowl
\end{lem} 

It seems that there is no natural, geometric, norm on $\sA(\G)$ but one can
introduce the family of useful norms ``indexed'' by $\om_0$. \cite{Ren} So
choose $\lo$ and write $\om=f\om_0$. We define quantities:
$$||\om||_l:=\sup_{a\in\G^0}\int_{F_l(a)}\lo^2\,|f|\,,\,\,\,
||\om||_r:=\sup_{a\in\G^0}\int_{F_r(a)}\ro^2 \,|f|\,,\,\,\, 
||\om||:=max\{||\om||_l,||\om||_r\}.$$ 
(We do not explicitely 
write the dependence from $\lo$ to make notation simpler.) 
The next lemma shows that the definitions are correct.
\begin{lem}\label{norm}
The functions $||.||_l\,,\,||.||_r\,,\,||.||$ are norms and give 
$\sA$ the structure of
normed algebra. Moreover $||\om^*||_l=||\om||_r$ so $||\om^*||=||\om||$. (in
fact $(\sA,*,||,||)$ is normed *-algebra as we will see later on.)\vs\\
Proof: {\em It is obvious that $||.||_l\,,\,||.||_r$ and consequently $||.||$
are norms. Also from lemma \ref{f-gw} we immediately have
$||\om||_l=||\om^*||_r$. \\
Let us show that 
$||\om_1\om_2||_l\leq ||\om_1||_l||\om_2||_l$. 
As usual we write: $\om_1=f_1\,\om_0\,,\,\om_2=f_2\,\om_0$. Then 
$||\om_1\om_2||_l=sup_{a\in\G^0}\, \int_{\wl{a}} \lo^2(x)|f_1*f_2|(x)$. 
$$\int_{\wl{a}} \lo^2(x)|f_1*f_2|(x)=\int_{\wl{a}}
\lo^2(x)\left|\int_{\wl{a}}\lo^2(y)f_1(y)f_2(s(y) x)\right|\leq \int_{\wl{a}}
\lo^2(x)\int_{\wl{a}}\lo^2(y)\,|f_1(y)||f_2(s(y)x)|=$$
$$=\int_{\wl{a}\times\wl{a}}\lo^2(x)\mt\lo^2(y)\,|f_1(y)||f_2(s(y) x)|=
\int_{\wl{a}\times\G_a}\lo^2(x)\mt\ro^2(y)\,|f_1(s(y))||f_2(y x)|=$$
$$=\int_{\G_a}\ro^2(y)\,|f_1(s(y))|\int_{\wl{a}}\lo^2(x)\,|f_2(y x)|= 
\int_{\G_a}\ro^2(y)\,|f_1(s(y))|\int_{F_l(y)}\lo^2(z)\,|f_2(z)|\leq$$
$$\leq ||\om_2||_l\int_{\G_a}\ro^2(y)\,|f_1(s(y))|=
||\om_2||_l\int_{\wl{a}}\lo^2(y) \,|f_2(y)|\leq ||\om_2||_l||\om_1||_l.$$
So also $sup_{a\in\G^0}\, \int_{\wl{a}} \lo^2(x)|f_1*f_2|(x)\leq
||\om_2||_l||\om_1||_l$. In the same way one can prove the inequality for
$||.||_r$ or use the fact that $(\sA,*)$ is a *-algebra which is proved in
lemma \ref{star}. }\\
\dowl
\end{lem}
\begin{re} 
\newcommand{\ag}{\vphantom{x}_a{\G}}
{\em 
We can try to define a ``geometric'' norm on $\sA(\G)$ as follows.
Recall that the orbit of a point $a\in E$ is a set: 
$\el(\G_a)=\er(\vphantom{x}_a{\G})$. 
It is known \cite{McK} that for each $a \in \G^0$ the set 
$\G_a\cap\vphantom{x}_a{\G}$ is a submanifold in $\G$ and a Lie group. 
Since right and left translations are diffeomorphisms of the fibers it is clear
that all sets $\vphantom{x}_a{\G}\cap\G_b$ for $a,b$ in the same orbit are 
diffeomorphic submanifolds. Also we have that $\el|_{\G_a}: \G_a\lra \G^0$ and
$\er|_{\vphantom{x}_a{\G}}:\vphantom{x}_a{\G}\lra \G^0$ are of constant rank so
orbits are immersed submanifolds.} Now suppose that each orbit in $\G^0$ is a 
submanifold {\em let us denote the orbit through $a$ by $O_a$. In this case 
$(\G_a,O_a,\el|_{\G_a})$ and 
$(\vphantom{x}_a{\G},O_a, \er|_{\vphantom{x}_a{\G}})$ are 
locally trivial differential fibrations, with the fibers diffeomorphic to a 
Lie group $\vphantom{x}_a{\G}\cap\G_a$. If $\lambda$ is a half density on
$\G$ along the left fibers then its restriction to  $\vphantom{x}_a{\G}$ 
 can be written as $\lambda(x)=\mu(\er(x))\mt\nu(x)$ for 
$\mu$ half density on $O_a$ and $\nu$ half 
density on $\vphantom{x}_a{\G}$ along the fibers of 
$\er|_{\vphantom{x}_a{\G}}$. In the same if $\rho$ is a  half density along
the right fibers than its restriction to $\G_a$ can be written as 
$\rho(x)=\mu_1(\el(x))\mt\nu_1(x)$ for $\mu_1$ - half density on the orbit and
$\nu_1$ - half density along the fibers of $\el|_{\G_a}$.
So  for $\om=\lambda\mt\rho\in \sA(\G)$ 
and $x\in\G$ with $\el(x)=a,\er(x)=b$ we have 
$\om(x)=\mu(\er(x))\mt\mu_1(\el(x))\mt \nu(x)\mt\nu_1(x)$ but since fibers of
$\er|_{\G_a}$ and $\el|_{\vphantom{x}_b{\G}}$ are the same $\nu(x)\mt\nu_1(x)$
is a density on $T_x(\G_a\cap\vphantom{x}_b{\G})$. Let $S$ denote the set of
orbits and define the following quantity:
$$||\om||_{geom}:=\sup_{s\in S}\sqrt{\int_{s\times s}|\mu_1|^2(a)\mt |\mu|^2(b)
(\int_{\vphantom{x}_a{\G}\cap\G_b}|\nu\nu_1|)^2}.$$
Let us show that this quantity is finite and for $\om\in\sA(\G)$: 
$||\om||_{geom}\leq ||\om||$ where the norm on the right side is introduced
above.\\
Let $\om=f\lo\mt\ro\,,\, \lo(x)=:\mu(b)\mt \nu(x)\,,\,
\ro(x)=:\mu_1(a)\mt\nu_1(x)\,,\,a:=\el(x)\,,\,b=\er(x)$. Choose an orbit $O_c$.
We have to estimate integral: 
$$I(c):=\int_{O_c\times O_c}(\mu_1)^2(a)\mt(\mu)^2(b)
\left[\int_{\vphantom{x}_a{\G}\cap\G_b}|f||\nu||\nu_1|\right]^2.$$ 
First we use the Schwarz
inequality for integral over $\ag\cap\G_b$: 
$$\int_{\ag\cap\G_b}|f||\nu||\nu_1|\leq
\left(\int_{\ag\cap\G_b}|f|\nu^2\right)^{1/2}
\left(\int_{\ag\cap\G_b}|f|\nu_1^2\right)^{1/2}.$$ 
In this way we get:  
$$I(c)\leq\int_{O_c\times O_c}(\mu_1)^2(a)\mt\mu^2(b)\,g(a,b)g_1(a,b),$$ 
where: $g(a,b):=\int_{\ag\times\G_b}|f|\nu^2$ and 
$g_1(a,b):=\int_{\ag\times\G_b}|f|\nu_1^2$. Then $g,g_1$ are compactly 
supported, smooth functions on $O_c\times O_c$. 
Let $g_1(a,b_0):=\sup_{b\in O_c}\{g_1(a,b)\}$. Now: 
$$\int_{O_c\times O_c}(\mu_1)^2(a)\mt\mu^2(b) g(a,b) g_1(a,b)=
\int_{O_c}(\mu_1)^2(a)\int_{O_c}\mu^2(b) g_1(a,b) g(a,b)\leq$$
$$\int_{O_c}(\mu_1)^2(a)g_1(a,b_0)\int_{O_c}\mu^2(b)g(a,b).$$ 
But $$\int_{O_c}\mu^2(b)g(a,b)=\int_{O_c}\mu^2(b)\int_{\ag\cap\G_b}|f|\nu^2=
\int_{\ag}|f|\lo^2\leq||\om||_l.$$ 
In this way:
$$I(c)\leq ||\om||_l\int_{O_c}(\mu_1)^2(a)g_1(a,b_0)
\leq ||\om||_l\int_{O_c}(\mu_1)^2(a)\int_{\ag\cap \G_{b_0}}|f|\nu_1^2=
||\om||_l\int{\G_{b_0}}|f|\ro^2\leq||\om||_l||\om||_r.$$
So for each orbit $s\in S$ we have 
$\sqrt{I(s)}\leq\sqrt{||\om||_l||\om||_r}\leq ||\om||$ and finally:
$||\om||_{geom}\leq||\om||$.
}
\end{re}
 
Now we show that $\sA$ is a normed *-algebra and morphisms define
*-homomorphisms. As in the Prop.\ref{funktor} the equality we prove is
taken from the definition of a conjugation of a linear mapping on 
$C^*$-algebra.
\begin{prop}\label{star}
Let $h:\G\rel\G'$ be a morphism of differential groupoids. Then
$$\om_3^*(\hat{h}(\om_1)(\om_2))=(\hat{h}(\om_1^*)\om_3)^*\om_2\,\,
{\rm for\,\, any\,\,} \om_1\in \sA(\G)\,,\,\om_2,\om_3\in\sA(\G').$$
Proof: {\em Choose $\om_0$ and $\om'_0$ - *-invariant bidensities on $\G$ and
$\G'$ and write:
$\om_1=f_1\,\om_0$, $\om_2=f_2\,\om'_0$, $\om_3=f_3\,\om'_0$. Then
$\om_3^*(\hat{h}(\om_1)\om_2)(z)= (f_3^**(f_1*_hf_2))(z)\,\om'_0(z)$, where 
$$(f_3^**(f_1*_hf_2))(z)= 
\int_{F_r(z)}\rho'^2_0(y)\,f_3^*(z s(y))\,(f_1*_hf_2)(y)= 
\int_{\G_a'} \rho'^2_0(y)\, f_3^*(s(y))\,(f_1*_hf_2)(y z)\,,\,a:=\el'(z)$$ 
in this equality we have used diffeomorphism $F_r(z)\simeq\G_a'$. Using the
definition of $\hat{h}$ we can write this expression as: 
$$\int_{\G_a'}\rho'^2_0(y)\,\overline{f_3(y)} 
\int_{F_l(f_h(b))}\lo^2(x)\,f_1(x)\,t_h(x,y_1)\,f_2(y_1),$$ 
where:
$b:=\el'(y)\,,\,y_1:=s'(h_b^L(x)) y z$.\\
This integral can be interpreted as integral over 
$\G*_h\G_a'$ of the density:
$$\Psi_z(x,y):= (\lo^2(x)\mt\rho'^2_0(y))
\,(\overline{f_3(y)}\, f_1(x)\,t_h(x,y_1)\,f_2(y_1)).$$
The right hand side:
$((\hat{h}(\om_1^*)\om_3)^*\om_2)(z)=((f_1^**_hf_3)^**f_2)(z)\,\om'_0(z)$.

$$((f_1^**_hf_3)^**f_2)(z)=
\int_{F_l(a)}\lo'^2(y_2)\,(f_1^**_hf_3)^*(y_2)\,f_2(s(y_2)z)=
\int_{F_l(a)}\lo'^2(y_2)\,\overline{f_1^**_hf_3(s(y_2))}\,f_2(s(y_2) z)=$$
$$=\int_{F_l(a)}\lo'^2(y_2) \int_{F_l(f_h(b_2))}\lo^2(x_2)
\,f_1(s(x_2))\,t_h(x_2,y_3)\,\overline{f_3(y_3)}\, f_2(s(y_2) z),$$ 
where: $b_2:=\er'(y_2)\,,\,y_3:=s'(h_{b_2}^L(x_2)) s'(y_2)$.\\
Now we interpret the integral as a integral over the submanifold
$$\G\vphantom{X}_{h}{\times} F_l(a):=
\{(x_2,y_2)\in\G\times\G'\,:\,\el'(y_2)=a\,,\, \el(x_2)=f_h(\er'(y_2))\}$$ 
of the density: 
$$\Phi_z(x_2,y_2):=(\lo^2(x_2)\mt\lambda'^2_0(y_2))\, 
f_1(s(x_2))\, t_h(x_2,y_3)\, \overline{f_3(y_3)}\, f_2(s(y_2) z).$$
The mapping $s\times s'$ is a diffeomorphism from $\G\vphantom{X}_{h}{\times}
F_l(a)$ onto $\G\timh \G_a'$, moreover 
$$(s\times s')(\lo^2(x_2)\mt
\lambda'^2_0(y_2))=\ro^2(s(x_2))\mt \rho'^2_0(s'(y_2)).$$ 
From this it follows
that the above integral is eqal to the integral over $\G\timh\G_a'$ of a
density: 
$$\tilde{\Phi}_z(x_3,y_4):=
(\ro^2(x_3)\mt\rho'^2_0(y_4))\, f_1(x_3)\,t_h(s(x_3),y_3)\,
\overline{f_3(y_3)}\,f_2(y_4 z),$$ where $y_3:=h_{b_2}^R(x_3) 
y_4\,,\, b_2:=\el'(y_4)$.\\
Now we use $\got_h:\G\timh\G_a'\lra \G *_h\G_a'$ and get:
$$(\got_h \tilde{\Phi}_z)(x,y):=(\lo^2(x)\mt\rho'^2_0(y))\, f_1(x)
\,\overline{f_3(y)}\,f_2(y_5 z)\, t_h^2(x,y_5)\, t_h(s(x),y),$$ 
where $(x,y)\in \G *_h\G_a'$, $y_5:=s'(h_{b}^L(x)) y\,,\, b:=\el'(y)$.\\ 
This is equal to $\Psi_z(x,y)$ provided that 
$t_h(x,y_5z)=t_h^2(x,y_5)t_h(s(x),y) $. From the lemma \ref{app}: 
$t_h(x,y_5)=t_h(x,y_5 z)$, so
it remains to show equality: $1=t_h(x,y_5) t_h(s(x),y)$. This is done in the
following: 
\begin{lem} \label{ths}Let $(x,y)\in \G*_h\G_a'$ and $b:=\el'(y)$. Then $t_h(s(x),y)
t_h(x,s'(h_b^L(x) y)=1.$\vs\\
Proof: {\em It is easy to check that $\got_h (s\times id) \got_h=(s\times
id)$ and 
$(s\times id) i_1(\ro(x)\mt \ro'(y))=i_2(\lo(s(x))\mt\ro'(y)).$\\
We compute:
$$\got_h(s\times id)\got_h i_1(\ro(x)\mt\ro'(y))=t_h(x,y)\got_h(s\times
id)i_2(\lo(x)\mt\ro'(y'))=$$
$$=t_h(x,y) \got_h i_1(\ro(s(x))\mt\ro'(y'))=t_h(x,y)
t_h(s(x),y') i_2 (\lo(s(x))\mt\ro'(y)).$$
So $t_h(x,y) t_h(s(x),y')=1$.} 
\end{lem}}
\dow
\end{prop}

\begin{col} *-algebra structure on $\sA$. {\em Take $h=id$ then we have
$\om_3^*(\om_1 \om_2)=(\om_3^*\om_1) \om_2=(\om_1^* 
\om_3 )^*\om_2$ for any bidensities $\om_1,\om_2,\om_3$. So $(\om_3^*\om_1)=
(\om_1^* \om_3)^*$ (multiplication is non degenerate), since * is an 
involution we have that $(\om_1 \om_2)^*=\om_2^* \om_1^*$. 
This together with lemma \ref{norm} shows that
$(\sA,||.||)$ is a normed *-algebra.}
\end{col}
{\bf Representation of *-algebra of groupoid associated with a morphism.}
\vs

Now we are going to show that any morphism $h:\G\rel\G'$ of differential
groupoids defines representation of *-algebra $\sA(\G)$ in the Hilbert space
$L^2(\G')$ of square integrable half densities on $\G'$. Again we use  lemma
\ref{rfakty}.
Let $\Psi$ be a smooth half density on $\G'$ with compact support and
$\om\in\sA(\G)$, $\om=\lambda\mt\rho$. 
Let $(x,y)\in \G\timh\G'$ and $\got_h(x,y)=:(x,z)$.
As in the definition of $\hat{h}$,  $\rho(x)\mt\Psi(y)$ can be viewed as a 
half density on
$T_{(x,y)}(\G\timh\G')$ and $\got_h(\rho(x)\mt\Psi(y))$ is a half density on
$T_{(x,z)}(\G*\G')$. Since
$\oml(x)\mt\Om^{1/2}T_z\G'\simeq\Om^{1/2}T_{(x,z)}(\G*_h\G')$ 
this half density can be written as $\tilde{\lambda_x}\mt\Psi_x(z)$ for
$\tilde{\lambda_x}$--a half density on $T_x(F_l(x))$ and
$\Psi_x(z)$--a half density on $T_z(\G')$. Then $\lambda(x)
\tilde{\lambda_x}\mt\Psi_x(z)$ is a 1-density on $T_x^r\G$
with values in half densities on $T_z(\G')$. Integrating it we get half
density on $T_z(\G')$.  Let us define:
$$(\pi_h(\om)\,\Psi)(z):=\int_{F_l(f(a))}
[\lambda(x)\tilde{\lambda}(x)]\mt\Psi_x(z).$$
Choose $\om_0$ and write $\om=f\,\om_0$. Since $\er'$ is a surjective
submersion we have: $\Om^{1/2}T_w\G'\simeq\omr(w)\mt\Om^{1/2}T_{\er'(w)}E'$
for any $w\in\G'$. In this way, if we choose $\ro'$ and
$\nu_0$--non vanishing, real half density on $E'$ then $\ro'\mt\nu_0$ defines
non -vanishing, real, half density on $\G'$. So any other smooth half density
with compact support $\Psi$ can be written as $\Psi=\psi
\,\ro'\mt\nu_0=:\psi\, \Psi_0$ for
some smooth, complex function $\psi$ with compact support. It is
easy to see that: $\got_h(\ro(x)\mt\ro'(y)\mt\nu_0(a))=t_h(x,y)
\lo(x)\mt\ro'(z)\mt\nu_0(a)$ where $t_h$ is as in the definition of $\hat{h}$.
So the  explicit formula is:
$$(\pi_h(\om)\Psi)(z)=
\left[\int_{\wl{b}}\lo^2(x)\,f(x)\,t_h(x,y)\,\psi(y)\right]\Psi_0(z),$$ 
where $b:=f_h(\el'(z))\,,\,\got_h(x,y)=(x,z)$. 
Note that formally the expression is
the same as in Prop. \ref{funktor}.

\begin{prop}\label{rep}
a) Let $h:\,G\rel\G'$ be a
morphism of 
differential groupoids. Let  $||.||$ be a norm on $\sA(\G)$ associated with 
choosen $\om_0$. 
The correspondance: $\sA(\G)\ni \om \mapsto 
\pi_h(\om)$ is representation of the normed *-algebra $\sA(\G)$ in
$L^2(\G')$.\\
b) If $k:\,\G'\rel\G''$ is a morphism then:
$\pi_k(\hat{h}(\om_1)\om_2)\Psi=\pi_{kh}(\om_1)\pi_k(\om_2)\Psi$ for any
$\om_1\in\sA(\G)\,,\,\om_2\in\sA(\G')\,,\,\Psi$-smooth half density on $\G''$
with compact support.\vs\\
Proof:\,
{\em  Let $\Psi=\psi\,\Psi_0$ for $\Psi_0:=\ro'\mt\nu_0$ and let
$\om=f\,\om_0$.\\
b) This follows directely from Prop. \ref{funktor}.\vs\\
a) Let $(\;,\;)$ be a scalar product in $L^2(\G')$.  
$$|(\Psi,\pi_h(\om)\Psi)|=
\left|\int_{\G'}\overline{\psi(z)}\Psi_0^2(z)\int_{\wl{b}}\lo^2(x)
f(x)\, t_h(x,y)\, \psi(y)\right|,$$ 
where $b:=f_h(\el'(z))$ and $\got_h(x,y)=(x,z)$.\\
Using the definition of $\Psi_0$ one can estimate:
$$\left|\int_{E'}\nu_0^2(a)\int_{\G_a'}\ro'^2(z)\,\overline{\psi(z)}
\int_{\wl{b}}\lo^2(x)\, f(x)\, t_h(x,y)\, \psi(y)\right|\leq 
\int_{E'}\nu_0^2(a)\int_{\G_a'}\ro'^2(z)\,|\psi(z)|
\int_{\wl{b}}\lo^2(x)\, |f(x)\, t_h(x,y)\, \psi(y)|.$$
For fixed $a\in E'$ the integral 
$$\int_{\G_a'}\ro'^2(z)\,|\psi(z)| \int_{\wl{b}}\lo^2(x)\, 
|f(x)|\, |t_h(x,y)|\, |\psi(y)|$$
can be viewed as the
integral over $\G*_h\G_a'$ of a 1-density 
$$(\lo^2(x)\mt\ro'^2(z))\, |\psi(z)\, f(x)\, t_h(x,y)\, \psi(y)|$$ 
and it is equal to 
the scalar product in $L^2(\G *_h\G_a')$ of half
densities $(\varphi_1,\varphi_2)$ for 
$$\varphi_1(x,z):= \sqrt{|f(x)|}\,
|\psi(z)|\,sgn(t_h(x,y))\,\lambda_0(x)\mt\ro'(z)\,,$$
$$\varphi_2(x,z):=\sqrt{|f(x)|}\, t(x,y)\,|\psi(y)|,(\lambda_0(x)\mt\ro'(z)).$$
(since $t_h(x,y)\neq 0\,$, $sgn(t_h)$ is well defined, smooth function). But
since $\got_h: \G\timh\G_a'\lra \G *_h \G_a'$ is a diffeomorphism, it defines
unitary operator $\got_h: L^2(\G\timh\G_a')\lra L^2(\G * \G_a')$ and
$\varphi_2=\got_h \tilde{\varphi_2}$ where $\tilde{\varphi_2}(x,y):=
\sqrt{|f(x)|}\,|\psi(y)|\,(\rho_0(x)\mt\ro'(y))$--is a half density on
$\G'\timh\G_a'$. So we have:
$|(\varphi_1,\varphi_2)|\leq
||\varphi_1||\,||\got_h\tilde{\varphi_2}||=
||\varphi_1||\,||\tilde{\varphi_2}||$.\\
$$|| \varphi_1 ||^2=\int_{\G *\G_a'}|\varphi_1|^2=
\int_{\G_a'}\ro'^2(z)\,|\psi(z)|^2\,\int_{\wl{b}}\lo^2(x)\,|f(x)| \leq
||\om||_L\,\int_{\G_a'}\ro'^2(z)\,|\psi(z)|^2.$$
$$||\tilde{\varphi_2}||^2=\int_{\G\timh \G_a'}|\tilde{\varphi_2}|^2=
\int_{\G_a'}\ro'^2(y)\,|\psi(y)|^2 \int_{\G_c}\ro^2(x)\,|f(x)| \leq
||\om||_R \int_{\G_a'}\ro'^2(y)\,|\psi(y)|^2,$$ where $c:=f_h(\el'(y)).$\\
So finally we get an estimate: 
$$|(\Psi,\pi_h(\om)\Psi)|\leq
\int_{E'}\nu_0^2(a)\,\sqrt{||\om||_L||\om||_R}\int_{\G_a'} \ro'^2(z)
\,|\psi(z)|^2=\sqrt{||\om||_L||\om||_R}\,||\Psi||^2\leq
||\om||\,||\Psi||^2.$$ 
This shows that the operator $\pi_h(\om)$ ( defined on smooth half densities
with compact support) is bounded. Since smooth half
densities with compact support are dense in $L^2(\G')$, $\pi_h(\om)$ can be
uniquely extended to a bounded operator on $L^2(\G')$.\\
Putting $h=id:\,\G\rel\G$ we get from b): 
$\pi_k(\om_1\,\om_2)\Psi=\pi_k(\om_1)\pi_k(\om_2)\Psi$. This shows that
$\pi_k$ is a representation.\\
Now we show that $\pi_h(\om^*)=(\pi_h(\om))^*.$
$$(\pi_h(\om^*)\Psi,\Psi)=
\int_{E'}\nu_0^2(a)\int_{\G_a'}\ro'^2(z)\,\psi(z)
\overline{\int_{\wl{c}}\lo^2(x) \,f^*(x)\,t_h(x,y)\,\psi(y)} =$$
$$=\int_{E'}\nu_0^2(a)\int_{\G_a'}\ro'^2(z)\,\psi(z)
\int_{\wl{c}}\lo^2(x) \,f(s(x))\,t_h(x,y)\,\overline{\psi(y)},$$ 
where $c:=f_h(\el'(z))\,,\,\got_h(x,y)=(x,z)$. \\
In the same way: 
$$(\Psi,\pi_h(\om)\Psi)=\int_{E'}\nu_0^2(a)\int_{\G_a'}\ro'^2(z)\,
\overline{\psi(z)} \int_{\wl{c}}\lo^2(x) \,f(x)\,t_h(x,y)\,\psi(y).$$
We will show that: 
$$\int_{\G_a'}\ro'^2(z)\,\psi(z)
\int_{\wl{c}}\lo^2(x) \,f(s(x))\,t_h(x,y)\,\overline{\psi(y)}=
\int_{\G_a'}\ro'^2(z)\,\overline{\psi(z)} \int_{\wl{c}}\lo^2(x) 
\,f(x)\,t_h(x,y)\,\psi(y).$$
In the Prop. \ref{star} the following equality was proved:
$$\int_{\G_a'}\ro'^2(y)
\overline{f_3(y)}\int_{\wl{c}}\lo^2(x)\,f_1(x)\,t_h(x,y_1)\, f_2(y_1)=
\int_{\vphantom{X}_{a}{\G'}}\lo'^2(y_2)
\int_{\wl{c_2}}\lo^2(x_2)\,f_1(s(x_2))\,t_h(x_2,y_3)\,\overline{f_3(y_3)}\, 
f_2(s(y_2)),$$
where
$c:=f_h(\el'(y))\,,\,\got_h(x,y_1)=(x,y)\,,\,
c_2:=f_h(\er'(y_2))\,,\,\got_h(x_2,y_3)=(x_2,s(y_2))\,$, 
$f_2,f_3$ are smooth function on $\G'$ with compact support and $f_1$ is
smooth function on $\G$ with compact support.  \\
Using $s'$ we can rewrite the right hand side of the equality as:
$$\int_{\G_a'}\ro'^2(y_4)\int_{\wl{c_4}}\lo^2(x_2)\,f_1(s(x_2))\,
t_h(x_2,y_3)\,\overline{f_3(y_3)}\,f_2(y_4),$$
where $c_4:=f_h(\el'(y_4))\,,\,\got_h(x_2,y_3)=(x_2,y_4)\,$. \\
Now put $f_2=f_3=\psi\,,\,f_1=f$. We get:
$$\int_{\G_a'}\ro'^2(y)\overline{\psi(y)}\int_{\wl{c}}\lo^2(x)\,f(x)\,
t_h(x,y_1)\, \psi(y_1)=$$
$$=\int_{\G_a'}\ro'^2(y)\,\psi(y)\int_{\wl{c}}\lo^2(x)\,f(s(x))\,
t_h(x,y_1)\,\overline{\psi(y_1)}.$$
And this is desired equality. }\\
\dow
\end{prop}
\begin{ex}\label{red-mod}
a)  Reduced $C^*$-algebra of a differential groupoid. 
{\em Let $l$ be the morphism from $\G$ to the pair groupoid $\G\times \G$
defined in Example \ref{exmor} f) i.e. $(x,y;z)\in Gr(l)\iff (x;z,y)\in
Gr(m)$. It is easy to see that in this case: $f_l=\el$,
$\G\times_l(\G\times\G)=\{(x,y,z)\in\G\times\G\times\G\,:\,\er(x)=\el(y)\}\,$,
$\G*_l(\G\times\G)=\{(x,y,z)\in\G\times\G\times\G\,:\,\el(x)=\el(y)\}$ and
$\got_l(x;y,z)=(x; x y, z)$. $\pi_l$ is a representation of $\sA(\G)$ in
$L^2(\G\times\G)=L^2(\G)\mt L^2(\G)$ and a short computation shows that
$\pi_l=\pi_{id}\mt I$. So
$||\pi_l(\om)||=||\pi_{id}(\om)||.$ Also from Lemma \ref{nondegmul} easily
 follows that
$\sA(\G)\ni\om\mapsto ||\pi_{id}(\om)||$ is a $C^*$ norm on $\sA(\G)$. The
completion of $\sA(\G)$ in this norm will be called} reduced
$C^*$-algebra of $\G$ {\em and denoted by $C^*_{red}(\G)$.}\vs\\
b) Modular function. {\em Let $\te$ be a morphism from $\G$ to the pair
groupoid $E\times E$ defined in
Example \ref{exmor} e), i.e. $Gr(\te)=\{(\el(x),\er(x);x)\,:\,x\in \G\}$.
It is easy to see that: $\G\times_{\te}(E\times E)=\{(x;\er(x),e)\,:\,x\in
\G\,,\,e\in E\}\,,\,\G*_{\te}(E\times E)=\{(x;\el(x),e)\,:\,x\in
\G\,,\,e\in E\}\,,\,m_{\te}(x;\er(x),e)=(\el(x),e)$ and
$\got_{\te}(x;\er(x),e)=(x;\el(x),e)$. Choose some $\om_0=\lo\mt\ro$ and some
real, non vanishing half density $\nu_0$ on $E$. Such choice defines function
$t_{\te}(x;\er(x),e)$. From Lemma \ref{app}, this function does not depend
from $e$ and if we define $\D(x):=t_{\te}(x,\er(x))$ then $\D(x y)=\D(x)
\D(y)$ for any composable $x,y\in \G$. $\D$ is called } modular function of
$\G$ {\em ( it depends from chosen $\lo,\nu_0$). 
The function $\D$ can also be described in the following way. 
When $\om_0\,,\,\nu_0$ 
are choosen, the expressions: $\psi_r(x):=\ro(x)\mt\nu_0(\er(x))$ and 
$\psi_l(x):=\lo(x)\mt\nu_0(\el(x))$ define smooth, non vanishing, real 
half densities on $\G$. Then $\D$ is defined by:  $\psi_l=:\D \psi_r$.}
\end{ex}
\begin{re} {\em 
Dependance of $\D$ on choice of $\lo$ and  $\nu_0$ can be described in the 
cohomological way. 
Let us define   $\G^{(0)}:=\G^0\,,\,\G^{(1)}:=\G$ and, for  
$n\geq 2$,
 $\G^{(n)}:=\{(x_0,...,x_{n-1})\in \G\times\dots\times \G :
\er(x_i)=\el(x_{i+1})\,,\,i=0,\dots,n-1\}$. \\
(Smooth) n-cochain is a smooth function 
 $f:\G^{(n)}\lra {\R}\setminus\{0\}$ which, for $n>0$, satisfies condition:
$$[\exists\,i\in\{0,\dots,n-1\}:x_i\in\G^0]\Rightarrow 
f(x_0,\dots,x_i,\dots,x_{n-1})=1.$$ 
Group of  n-cochains  (with a pointwise multiplication)
we denote by  $C^n(\G)$. 

Define coboundary operators 
$\delta^n : C^n(\G)\lra C^{n+1}(\G)$:  
$$(\delta^0 f)(x):=\frac{f(\el(x))}{f(\er(x))}\,\,{\rm and \,,\,\, for \,\,} 
n>0,$$ 
$$(\delta^n f)(x_0,x_1,\dots,x_n):=$$
$$=f(x_1,\dots,x_n)\prod_{i=1}^n(f(x_0,\dots,x_{i-1}x_i,\dots,x_n))^{s(i)} 
(f(x_0,\dots,x_{n-1}))^{s(n+1)},$$ where  $s(i):=(-1)^i$.\\
It is easy to check that  $\delta^{n+1}\delta^n=1$. In this way we get complex
and cocycles, coboundaries and cohomology groups \cite{Ren}.

Now, let  $\tilde{\lo},\tilde{\nu_0}$ be other half densities.
We have  $\tilde{\lo}(x)=f(\er(x))\lo(x)$ and  
$\tilde{\nu_0}(a)=g(a)\nu_0(a)$ for some smooth, non vanishing, real
functions on $\G^0$. Then from the equality 
$\tilde{\psi_l}=\tilde{\D}\tilde{\psi_r}$ we get 
 $$\tilde{\D}(x)=\D(x)\frac{f(\er(x)) g(\el(x))}{f(\el(x)) g(\er(x))}=
\D(x)(\delta^0 \frac{g}{f})(x).$$ 
So  $\D$ i $\tilde{\D}$ are in the same cohomology class.}
\end{re}
{\bf Nondegeneracy of morphisms.}

Now we are going to show that the action of morphisms on bidensities is
nondegenerate (i.e. if $\hat{h}(\om)\om'=0$ for any $\om$ then $\om'=0$.)
This is important in the context of morphisms of $C^*$-algebras. In
fact we prove more general: 
\begin{prop}\label{apjed}
Let $h:\G\rel\G'$ be a morphism of differential groupoids. Then for any
$\om'\in\sA(\G')$ there exists a sequence $\om_n\in\sA(\G)$ such that:
$\lim_{n\rightarrow\infty} \hat{h}(\om_n)\om'=\om'$. (The limit is in
topology defined by some $||.||$ of above defined type on $\G'$.)\vs\\
{\em The proof is based on several lemmas and is slightly modified version
of the proof given in \cite{Ren}:
\begin{lem}\label{apjed1}
Let $h:\G\rel\G'$ be a
morphism of differential groupoids and let 
$K\subset\G'$ be a compact subset. Then  exists $U_K$ -- open
neighbourhood of $\G^0$ in $\G$ such that $\overline{m_h((U_K\times K)\cap
(\G\timh\G'))}$ is compact.\\ 
Proof: {\em 
\begin{enumerate}
\item $m_h:\G\timh\G'\lra \G'$ is a smooth mapping and
$m_h(a,z)=z $ for $z\in \G'\,,\,a:=f_h(\el'(z))$, so for any neighbourhood
$O_z\ni z$ there exist neighbourhoods $\G\supset V_a\ni a$ and $O_z'\ni z$
such that $m_h((V_a\times O_z')\cap(\G\timh\G'))\subset O_z$. Moreover we can
assume that $\er^{-1}(V_a\cap\G^0)\cap V_a=V_a\,,\,O_z'\subset
O_z\,,\,f_h(\el'(O_z'))\subset V_a$. 
\item Let $\{O_z\,,\,z\in K\}\,,\,z\in O_z\,$ be an open covering of $K$ such
that $\overline{O_z}$ is compact.
From the previous point the family $\{O_z'\,,\,z\in K\}$ is an open covering
of $K$ and $\{V_a\,,\,a=f_h(\el'(z))\,,\,z\in K\}$ an open covering of a 
compact set $H:=f_h(\el'(K))$. So one can choose finite covering
$O_{z_1}'\cup ... \cup O_{z_m}'\supset K$ with the corresponding
$V_{a_1},...,V_{a_m}$ and $O_{z_1},...,O_{z_m}$. Then $H\subset V_1\cup
...\cup V_m$, where $V_i:=V_{a_i}$. 
\item For $x\in H$ let $W_x$ be an open (in $\G^0$) neighbourhood of $x$ 
contained in $V_1\cup ...\cup V_m$. Define $U_x:=\er^{-1}(W_x)\cap
V_{i_1}\cap ...\cap V_{i_l}$ where the intersection is with these sets $V_i$
which contain $x$. The family $\{U_x\,,\,x\in H\}$ is an open (in $\G$)
covering of $H$ so we can choose $U_1\cup ...\cup U_n=:U\supset H$.
\item Now let $(x,z)\in (U\times K)\cap (\G\timh\G')$. $z$ is contained in
some $O_{z_i}'$ and $x$ in some $U_j$. Since $\er(x)=f_h(\el'(z))$,
$\er(x)\in V_i$ and from the construction also $x\in V_i$ and then
$m_h(x,z)\in O_{z_i}$. So $m_h((U\times K)\cap(\G\timh\G'))\subset
O_{z_1}\cup ...\cup O_{z_m}\subset \overline{O_{z_1}}\cup ...\cup
\overline{O_{z_m}}$ and this is compact set.
\item Finally one defines $U_K:=U\cup \bigcup_{x\in \G^0\setminus H}
\er^{-1}(O_x)$ where $O_x$ is an open (in $\G^0$) neighbourhood of $x$ and
$O_x\cap H=\emptyset$.  
\end{enumerate}}
\dowl
\end{lem}
\begin{lem}\label{apjed2} 
Let $h:\,\G\rel\G'$ be a
morphism of differential groupoids. If 
$g$ is a continous function on $\G'$ and $K\subset\G'$ is compact then for
any $\delta>0$ there exists $U^{\delta}$ -- open neighbourhood of $\G^0$ in
$\G$ such that:
$$(x,z)\in (U^{\delta}\times K)\cap (\G *_h\G') \Rightarrow
|g(s'(h_a^L(x)) z)-g(z)|\leq \delta\,,\,\,{\rm  where\,\,} a:=\el'(z).$$
Proof: {\em 
Choose $\delta>0$. Then for any $z\in K$ there exists
$O^{\delta}_z$ -- open neighborhood of $z$ such that: $y\in
O^{\delta}_z\Rightarrow |g(y)-g(z)|\leq\delta/2.$ The family
$\{O_z^{\delta}\,,\,z\in K\}$ is a covering of $K$. Put $U_1^{\delta}:=U_K$
where $U_K$ is constructed from $\{O_z^{\delta}\}$ as in the previous lemma.
Then for $(x,z)\in (U_1^{\delta}\times K)\cap (\G \timh\G')$ we have $z\in
O^{\delta}_{z_i}$ for some $z_i\in K$ and $m_h(x,z)\in O^{\delta}_{z_i}$, then
$|g(m_h(x,z))-g(z)|=|g(m_h(x,z))-g(z_i)+g(z_i)-g(z)|\leq
\delta/2+\delta/2=\delta$. Finally we define $U^{\delta}:=s(U_1^{\delta})$ and
the result follows from $(x,z)\in (\G*_h\G')\iff (s(x),z)\in \G\timh\G'$ and
$s'(h_a^L(x))=h_a^R(s(x))$.}\\
\dowl
\end{lem}
\begin{lem}\label{apjed3} 
Let $h\,:\G\rel\G'$ be a 
morphism of differential 
groupoids. Choose $\om_0,\om_0'$ -- this defines function $t_h$. Let
$K\subset \G'$ be compact 
and let $1>\delta>0$. Then there 
exists $U_1^{\delta}\subset\G$ -- open neighbourhood of $\G^0$ such that:
$$(x,z)\in (U_1^{\delta}\times K)\cap (\G*_h\G')
\,\Rightarrow\,|t_h(x,s'(h_a^L(x)) z)-1|\leq \delta\,,\,\,
{\rm where\,\,} a:=\el'(z).$$
Proof: {\em We begin by showing that $t_h(f_h(b),b)=1$ for any $b\in E'$.
Let $a:=f_h(b)$, $V:=T_{(a,b)}(\G\timh \G_b')\,,\,W:=T_{(a,b)}(\G*_h\G_b')$. 
Let $(X_1,...X_k)$ be a basis in $T_a^r\G$ and $(Y_1,...Y_m)$ basis in 
$T_b^r\G'$. Then $(s(X_1),...,s(X_k))$ is a basis in $T_a^l\G$,  
$(\hat{X_1},...,\hat{X_k},\hat{Y_1},...\hat{Y_m})$ is a basis in $V$ and 
$(\tilde{X_1},...\tilde{X_k},\hat{Y_1},...,\hat{Y_m})$ is a basis in $W$.
Where $\hat{X_i}:=(X_i,0_b)\,,\,\tilde{X_i}:=(s(X_i),0_b)\,,\,
\hat{Y_i}:=(f_h\el'(Y_i),Y_i)$. 

The isomorphism 
$i_1:\omr(a)\mt\omr(b)\lra\omh V$ is given by the formula:
$$i_1(\ro(a)\mt\ro'(b))(\hat{X_1}\dz...\dz\hat{Y_m}):=\ro(a)(X_1\dz...\dz X_k)
\ro'(b)(Y_1\dz...\dz Y_m)$$
and $i_2:\oml(a)\mt\ro(b)\lra\omh W$ by: 
$$i_2(\lo(a)\mt\ro'(b))(\tilde{X_1}\dz...\dz \hat{Y_m})
:=\lo(a)(s(X_1)\dz... \dz s(X_k))\ro'(b)(Y_1\dz...\dz Y_m).$$
Moreover $\got_h(\hat{Y_i})=\hat{Y_i}$ and 
$\got_h(\hat{X_i})=(X_i,h_b^R(X_i)b)
=:\alpha_{il}\tilde{X_l}+\beta_{ij}\hat{Y_j}.$
So 
$$i_2(\lo(a)\mt\ro'(b))(\got_h\hat{X_1}\dz...\dz \got_h\hat{Y_m})=
|\det \alpha|^{1/2}\lo(a)(s(X_1)\dz...\dz s(X_k))\ro'(b)(Y_1\dz...\dz Y_m)=$$
$$=|\det \alpha|^{1/2}\ro(a)(X_1\dz...\dz X_k)\ro'(b)(Y_1\dz...\dz Y_m)$$ 
and
$$(\got_h i_1)(\ro(a)\mt\ro'(b))(\got_h \hat{X_1}\dz...\dz \got_h\hat{Y_m})=
\ro(a)(X_1\dz...\dz X_k)\ro'(b)(Y_1\dz..\dz Y_m).$$
In this way $t_h(a,b)=|\det \alpha |^{1/2}.$
From the definition of $\alpha$: 
$$(X_i,h_b^R(X_i)b)=\alpha_{il}\tilde{X_l}+
\beta_{ij}\hat{Y_j}=\alpha_{il}(s(X_l),0_b)+\beta_{ij}(f_h\el'(Y_j),Y_j).$$
So $X_i=\alpha_{il}s(X_l)+\beta_{ij}f_h\el'(Y_j)$ and this is decomposition
of $X_j$ with respect to the direct sum $T_a\G=T_a^l\G\oplus T_aE$. Applying
$s$ to this decomposition we easy get $\alpha^2=1$, so $|\det \alpha|=1$.

Let $1>\delta>0$ be given and let $H\in E'$ be compact. Arguing as in the proof
of Lemma \ref{apjed1} we can find $H\subset O_H$ -- open in $\G'$ and $U_H$ --
an open neighbourhood of $\G^0$ in $\G$ such that:
$$(x,y)\in (U_H\times O_H)\cap (\G\timh\G')\Rightarrow |t_h(x,y)-1|\leq 
\delta/2.$$ 
From Lemma \ref{app} $t_h(x,yz)=t_h(x,y)$ so above estimate is valid for 
$y\in\G'$ with $\el'(y)\in H$. So putting $H:=\el'(K)$ we get:
$$(x,y)\in (U_H\times K)\cap (\G\timh\G')\Rightarrow |t_h(x,y)-1|\leq
\delta/2.$$
Define $U_1^{\delta}:=s(U_H)$. 
Then: 
$$(x,y)\in (U_1^{\delta}\times K)\cap(\G*_h\G')\Rightarrow (s(x),y)\in
(U_H\times K)\cap(\G\timh\G')\Rightarrow |t_h(s(x),y)-1|\leq\delta/2.$$
But from Lemma \ref{ths} $t_h(s(x),y)=\frac{1}{t_h(x,s'(h_a^L(x))y)}\,,\,
a:=\el'(y)$, so finally we have:
$$|t_h(x,s'(h_a^L(x))y)-1|=\frac{|t_h(s(x),y)-1|}{|t_h(s(x),y)|}\leq\delta.$$}
\dowl
\end{lem}

\begin{lem} For any compact $L\subset\G^0$, there exists 
sequence $U_{n+1}\subset U_n$ of open 
subsets of $\G$ with the following properties:\\
a) $\forall\,n\in N\,,\,L\subset U_n $\,, b) $\overline{U_1}$ is compact\,,\,
c) $s(U_n)=U_n$\,,\, d) $\bigcap_{n\in N}U_n\subset\G^0$\,,\, e) for any open
$V\supset\G^0\,$ there exists $N_V$ such that $\,U_n\subset V$ for all
$n>N_V$.\vs\\ 
Proof: {\em We begin with the following observation: for any  $a\in \G^0\,$ 
there exist neighbourhoods $U_a\subset\overline{U_a}\subset U_a'$ and
mappings $\phi_a,\phi_a'$ with properties:\\
a) $(U_a,\phi_a),\,(U_a',\phi_a')$
- are maps submitted to submersion $\el$ and $\phi_a=\phi_a'|_{U_a}$, \\
b) $\overline{U_a'}$ is compact,\\
c) $s(U_a)=U_a$,\\
d) $\phi_a : U_a \lra I^m\times I^n$ - where  $I^k:=]-1,1[^k$, \\
e) $\phi(U_a\cap\G^0)\subset\{0\}\times I^n$. \\
From the open covering of $L$ $\{U_a\,,\,a \in L\}$ ($U_a$ - as
above) we choose finite: $U_1:=U_{a_1}\cup...\cup U_{a_m}$,
then we define for $k\in N$: 
$\tilde{U}^k_a:=\phi_a^{-1}(\phi_a(U_a)\cap (]\frac{-1}{k},\frac{1}{k}[^m
\times I^n))\,,\,U^k_a:=s(\tilde{U}^k_a)\cap \tilde{U}^k_a$ and
$U_n:=U^n_{a_1}\cup...\cup U^n_{a_m}$. Then the family 
$\{U_n\}_{n\in N}$ has the desired properties.}\\
\dowl
\end{lem}

\noindent
Let $L\subset\G^0$ be compact. 
Let $\{h_n\}_{n\in N}$ be a sequence of smooth functions on $\G$ satisfying
conditions:\\
a) $0\leq h_n\leq 1\,$,\,\,\, b) supp$\,h_n\subset U^n\,$,\,\,\, 
c) $h_n\equiv1$ on $L\,$, \,\,\,d) $h_{n+1}\leq h_n$.\\
Then the functions $\G^0\ni a \mapsto \int_{F_l(a)} h_n(x) \lambda_0^2(x)$ 
are smooth, nonnegative and separated from 0 on $L$. So one can find $g_n$ --
smooth, nonnegative functions on $\G^0$ such that: 
$f_n(x):=h_n(x) g_n(\el(x))$ is smooth, nonnegative function with compact
support contained in  $U^n$ and $\int_{F_l(a)} f_n(x) \lambda_0^2(x) =1$ for
$a\in L$. Let us define $\om_n(x):=f_n(x)\, \om_0(x).$ }\vs\\
Proof of the proposition: {\em 
Let $\om'\in \sA(\G')\,\,\,\om'=f'\,\om_0'$ has the support contained in
$K$ and let $U_K$ be as in the lemma \ref{apjed1}. Then for any
$\om\in\sA(\G)$ with support in $U_K\,,\, \hat{h}(\om)(\om')$ has support
contained in $H$ -- the fixed compact subset of $\G'$. Take $\delta>0$ then
from lemmas \ref{apjed2} and \ref{apjed3} we get $U^{\delta}$ and
$U_1^{\delta}$ - open neighbourhoods of $\G^0$ in $\G$ 
such that for any $(x,z)\in ((U^{\delta}\cap U_1^{\delta})\times H)\cap (\G
*_h\G')$ we have: $|f'(s'(h_a^L(x)) z)-f'(z)|\leq \delta$ and
$t_h(x,s'(h_a^L(x)) z)-1|\leq\delta $. Let $U:=U_K\cap U^{\delta}\cap
U_1^{\delta}$ and $L:=f_h(\el'(H))$.  Let $\om_n=f_n\,\om_0$ be as above. 
Then for $n>N_0$ support of $\om_n$ is contained in in $U$.
Then the support of $\hat{h}(\om_n)\om'$ 
is contained in $H$. \\
For $z\in H$ we have: $\hat{h}(\om_n)\om'=(f_n*_h
f')\,\om_0'$ and $f_n*_h f'(z):=\int_{F_l(a)} \lo^2(x) f_n(x)
t_h(x,y z)f'(y z)\,\, a:=f_h(\el'(z))$ and we put for simplicity of the
notation $y:=s'(h_a^L(x))$.
$$|f_n*_h f'(z)-f'(z)|=\left|\int_{F_l(a)} \lo^2(x) f_n(x)
t_h(x,y z)f'(y z)-f'(z)\right|=$$
$$=\left|\int_{F_l(a)} \lo^2(x) [f_n(x)(t_h(x,y z)f'(y z)-f'(z))]\right|
\leq \int_{F_l(a)} \lo^2(x) |f_n(x)||t_h(x,y z)f'(y z)-f'(z)|=$$
$$=\int_{F_l(a)} \lo^2(x) |f_n(x)||f'( y z)-f'(z)|+
\int_{F_l(a)}\lo^2(x) |f_n(x)| |t_h(x,y z)-1||f'(y z)|]\leq $$
$$\leq\delta \left(\int_{F_l(a)} \lo^2(x) |f_n(x)|+\int_{F_l(a)} \lo^2(x)
|f_n(x)||f'( y z)|\right) \leq \delta( 1+M)\,,\,\,
{\rm  for\,\,} M:=\sup\,|f'|.$$
In this way $\sup |f_n*_hf'-f'|\leq \delta(1+M)$ for $n>N_0$. \\
Let $g$ be a smooth function on $\G'$ with compact support $0\leq g(x)\leq 1$
and $g|_H=1$. \\
Then 
$$||\hat{h}(\om_n)\om'-\om'||_L=
\sup_{b\in\G'^0}\int_{F_l(b)}\lo'^2(z)|f_n*_h f'(z)-f'(z)|\leq$$
$$\leq \sup_{b\in\G'^0}\int_{F_l(b)}\lo'^2(z)g(z)|f_n*_h f'(z)-f'(z)|\leq$$
$$\delta(1+M) \sup_{b\in\G'^0}\int_{F_l(b)}\lo'^2(z)\,g(z)
\leq\delta(1+M)M_g\,,\,\,{\rm  where\,\,}
M_g:=\sup_{b\in\G'^0}\int_{F_l(b)}\lo'^2(z)\,g(z).$$ 
In the same way we have $||\hat{h}(\om_n)\om'-\om'||_R\leq \delta (1+M) m_g$,
where $m_g:=\sup_{b\in\G'^0}\int_{F_r(b)}\ro'^2(z)\,g(z)$.
This proves that $\lim_{n\rightarrow\infty}\hat{h}(\om_n)\om'=\om'$.
} \\
\dow
\end{prop}

\begin{col}\label{nondeg1}
The action of morphism
on bidensities is nondegenerate. 
{\em Indeed, if $\hat{h}(\om)\om'=0$ for any $\om$ then taking $\om_n$ as
above we have: $0=\lim_{n\rightarrow \infty} \hat{h}(\om_n)\om'=\om'$.
}\end{col}
\begin{col}\label{nondeg2}
The representation
$\pi_h$ assosiated with morphism $h$ is nondegenerate.\vs\\
Proof: {\em Choose $\Psi_0=\ro'\mt \nu_0$ and let $\Psi:=\psi\,\Psi_0$ be a
smooth half density on $\G'$ with compact support. Let $\om_n$ be as above.
Then sup\,$|f_n*_h\psi-\psi|\leq \delta(1+M)$ for $n>N_0$ and $M:=sup|\psi|.$
Let $g$ be a smooth function on $\G'$ with compact support $0\leq g(x)\leq 1$
and $g|_H=1$. 
$$||\pi_h(\om_n)\Psi-\Psi||^2=\int_{E'}\nu_0^2(a)\int_{\G_a'}\ro'^2(x)
|f_n*_h\psi-\psi|^2\leq \int_{E'}\nu_0^2(a)\int_{\G_a'}\ro'^2(x)
|f_n*_h\psi-\psi|^2|g(x)|^2\leq $$
$$ \leq
\delta^2(1+M)^2\int_{E'}\nu_0^2(a)\int_{\G_a'}\ro'^2(x)|g(x)|^2=
\delta^2(1+m)^2||\Psi_g||^2,$$
where $\Psi_g:=g\,\Psi_0$ and $n>N_0$.\\
So for any $\Psi$-- smooth half density on $\G'$ with compact support there
exists a sequence $\{\om_n\}$\,,\,$\om_n\in \sA(\G)$ such that 
$\lim_{n\rightarrow\infty}\pi_h(\om_n)\Psi=\Psi$.\\
Let $\Phi\in L^2(\G')$ be such that $\pi_h(\om)\Phi=0$ for any
$\om\in\sA(\G)$. Then for any $\Psi\in L^2(\G')$--smooth with compact support
we have $(\Psi,\Phi)=\lim (\pi_h(\om_n)\Psi,\Phi)=\lim
(\Phi,\pi_h(\om_n^*)\Psi)=0$. So $\Phi=0$ and this is nondegeneracy
condition. }\end{col}
\begin{re} Note that in the statements above, we don't claim that 
$\hat{h}\om_n$ is an approximate identity for $(\sA(\G'),||.||)$, since $\om_n$
depend from chosen $\om'$. Also $\pi_h(\om_n)$ does not converge strongly or 
weakly to identity on $L^2(\G')$.
\end{re}

\noindent
{\bf Proof of the Lemma.\ref{app}}

\noindent
{\bf A}. The set $\G\timh\G'\times_k\G_a'':=\{(x,y,z)\in
\G\times\G'\times\G_a''\,:\, \er(x)=f_h(\el(y))\,,\,\er(y)=f_k(\el(z))\}$ is 
a submanifold in  $(\G\timh\G')\times \G_a''$ and in 
$\G\times(\G'\times_k\G_a'')$.\\
Indeed $\G\timh\G'\times_k\G_a''=(\er'\pi_2\times f_k \el'')^{-1}
(diag(E'\times E'))$ where $\pi_2:\,\G\timh\G'\lra\G'$ is a projection onto 
the second factor. It is easy to see that 
$(\er'\pi_2\times f_k \el'')\tran diag(E'\times E')$. 
We also have:  $\G\timh\G'\times_k\G_a''=(\er\times f_h \el'\pi_1)^{-1}
(diag(E\times E))$ where $\pi_1:\,\G'\times_k\G_a''\lra\G'$ is a projection 
onto the first factor and $(\er\times f_h \el'\pi_1)\tran diag(E\times E)$.\\
The mappings: $\G\timh\G'\times_k\G_a''\ni (x,y,z)\mapsto z\in\G_a''$ and
 $\G\timh\G'\times_k\G_a''\ni(x,y,z)\mapsto (y,z)\in \G'\times_k\G_a''$ 
are surjective submersions --- 
this is due to the fact that for any morphism $h:\,\G\rel\G'$ the mappings 
$\G\timh\G'\ni(x,y)\mapsto y\in\G'$ and $\G\timh\G_a'\ni(x,y)\mapsto y\in\G_a'$
are surjective submersions.\\
From this it follows that we have the following isomorphisms:\\
$i_1:\,\omr(x)\mt\omr(y)\lra \omh T_{(x,y)}(\G\timh\G_b')$ where $b:=\er'(y)$
.\\
$i_2:\,\omh T_{(x,y)}(\G\timh\G_b')\mt\omr(z)\lra\omh T_{(x,y,z)}
(\G\timh\G'\times_k\G_a'')$.\\
$i_3:\,\omr(y)\mt\omr(z)\lra\omh T_{(y,z)}(\G'\times_k\G_a'')$.\\
$i_4:\,\omr(x)\mt\omh T_{(y,z)}(\G'\times_k\G_a'')\lra\omh T_{(x,y,z)}
(\G\timh\G'\times_k\G_a'')$.
\begin{lem}
$i_2(i_1\mt id)=i_4(id\mt i_3)$\vs\\
Proof: {\em Let $V:=T_{(x,y,z)}(\G\timh\G'\times_k\G_a'')$. Let also choose
$B_x,B_y$ such that $T_x\G=T_x^r\G\oplus B_x$ and $T_y\G'=T_y^r\G'\oplus
B_y$. Define: $V_1:=\{(\dot{x},\dot{y},\dot{z})\in V\,:\,\dot{x}\in
B_x,\dot{y}\in B_y,\dot{z}\in T_z^R\G''\}$ and $W_1:=\{(\dot{x},\dot{y},0)\in
V\,:\,\dot{y}\in T_y^R\G'\}$. It is clear that $V=V_1\oplus W_1$, $W_1$ is a
kernel of the projection onto $\G_a''$ and $W_1$ is isomorphic to
$T_{(x,y)}(\G\timh\G_b')$ -- the isomorphism is given by (tangent to) 
projection $\pi_{12}:(x,y,z)\mapsto (x,y)$. The isomorphism $i_2$ is now 
given by:
$i_2(\om_{xy}\otimes \rho_z)(v_1\dz w_1):=
\om_{xy}(\pi_{12}w_1)\rho_z(\pi_3 v_1),$
where $\om_{xy}\in\omh T_{(x,y)}(\G\timh\G_b')\,,\,\rho_z\in\omh
T_z^r\G''\,,\, v_1\in\Lambda^{max}V_1\,,\,w_1\in \Lambda^{max}W_1$ and
$\pi_3$ is the projection onto the third factor.  \\
Let $V_2:=\{(\dot{x},\dot{y},0)\in W_1\,:\,\dot{x}\in B_x\}$ and
$V_3:=\{(\dot{x},0,0)\in W_1\,:\,\dot{x}\in T_x^r\G\}$. We have
$W_1=V_2\oplus V_3$ and the isomorphism $i_1$ is given by
$i_1(\rho_x\otimes\rho_y)(v_2\dz v_3):=\rho_x(\pi_1 v_3)\rho_y(\pi_2
v_2)$ with the obvious notation. So we have $i_2(i_1\otimes
id)(\rho_x\otimes\rho_y\otimes\rho_z)(v_1\dz v_2\dz v_3)=\rho_x(\pi_1
v_3)\rho_y(\pi_2 v_2)\rho_z(\pi_3 v_1)$.\\
It is clear that $V_3$ is the kernel of the projection
$\pi_{23}:\,(x,y,z)\mapsto(y,z)$. So $i_4$ is given by:
$i_4(\rho_x\otimes\om_{yz})(v_1\dz v_2\dz v_3):=\rho_x(\pi_1
v_3)\om_{yz}(\pi_{23} v_2\dz \pi_{23} v_1)$. Also
$i_3(\rho_y\otimes\rho_z)(\pi_{23} v_2\dz\pi_{23}
v_1):=\rho_y(\pi_2\pi_{23} v_2\dz \pi_3\pi_{23} v_1)=\rho_y(\pi_2
v_2)\rho_z(\pi_3 v_1)$. And this is desired equality.}\\
\dowl
\end{lem}

\noindent
{\bf B.} The set $\G*_h\G'\times_k\G_a'':=\{(x,y',z)\in
\G\times\G'\times\G_a''\,:\, \el(x)=f_h(\el(y'))\,,\,\er(y')=f_k(\el(z))\}$ 
 is a submanifold in $\G\times(\G'\times_k\G_a'')$ and in $(\G*_h\G')\times\G_a''$.\\
The argument is as above. Write $\G*_h\G'\times_k\G_a''=
(\el\times f_h \el'\pi_1)^{-1}(diag(E\times E))=
(\er'\pi_2\times f_k\el'')^{-1}(diag(E'\times E'))$ and use the transversality.
\\
The mappings: $\G*_h\G'\times_k\G_a''\ni(x,y',z)\mapsto z\in\G_a''$ and 
 $\G*_h\G'\times_k\G_a''\ni(x,y',z)\mapsto (y',z)\in\G'\times_k\G_a''$ are 
surjective submersions. \\
As above this provides us with the following isomorphisms:\\
$i_5:\,\oml(x)\mt\omr(y')\lra \omh T_{(x,y')}(\G*_h\G_b')$.\\
$i_6:\, \omh T_{(x,y')}(\G*_h\G_b')\mt\omr(z)\lra \omh T_{(x,y',z)}
(\G*_h\G'\times_k\G_a'')$.\\
$i_7:\,\omr(y')\mt\omr(z)\lra \omh T_{(y',z)}(\G'\times_k\G_a'')$.\\
$i_8:\,\oml(x)\mt \omh T_{(y',z)}(\G'\times_k\G_a'')\lra \omh T_{(x,y',z)} 
(\G*_h\G'\times_k\G_a'')$.
\begin{lem}
$i_8(id\mt i_7)=i_6(i_5\mt id)$\vs\\
Proof: {\em As above choose $B_x,B_{y'}$ such that: $T_x\G=T_x^l\G\oplus B_x$
and $T_y\G'=T_y^r\G'\oplus B_{y'}$. Then
$V:=T_{(x,y',z)}(\G*_h\G'\times_k\G_a'')$ has the decomposition: $V=V_1\oplus
V_2\oplus V_3$ for $V_1:=\{(\dot{x},\dot{y}',\dot{z})\in V\,:\,\dot{x}\in
B_x,\dot{y}\in B_{y'},\dot{z}\in T_z^r\G''\}\,,\,
V_2:=\{(\dot{x},\dot{y},0)\in V\,:\,\dot{x}\in 
B_x,\dot{y}'\in T_{y'}^r\G'\}$ and $V_3:=\{(\dot{x},0,0)\in V\,:\,\dot{x}\in
T_x^l\G\}$. And the lemma follows in the same way as the previous one.}\\
\dowl
\end{lem}

\noindent
{\bf C.} The set $\G*_h\G'*_k\G_a'':=\{(x,y',z')\in
\G\times\G'\times\G_a''\,:\, \el(x)=f_h(\el(y'))\,,\,\el(y')=f_k(\el(z'))\}$ 
 is a submanifold in $\G\times\G'\times\G_a''$.\\
The mappings: $\G*_h\G'*_k\G_a''\ni(x,y',z')\mapsto (y',z')\in \G'*_k\G_a''$ 
and $\G*_h\G'*_k\G_a''\ni(x,y',z')\mapsto (x,z')\in\G*_{k h}\G_a''$ are
surjective submersions. These fact can be seen in the same way  as above.  
Again we have isomorphisms:\\
$i_9:\,\oml(y')\mt\omr(z')\lra\omh T_{(y',z')}(\G*_k\G_a'')$.\\
$i_{10}:\,\oml(x)\mt\omh T_{(y',z')}(\G*_k\G_a'')\lra \omh T_{(x,y',z')}
(\G*_h\G'*_k\G_a'')$.\\
$i_{11}:\,\oml(x)\mt\omr(z')\lra \omh T_{(x,z')}(\G*_{k h}\G_a'')$.\\
$i_{12}:\, \omh T_{(x,z')}(\G*_{k h} \G_a'')\mt\oml(y')\lra \omh
T_{(x,y',z')} (\G*_h\G'*_k\G_a'')$.
\begin{lem}
$i_{10}(id\mt i_9)=i_{12}(i_{11}\mt id)(id\mt\sim)$, where 
$\sim: \oml(y')\mt\omr(z')\lra \omr(z')\mt\oml(y')$ is the flip.\vs\\
Proof: {\em The proof is based on the same arguments as above.}\\
\dowl
\end{lem}

\noindent
{\bf D.} $\G\timh\G'*_k\G_a'':=\{(x,y,z'')\in
\G\times\G'\times\G_a''\,:\, \er(x)=f_h(\el(y))\,,\,\el(y)=f_k(\el(z''))\}$ --
this is submanifold in $\{(x,y,z'')\in\G\times\G'\times\G_a''\,:\,(x,z'')\in
\G\times_{kh}\G_a''\}$ and in $\G\times(\G'*_k\G_a'')$. \\
The mappings: $\G\timh\G'*_k\G_a'' \ni(x,y,z'')\mapsto (y,z'')\in
\G'*_k\G_a''$ and $\G\timh\G'*_k\G_a''\ni(x,y,z'')\mapsto (x,z'')\in
\G\times_{k h}\G_a''$ are surjective submersions.  
And again we have isomorphisms:\\
$i_{13}:\,\omr(x)\mt\omh T_{(y,z'')}(\G'*_k\G_a'')\lra \omh
T_{(x,y,z'')}(\G\timh\G'*_k\G_a'')$.\\ 
$i_{14}:\,\oml(y)\mt\omr(z'')\lra \omh T_{(y,z'')}(\G'*_k \G_a'')$.\\
$i_{15}:\,\oml(y)\mt\omh T_{(x,z'')}(\G\times_{k h}\G_a'')\lra \omh
T_{(x,y,z'')}(\G\timh\G'*_k\G_a'')$.\\ 
$i_{16}:\, \omr(x)\mt\omr(z'')\lra\omh T_{(x,z'')} (\G\times_{k h}\G_a'')$.
\begin{lem} $i_{13}(id\mt i_{14})=i_{15}(id\mt i_{16})(\sim\mt id)$\vs\\
Proof: {\em As above.}\\
\dowl
\end{lem}

\noindent
{\bf E.} The mapping $(\got_h\times id): \G\timh\G'\times_k\G_a''\lra
\G*_h\G'\times_k\G_a''$ is a diffeomorphism.\\
The mapping $(id\times \got_k): \G*_h\G'\times_k\G_a''\lra \G*_h\G'*_k\G_a''$
is a diffeomorphism.\\
The mapping $(id\times \got_k): \G\timh\G'\times_k\G_a''\lra
\G\timh\G'*_k\G_a''$ is a diffeomorphism.\\
Define the mapping $\tilde{\got}: \G\timh\G'*_k\G_a''\ni(x,y,z)\mapsto
(x,m_h(x,y),m_{k h}(x,z))\in  \G*_h\G'*_k\G_a''$ .\\
Then: $(id\times \got_k)\,(\got_h\times id)=\tilde{\got}\,(id\times \got_k)\,:
\G\timh\G'\times_k\G_a''\lra \G*_h\G'*_k\G_a''$. So $\tilde{\got}$ is a 
diffeomorphism.

\begin{figure}[tbhp]
\centering
\fbox{\includegraphics[height=0.25\textheight,width=0.9\textwidth]
{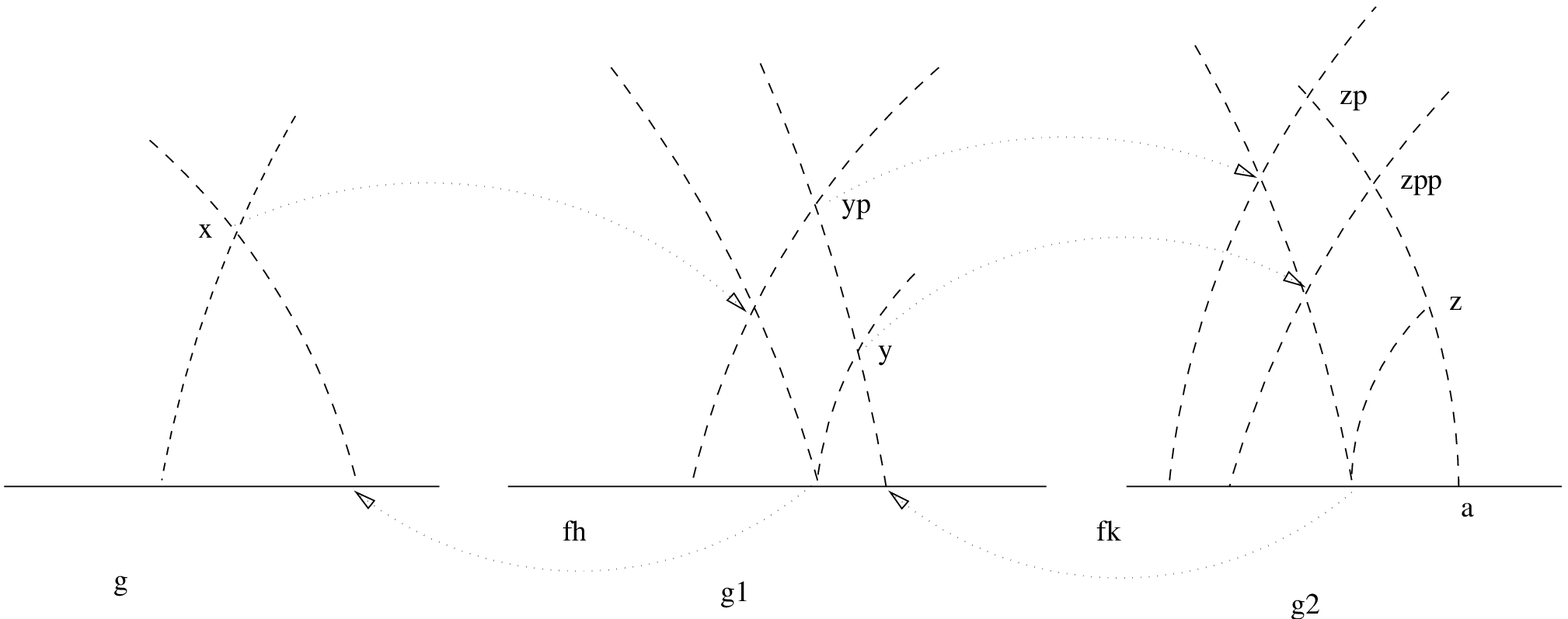}}
\caption{\label{propfig}}
\end{figure}

\noindent
{\bf F.} Let $(x,y,z)\in \G\timh\G'\times_k\G_a''$ and let $y',z',z''$ be 
defined by: $\got_h(x,y)=:(x,y')\,,\,\got_k(y',z)=:(y',z')\,,\,
\got_k(y,z)=:(y,z'')$. The situation is illustrated on the figure 
\ref{propfig}.

\noindent
From the definition we have: $\got_h i_1(\ro(x)\mt\ro(y))=:
t_h(x,y)i_5(\lo(x)\mt\ro(y'))\,$, \\
$\got_k i_7(\ro(y')\mt\ro(z))=:t_k(y',z)i_9(\lo(y')\mt\ro(z'))$, \\
$\got_k i_3(\ro(y)\mt\ro(z))=:t_k(y,z)i_{14}(\lo(y)\mt\ro(z''))$ 
and $\got_{k h}i_{16}(\ro(x)\mt\ro(z''))=:t_{k h}(x,z'')
i_{11}(\lo(x)\mt\ro(z'))$. 

\begin{lem} \,a)\,\,$(\got_h\times id) i_2 (\alpha\mt \rho_z)=
i_6(\got_h\alpha\mt\rho_z)$ \hspace{0.1\textwidth} b)\,\,
$(id\times \got_k)i_8(\lambda_x\mt\beta)=i_{10}(\lambda_x\mt \got_k\beta)$\\
c)\,\,$(id\times \got_k)i_4(\rho_x\mt\beta)=i_{13}(\rho_x\mt \got_k\beta)$
\hspace{0.1\textwidth}
d)\,\,$\tilde{\got}i_{15}(\lambda_y\mt\beta)=
i_{12}\sim(h_{\el'(y)}^R(x)\lambda_y\mt \got_{k h}\beta)$.\vs \\
Proof: {\em a)\,\,Let $V:=T_{(x,y,z)}(\G\timh\G'\times_k\G_a'')$ and
$V_1,W_1\subset$ be as in {\bf A}. Then $(\got_h\times id) W_1$ is a kernel of
the projection $(x,y',z)\mapsto z$. Since $\got_h\times id$ is a
diffeomorphism, 
$(\got_h\times id) V_1$ is complementary to this kernel. 
Compute the left hand side: 
$$(\got_h\times id)i_2(\alpha\mt\rho_z)((\got_h\times id)v_1\mt (\got_h\times
id) w_1)=i_2(\alpha\mt\rho_z)(v_1\mt w_1)=\alpha(\pi_{1 2}w_1)\rho_z(\pi_3
v_1).$$ And the right hand side: 
$$i_6(\got_h\alpha\mt\rho_z)((\got_h\times id)v_1\mt (\got_h\times id)
w_1)=\got_h\alpha(\pi_{1 2} (\got_h\times id)w_1)\rho_z(\pi_3(\got_h\times
id)v_1)=$$
$$=\got_h\alpha (\got_h\pi_{1 2} w_1) \rho_z(\pi_3 v_1)=\alpha(\pi_{ 1
2}w_1)\rho_z (\pi_3 v_1).$$ This proves a). 
Statements b) and c) can be proved in the same way.

Let us prove d). Let $V:=T_{(x,y,z'')}(\G\timh\G'*_k\G_a'')\,$,   
$V_2:=\{(0_x,\dot{y},0_{z''})\in V\,:\,\dot{y}\in T_y^l\G'\}$ and $V_1$ be 
such that $V=V_1\oplus V_2$. The isomorphism $i_{15}: 
\oml(y)\mt\omh T_{(x,z'')} (\G\times_{kh}\G_a'')\lra 
\omh T_{(x,y,z'')}(\G\timh\G'*_k\G_a'')$ is given by: 
$i_{15}(\lambda_y\mt\beta)(v_1\dz v_2):=
\lambda_y(\pi_2 v_2)\beta(\pi_{13}v_1)$.\\
Let also $\tilde{V}:=T_{(x,y',z')}(\G*_h\G'*_k\G_a'')\,,\,
\tilde{V_2}:=\{(0_x,\dot{y'},0_{z'})\in\tilde{V}\,:\,\dot{y'}\in 
T_{y'}^l\G'\}$ and $\tilde{V_1}$ be complementary to $\tilde{V_2}$. Then 
$i_{12}:\omh T_{(x,z')}(\G*_{kh}\G_a'')\mt\oml(y')\lra\omh T_{(x,y',z')}
(\G*_h\G'*_k\G_a'')$ is given by: $i_{12}(\beta\mt\lambda_{y'})(\tilde{v_1}\dz
\tilde{v_2}):=\beta(\pi_{13}\tilde{v_1})\lambda_{y'}(\pi_2\tilde{v_2})$.\\
We have: 
$\tilde{\got}i_{15}(\lambda_y\mt\beta)(\tilde{\got}v_1\dz\tilde{\got}v_2):=
\lambda_y(\pi_2 v_2)\beta(\pi_{13}v_1)$.\\
 It is easy to see that 
$\tilde{\got}V_2=\tilde{V_2}$ so $\tilde{\got}V_1$ is complementary to $\tilde{V_2}$.
From this we have:  
$$i_{12}(\got_{kh}\beta \mt h_{\el'(y)}^R\lambda_y)
(\tilde{\got}v_1\dz\tilde{\got}v_2)=
(h_{el'(y)}^R\lambda_y)(\pi_2\tilde{\got}v_2)
(\got_{kh}\beta)(\pi_{13}\tilde{\got} v_1)=$$
$$=(h_{el'(y)}^R\lambda_y)(h_{el'(y)}^R\pi_2v_2)
(\got_{kh}\beta)(\got_{kh}\pi_{13}v_1)=\lambda_y(\pi_2 v_2)\beta(\pi_{13}v_1).
$$
 Where we used:  $\pi_{13}\tilde{\got}=\got_{kh}\pi_{13}$ and
$\pi_2\tilde{\got}v_2=h_{\el'(y)}^R\pi_2 v_2$.}\\
\dowl
\end{lem}
\FloatBarrier
\noindent
{\em Proof of the Lemma \ref{app}:}
$$(id\times \got_k)(\got_h\times id)i_2(i_1\mt id)(\ro(x)\mt\ro(y)\mt\ro(z))=
(id\times \got_k)i_6(\got_h i_1(\ro(x)\mt\ro(y))\mt\ro(z))=$$
$$=t_h(x,y)(id\times \got_k)i_6(i_5(\lo(x)\mt\ro(y'))\mt\ro(z))=
t_h(x,y)(id\times \got_k)i_8(id\mt i_7)(\lo(x)\mt\ro(y')\mt\ro(z))=$$
$$=t_h(x,y)i_{10}(\lo(x)\mt \got_k i_7(\ro(y')\mt\ro(z)))=
t_h(x,y)t_k(y',z)i_{10}(id\times i_9)(\lo(x)\mt\lo(y')\mt\ro(z')).$$
From the other side:
$$\tilde{\got}(id\times \got_k)i_2(i_1\mt id)(\ro(x)\mt\ro(y)\mt\ro(z))=
\tilde{\got}(id\times \got_k)i_4(id\mt i_3)(\ro(x)\mt\ro(y)\mt\ro(z))=$$
$$=\tilde{\got}i_{13}(\ro(x)\mt \got_k i_3(\ro(y)\mt\ro(z)))=
t_k(y,z)\tilde{\got} i_{13}(\ro(x)\mt i_{14}(\lo(y)\mt\ro(z'')))=$$
$$=t_k(y,z)\tilde{\got} i_{15}(\lo(y)\mt i_{16}(\ro(x)\mt \mt\ro(z'')))= 
t_k(y,z) i_{12}\sim (h_{\el'(y)}^R(x)\lo(y)\mt \got_{k h} 
i_{16}(\ro(x)\mt\ro(z'')))=$$
$$=t_k(y,z)t_{k h}(x,z'') i_{12}\sim (\lo(y')\mt i_{11}(\lo(x)\mt\ro(z')))=
t_k(y,z)t_{k h}(x,z'') i_{12}(i_{11}(\lo(x)\mt\ro(z')\mt \lo(y'))=$$
$$=t_k(y,z)t_{k h}(x,z'') i_{12}(i_{11}\mt id)(id\mt
\sim)(\lo(x)\mt\lo(y')\mt\ro(z'))= t_k(y,z)t_{k h}(x,z'') i_{10}(id\mt i_9)
(\lo(x)\mt\lo(y')\mt\ro(z')).$$
So $t_h(x,y)t_k(y',z)=t_k(y,z)t_{k h}(x,z'')$.\\
\dowl

\section{$C^*$-algebra of a differential groupoid.}
\newcommand{\cred}{\mbox{$C^*_{red}(\G)$}}
\newcommand{\hfi}{\mbox{$\hat{\varphi}$}}
In this section we define the $C^*$-algebra of a differential groupoid and
show that this correspondance is a covariant functor to the category of $C^*$
-algebras.

Let $\G$ be a differential groupoid. In the previous section we have
established the following facts:\\
1)  With any morphism $h:\G\rel\G'$ one can associate two objects: the
mapping $\hat{h}: \sA(\G)\lra LM\sA(\G')$ and the non degenerate
representation $\pi_h$ of the *-algebra $\sA(\G)$ in $L^2(\G')$. \\
2) $\om_3^*(\hat{h}(\om_1)\om_2)=(\hat{h}(\om_1^*)\om_3)^*\om_2$ for any
$\om_1\in \sA(\G)\,,\,\om_2,\om_3\in\sA(\G')$.\\
3) If $k:\G'\rel \G''$ is a morphisms of differential groupoids. 
Then: \\
a) $\hat{k}(\hat{h}(\om_1)\om_2)\om_3=\hat{kh}(\om_1)(\hat{k}(\om_2)\om_3) 
\,,\, \om_1\in\sA(\G)\,,\,\om_2\in\sA(\G')\,,\, \om_3\in \sA(\G'')$.\\
b) $\pi_k(\hat{h}(\om_1)\om_2) \psi=
\pi_{kh}(\om_1)(\pi_k(\om_2)\psi)$, for
$\om_1\in\sA(\G)\,,\,\om_2\in\sA(\G')\,,\, \psi\in L^2(\G'')$.\\
4) We can choose norm on $\sA(\G)$ which agrees with *-algebra structure,
such that:  $||\pi_h(\om)||\leq ||\om||$. \\
5) Exists morphism $l$ from $\G$ to the pair groupoid $\G\times \G$ such that
$\pi_{l}$ is a faithfull representation of $\sA(\G)$.

From 4) and 5) follows that the following definition is meaningfull: 
\begin{defi} The $C^*$-algebra of a differential groupoid $\G$ {\em is a
completion of $\sA(\G)$ with respect to the norm:
$||\om||:=sup_h||\pi_h(\om)||$, where the supremum is taken over all morphisms
$h:\G\rel \G'$.}\end{defi} 

\begin{prop}
For any morphism $h:\G\rel\G'\,$, $\pi_h$ extends to a nondegenerate 
representation of $C^*(\G)$ and $\hat{h}$ extends to $C^*(h)\in
Mor(C^*(\G),C^*(\G')).$\vs \\  
Proof: {\em Let us recall (cf. appendix E) that $\varphi$ is a morphism of
$C^*$-algebras $A,B$ iff $\varphi: A\lra M(B)$ is a *-algebra homomorphism and
the set $\varphi(A)B$ is dense in $B$. Let $k:\G'\rel\G''$ be any morphism. 
From  prop. \ref{rep} one has: 
$$||\pi_k(\hat{h}(\om_1)\om_2)||\leq ||\pi_{kh}(\om_1)||||\pi_k(\om_2)|| \leq
||\om_1||_{C^*(\G)}||\om_2||_{C^*(\G')},$$ so also:
$||\hat{h}(\om_1)\om_2||_{C^*(\G')}\leq
||\om_1||_{C^*(\G)}\,||\om_2||_{C^*(\G')}$. \\
Since $\sA(\G')$ is dense in $C^*(\G')$ 
$\hat{h}(\om_1)$ can be extended to a continous linear mapping on $C^*(\G')$
and by the density of $\sA(\G)$ in $C^*(\G)$ $\hat{h}$ defines a bounded
algebra homomorphism $C^*(h):C^*(\G)\lra BC^*(\G')$. \\
From  \ref{star} :
$\om_3^*(\hat{h}(\om_1)\om_2)=(\hat{h}(\om_1^*)\om_3)^*\om_2$ for any 
$\om_1\in \sA(\G)\,,\,\om_2,\om_3\in\sA(\G')$.\\ 
By density of $\sA(\G')$ and continouty this equality extends to
$\om_2,\om_3\in C^*(\G')$. So for any $\om\in\sA(\G)$  $\hat{h}(\om)$ has the
hermitian conjugate (see appendix E) equal to $\hat{h}(\om^*)$ and
$C^*(h)(\om)\in MC^*(\G')$. It follows that $C^*(h)$ is a continous
*-algebra homomorphism $\sA(\G)\lra MC^*(\G')$. By continouty it is also true
for $\om\in C^*(\G)$.\\
From the Prop. \ref{apjed} for any $\om'\in\sA(\G')$ one can find $\om_n\in
\sA(\G)$ such that $||\hat{h}(\om_n) \om'- \om'||_0\lra 0$ for *-norm on
$\sA(\G')$ given by some $\om_0'$. So also 
$||\hat{h}(\om_n) \om'- \om'||_{C^*(\G')}\lra 0$ and
$\sA(\G')\subset\overline{C^*(h)(\sA(\G))\sA(\G'))}$. Since $\sA(\G')$ is
dense in $C^*(\G')$ this is nondegeneracy condition. \\
In this way $C^*(h)\in Mor(C^*(\G),C^*(\G')).$\\
From the definition $||\pi_h(\om)||\leq ||\om||_{C^*(\G)}$ for $\om\in
\sA(\G)$. So $\pi_h$ defines representation of $C^*(\G)$. The nondegeneracy
condition is just the Col. \ref{nondeg2}.}\\
\dow
\end{prop}
The next proposition shows the functoriality of $C^*$.
\begin{prop}
$C^*(k h)=C^*(k) C^*(h)$.\vs\\
Proof: {\em Let $\phi_1:=C^*(h)\,,\,\phi_2:=C^*(k)\,,\,\phi_3:=C^*(k h)$. Let
$\hat{\phi_2}: M C^*(\G')\lra M C^*(\G'')$ denote the unique extension of
$\phi_2$ (see appendix E). We have to show that $\phi_3=\hat{\phi_2}\phi_1$.
In fact it is enaugh to show the equality for all $\om_1\in \sA(\G)$: 
$\phi_3(\om_1)=\hat{\phi_2}\phi_1(\om_1)$. Since $\hat{k}(\sA(\G'))\sA(\G'')$
is dense in $C^*(\G'')$ it is enaugh to show that: 
$\hat{k h}(\om_1)\hat{k}(\om_2)\om_3=
\hat{\phi_2}\hat{h}(\om_1)\hat{k}(\om_2)\om_3= 
\hat{k}(\hat{h}(\om_1)\om_2)\om_3$. But this is Prop. \ref{funktor}.}\\
\dow
\end{prop}
The definition of the $C^*$ norm given above is rather abstract and 
obtained $C^*$ algebra seems untreatable. However, as is shown below, we
can restrict ourself to a smaller class of morphisms--namely to the morphisms
to the pair groupoids.\\
Let $h:\G\rel\G'$ be a morphism of differential groupoids and let
$l':\G'\rel\G'\times\G'$ be a left regular representation as defined in 
\ref{exmor} f) ($\G'\times\G'$ is a pair groupoid, not the product of 
groupoids). Then $\tilde{h}:=l'h$ is a morphism from $\G$ to the pair groupoid
$\G'\times\G'$. Then it is easy to see that: $f_{\tilde{h}}=f_h\el'\,,\,
\G\times_{\tilde{h}}(\G'\times\G')=(\G\timh\G')\times\G'\,,\,
\G*_{\tilde{h}}(\G\times\G')=(\G*_h\G')\times \G'\,,\,
m_{\tilde{h}}=m_h\times id$ and $\got_{\tilde{h}}=\got_h\times id$. Taking
into account these formulas and equality 
$L^2(\G'\times\G')=L^2(\G')\mt L^2(\G')$ one can show that:
 $\pi_{\tilde{h}}(\om)=\pi_h(\om)\mt I$ so (semi)norm coming from $\pi_h$ and
$\pi_{\tilde{h}}$ are equal. Since all manifolds are second countable, our
$C^*$ algebra is separable.
\begin{ex}
 a) Pair groupoids. {\em Let $\G:=X\times X\,,\,\G':=Y\times Y$ be  pair 
groupoids. Due to the structure of morphisms in this case as explained in 
Example \ref{exmor2}b) we can assume that $Y=X\times Z$ for some manifold $Z$
and $h:\G\rel\G'$ is given by $Gr(h):=\{(x,z,x',z;x,x')\,:\,x,x'\in X\,,\,
z\in Z\}$. Choose $\varphi_0$--smooth, real, non vanishing half density on $X$.
Since $T_{(x,x')}^l\G=T_{x'} X$ and $T_{(x,x')}^r\G=T_x X$ it defines $\lo$ by
the formula $\lo(x,x')=\varphi_0(x')$. Then the corresponding $\om_0$ is 
given by $\om_0(x,x')=\varphi_0(x')\mt\varphi_0(x)$. Choose also 
$\mu_0$--smooth, real, non vanishing half density on $Z$. Then 
$\psi_0(x,z):=\varphi_0(x)\mt\mu_0(z)$ is smooth, real,non vanishing half 
density on $Y$. Again $\tilde{\psi}_0(y_1,y_2):=
\psi_0(y_1)\mt\psi_0(y_2)$ is smooth, real,
non vanishing half density on $\G'$, moreover $\psi$ defines right invariant, 
half density on $\G'$ along the right fibers, so this is decomposition of 
$\tilde{\psi}_0$ into the form $\tilde{\psi}_0=\ro\mt\nu_0$ -- $\nu_0=\psi_0$ 
is half density on the set of identities of $\G'$. It is easy to see that:
$$\G\timh\G'=\{(x,x';x',z',x'',z''):x,x',x''\in X\,z',z''\in Z\},$$
$$\G*_h\G'=\{(x,x';x,z,x'',z''):x,x',x''\in X\,z,z''\in Z\},$$
$$\got_h(x,x';x',z',x'',z'')=(x,x';x,z',x'',z'').$$ 
Also a short computation shows
that with the above choice of $\om_0,\psi_0$ the formula function $t_h$ is 
is constant and equal to $1$. So for $\om=f\,\om_0\,,\,
\Psi=\psi\tilde{\psi}_0$
the representation associated with the morphism $h$ is given by:
$$(\pi_h(\om)\Psi)(y_1,y_2)=\left[\int_{X}\varphi_0^2(x) f(x_1,x)\psi(y,y_2)
\right]\tilde{\psi}_0(y_1,y_2),$$ 
where $y:=(x,z_1)\,,\,y_1=:(x_1,z_1)$.\\
We have $ L^2(\G')=L^2(Y\times Y)=L^2(Y)\mt L^2(Y)$ and from the formula it is
clear that $\pi_h=\pi_1\mt I$ for 
$$(\pi_1(\om)\Psi)(y_1):=\left[\int_{X}\varphi_0^2(x) f(x_1,x)\psi(y)\right]
\psi_0(y_1)$$ for $\Psi=\psi\psi_0$--smooth, half density on $Y$ 
with compact support. Again since $Y=X\times Z$ this representation is of 
the form $\pi_1=\pi_0\mt I$ for 
$$(\pi_0(\om)\Phi)(x_1):=\left[
\int_{X}\varphi_0^2(x) f(x_1,x) \phi(x) \right]\varphi_0(x_1),$$ where 
$\Phi=:\phi\varphi_0$ is smooth, half density on $X$ with compact support.\\
In this way we have shown, that for pair groupoids  the $C^*$ norm on $
\sA(\G)$ is equal to the norm coming from the left regular representation. 
The completion of $\sA(\G)$ in this norm is the algebra of compact operators.}
\vs\\
Groups. {\em Let $\G=G$ be a Lie group. $\sA(G)$ is by definition *-algebra of
compactly supported, smooth, complex densities on $G$.  For any $U$-strongly 
continous, unitary representation of $G$ on the Hilbert space $H$ the formula 
$$(x,\pi_U(\nu) x):=\int_G \nu(g)(x,U(g) x)\,,\,\nu\in\sA(G),x\in H$$ 
defines nondegenerate *-representation of $\sA(G)$. $C^*(G)$ is a 
completion of 
$\sA(G)$ with respect to the norm $||\nu||:=sup_{\pi_U}||\pi_U(\nu)||.$

Since for a group, left and right fibers are equal it is clear that any 
$\om\in\sA(\G)$ is  a density on $G$ by an assigment: 
$\sA(\G)\ni\lambda\mt\rho\mapsto \lambda \rho\in \sA(G)$. Let $X$ be a manifold
and $\G':=X\times X$ corresponding pair groupoid. By a slight modification of
arguments used in Example. \ref{exmor2} a) morphisms $h:\G\rel\G'$ are in one
to one correspondance with smooth actions of $G$ on $X$: $Gr(h):=\{(gx,x;g)\,
:\,x\in X,G\in G\}$. It is easy to see that $\G\timh\G'=\G*_h\G'=\{(g;x_1,x_2)
\,:\,x_1,x_2\in X,g\in G\}$ and $\got_h(g;x_1,x_2)=(g;gx_1,x_2)$. Moreover
$T_{(g;x,x_0)}(\G\timh \G_{x_0}')=T_gG\ms T_xX$. Choose $\lo$ -- real, smooth,
non vanishing left invariant half density on $G$ and $\psi_0$--smooth, real, 
non vanishing half density on $X$. Then a short computations show that 
$t_h(g;g^{-1}x,x_0)=\frac{\ro}{\lo}(g)\frac{g \psi_0}{\psi_0}(x)$ and for
$\om=f\lo\mt\ro\,,\,\Psi=\psi\psi_0\mt\psi_0$ we have:
$$(\pi_h(\om)\Psi)(x_1,x_2)=\left[\int_G \lo^2(g) f(g)t_h(g;g^{-1}x_1,x_2)
\psi(g^{-1}x_1,x_2)\right]\psi_0(x_1)\mt\psi_0(x_2).$$ 
Since $L^2(\G')=L^2(X)\mt L^2(X)$ this representation is of form $\pi_h(\om)=
\tilde{\pi}_h(\om)\mt I$ for
$$(\tilde{\pi}_h(\om)\Psi)(x):=
\left[\int_G\lo^2(g) f(g) \frac{\ro}{\lo}(g)\frac{g \psi_0}{\psi_0}(x)
\psi(g^{-1}x)\right] \psi_0(x),$$ 
where $\Psi=\psi\psi_0$ is a smooth half  density on 
$X$ with compact support.

From the other side, the action of $G$ on $X$ defines strongly continous 
unitary representation of $G$ on $L^2(X)$ by the formula: $U_g\Psi:=g \Psi$ 
for $\Psi$-- smooth, compactly supported half density on $X$. If $\Psi=\psi
\psi_0$ then $(U_g\Psi)(x)=\psi(g^{-1}x) \frac{g \psi_0}{\psi_0}(x)
\psi_0(x)$. If $\om=f\lo\mt\ro$ then $\nu:=f\lo\ro=f\frac{ro}{\lo}\,\lo^2$ and
$\tilde{\pi}_h(\om)\Psi=\pi_U(\nu)\Psi$. In this way 
$||\om||_{C^*(\G)}\leq||\nu||_{C^*(G)}$. It is clear that if $\pi_l$ is a left
regular representation then $||\pi_l(\om)||\leq||\om||_{C^*(\G)}$ so the 
$C^*$ algebra of a Lie group $G$ viewed as a differential groupoid is something
``between'' the reduced $C^*$ algebra of $G$ and the algebra $C^*(G)$ where 
$G$ is treated as locally compact topological group.}\vs\\
Transformation groupoids. {\em Let $\G:=G\times X$ be a transformation 
groupoid. By $C_0(X)$ we denote the $C^*$ algebra of complex, continous,
vanishing at infinity functions on $X$. The action of $G$ on $X$ induces a
strongly continous homomorphism $\alpha: G\ni g\mapsto \alpha_g\in 
Aut(C_0(X))$, where $Aut(C_0(X))$ is a group of *-isomorphism of $C_0(X))$,
namely $(\alpha_g f)(x):=f(g^{-1}x)$. 
So $(G,C_0(X),\alpha)$ is a $C^*$ dynamical system (see appendix). Let 
$Y$ be a manifold and $\G':=Y\times Y$ the corresponding pair groupoid.
By Example \ref{exmor2} a) morphisms $h:\G\rel\G'$ are in one to one 
correspondance with smooth actions $G\times Y\ni(g,y)\mapsto g y\in Y$ together
with smooth equivariant mapping $F:Y\lra X$. The graph of $h$ is then equal: 
$Gr(h)=\{(gy,y;g,F(y))\,:,y\in Y,\,g\in G\}$. The action of $G$ on $Y$ induces
strongly contionous, unitary representation $G\ni g\mapsto U_g\in B(L^2(Y))$.
The mapping $F$ defines nondegenerate representation $\pi$ of $C_0(X)$ on 
$L^2(Y)$  by the formula: $(\pi(f)\psi)(y):=f(F(y))\psi(y)\,,\,f\in C_0(X)$ 
and $\psi$-smooth, compactly supported half density on $Y$. The pair $(\pi,U)$
is a covariant representation of $(G,C_0(X),\alpha)$. Indeed 
$$(\pi(\alpha_g f)U_g \psi)(y)=(\alpha_g f)(F(y))(g\psi)(y)=
f(g^{-1}F(y))(g\psi)(y)=$$
$$=f(F(g^{-1}y))(g\psi)(y)=(\pi(f))(g^{-1}y)(g\psi)(y)=
(U_g\pi(f)\psi)(y).$$

Now we go back to groupoid $\G$. Choose $\mu_0\neq 0$-real, half density on
$T_eG$. Since $T_{(g,x)}^r\G=T_gG$, it defines right invariant, non vanishing,
half density on $\G$ by the formula $\ro(g,x)(v_g):=\mu_0(v_gg^{-1})$ where 
$v_g\in \lma T_gG=\lma T_{(g,x)}^r\G$. The corresponding left invariant
half density is given by $\lo(g,x)(v)=\mu_0(g^{-1}\pi(v))$ where 
$v\in \lma T_{(g,x)}^l\G$ and $\pi:\G\times X\lra G$ is a projection.
Let $\lambda$ be a corresponding left invariant density on $G$ i.e. 
$\lambda(g):=g\mu_0^2$ and let $\D$ be a corresponding modular function.
Then for $\om_1=f_1\om_0,\,\om_2=f_2\om_0 \in \sA(\G)$ we have:
$$(\om_1\om_2)(g,x)=\left[\int_{G}\lambda(g_1)f_1(g_1,g_1^{-1}gx)
f_2(g_1^{-1}g,x)\right] \om_0(g,x)\,\,{\rm and\,\,} 
\om^*(g,x)=\overline{f(g^{-1},gx)}\om_0(g,x).$$
Let $K(G,C_0(X))$ be a *-algebra of compactly supported, continous functions 
from  $G$ to $C_0(X)$ with the usual structure.(see appendix). 
We define the mapping
$\sA(\G)\ni \om \mapsto \hat{\om} \in K(G,C_0(X))$ by 
$\hat{\om}(g)(x):=\D(g)^{-1/2}f(g,g^{-1}x)$ for $\om=f\om_0$. Straightforward
computations show that this is injective *-homomorphism.\\
Let $h:\G\rel\G'$ be a morphism, $Gr(h):=\{(g y,y;g,F(y))\}$ Choose $\psi_0$ -
smooth, real, non vanishing half density on $Y$. Then 
$(\psi_0\mt\psi_0)(y_1,y_2):=\psi_0(y_1)\psi_0(y_2)$ is real, smooth, non 
vanishing half density on $\G'$. We have: $\G\timh\G_y'=\{(g,F(y_1);y_1,y)\,:\,
g\in G,\,y_1,y\in Y\}$, $\G*_h\G_y'=\{(g,F(y_1);g y_1,y)\,:\,g\in G,\,y_1,
y\in Y\}$ and $\got_h(g,F(y_1);y_1,y)=(g,F(y_1),g y_1,y)$.
Simple computations show that the function $t_h$ is given by:
$t_h(g,F(y);y,y_1)=\D(g)^{-1/2}\frac{g \psi_0}{\psi_0}(g y)$. So the 
representation $\pi_h$ is given by:
$$(\pi_h(\om)\Psi)(y_1,y_2)=\left[\int_G\lambda(g) f(g,g^{-1}F(y))\D(g)^{-1/2}
\frac{g \psi_0}{\psi_0}(y_1)\psi(g^{-1}y_1,y_2)\right]
\psi_0(y_1)\mt\psi_0(y_2).$$ 
Again we see that the norm of $\pi_h(\om)$ is equal to the norm of 
$\tilde{\pi}_h(\om)$ where $\tilde{\pi}_h$ is a representation on $L^2(Y)$ 
given by: 
$$(\tilde{\pi}_h(\om)\Psi)(y):=\left[\int_G\lambda(g) f(g,g^{-1}F(y)) 
\D(g)^{-1/2} \psi(g^{-1}y)\frac{g \psi_0}{\psi_0}(y)\right]\psi_0(y).$$ 
Let $\rho$ be a representation on $K(G,C_0(X))$ associated with a covariant
representation $(\pi,U)$ defined by morphism $h$. Then 
$$(\rho(\hat{\om})\Psi)(y)=\int_G \lambda(g) (\pi(\hat{\om}(g))U(g)\Psi)(y)=
\int_G\lambda(g)\hat{\om}(g)(F(y)) (g\Psi)(y)=$$
$$=\left[\int_G\lambda(g) \D(g)^{-1/2}
f(g,g^{-1}F(y))\frac{g \psi_0}{\psi_0}(y)\psi(g^{-1}y)\right]\psi_0(y).$$
So $||\rho(\hat{\om})||=||\pi_h(\om)||.$ And again $C^*(G\times X)$ is 
a kind of ``smooth''  crossed product, which is ``smaller'' then universal
crossed product $C^*(G,C_0(X),\alpha)$. }\end{ex}

\noindent
{\bf Bissections as multipliers.}\vs

Let $\G$ be a differential groupoid, $\omega=\lambda\mt\rho$ be a bidensity
and $B$ bissection. We define action of $B$ on $\omega$ by:
$$(B\omega)(Bx)(Bv\mt Bw):=\omega(x)(v\mt w),$$ 
where $v\in\lma T^l_x\G\,,\,w\in
\lma T^r_x\G$. \\
This is well defined since $BT^l_x\G=T^l_{Bx}\G$ and
$BT^r_x\G=T^r_{Bx}\G$. Choose some $\om_0$ and for a bissection $B$ define
function $b:\G\lra R$ by $B\om_0=:b \om_0$. Then $b$ is nonvanishing and
smooth. If $\om=f \om_0$ then $B\om=B(f) \om_0$ for $(Bf)(x)=f(B^{-1}x) b(x)$.
\\ Since $B$ acts  as a diffeomorphism of $\G$, it defines an unitary operator
on $L^2(\G)$ by: $(B\psi)(B x)(B v):=\psi(x)(v)\,,v\in\lma T_x\G$ and 
$\psi$- smooth, half density on $\G$. Choosing $\Psi_0:=\ro\mt\nu_0$ as in 
Prop. \ref{rep} we have $B\Psi_0=b \Psi_0$ for the same function $b$.
\begin{lem} \label{bis-gw}
Let $B$ be a bissection and $\om_1,\om_2\in\sA(\G)$, then
$\om_1^*(B\om_2)=(s(B)\om_1)^*\om_2$.\vs\\
Proof: {\em From $(B_1 B_2)\om=B_1(B_2 \om)$ and $B s(B)=\G^0$ we have
$s(b)(x)=\frac{1}{b(Bx)}$ where $s(b)$ is defined by $s(B)\om_0=:s(b) \om_0$.
Also from definition of $\om_0$ and $b$ is easy to see that in fact $b$ is
defined by: $$b(Bx)\ro(Bx)(B w)=\ro(x) (w)\,,\,w\in \lma T^r_x\G,$$ so $b(x
y)=b(x)$. With the usual notation, the LHS is equal:\\
$(\om_1^*(B\om_2))(x)=\left[\int_{F_l(x)}\lo^2(y)\,f_1^*(y)
(Bf_2)(s(y)x)\right]\om_0(x).$
We compute:
$$\int_{F_l(x)}\lo^2(y)\,f_1^*(y) (Bf_2)(s(y)x)=
\int_{F_l(x)}\lo^2(y)\,\overline{f_1(s(y))}\, b(s(y)x)\,f_2(B^{-1}(s(y)x))$$
And the RHS: $((s(B)\om_1)^*\om_2)(x)=
\left[\int_{F_r(x)}\ro^2(z)\,(s(B)f_1)^*(xs(z))\,f_2(z)\right]\om_0(x)$. 
$$\int_{F_r(x)}\ro^2(z)\,(s(B)f_1)^*(xs(z))\,f_2(z)=
\int_{F_r(x)}\ro^2(z)\,\overline{(s(B)f_1)(zs(x))}\,f_2(z)=$$
$$=\int_{F_r(x)}\ro^2(z)\,s(b)(zs(x))\overline{f_1(s(B)^{-1}zs(x))}\,f_2(z)=
\int_{F_r(a)}\ro^2(z')\,\frac{1}{b(Bz')}\overline{f_1(B z')}\,f_2(z'
x),$$
 where $a:=\el(x)$.\\
 Now since $B$ induces diffeomorphism $F_r(a)\lra
F_r(a)$ we can rewrite this integral as: 
$$\int_{F_r(a)}\ro^2(z')\,b(z')\overline{f_1(z')}\,f_2(B^{-1}z'x)=
\int_{F_l(x)}\lo^2(y)\,b(s(y))\overline{f_1(s(y))}\,f_2(B^{-1}s(y)x)$$ 
and since $b(s(y))=b(s(y)x)$ this is equal to LHS.}\\
\dowl
 \end{lem}
\begin{prop} Let $\G$ be a differential groupoids and $B$ a bissection of
$\G$. Then  for any morphism $h:\G\rel\G'$:\\
a) $(\hat{h}(B \om_1))\om_2=h(B)(\hat{h}(\om_1)\om_2)\,,\,\om_1\in\sA(\G),
 \,\om_2\in\sA(\G')$\\ 
b) $\pi_h(B\om_1)=h(B)\pi_h(\om_1)\,,\,\om_1\in\sA(\G)$.\vs\\
Proof: 
{\em Choose $\om_0\,,\,\om'_0$ and let $\om_1=f_1\om_0\,,\,\om_2=f_2\om_0'$. 
Let $B$ be a bissection. $B\om_1=(Bf_1)\om_0\,,\,(Bf_1)(x)=f_1(B^{-1}x) b(x)$
and $h(B)\om_2=(h(B)f_2)\om'_0\,,\,(h(B)f_2)(z)=f_2((h(B))^{-1}z)h(b)(z)$.\\
The left hand side of a): 
$$(\hat{h}( B\om_1)\om_2)(z)=
\left[\int_{F_l(c)}\lo^2(x) (B f)(x) t_h(x,y) f_2(y)\right]\,\om'_0(z),$$
$$\int_{F_l(c)}\lo^2(x) (B f)(x) t_h(x,y) f_2(y)=
\int_{F_l(c)}\lo^2(x) b(x) f(B^{-1}x) t_h(x,y) f_2(y),$$
where $\got_h(x,y)=(x,z)$.\\
And the right hand side: $(h(B)(\hat{h}(\om_1)\om_2))(z)=
h(b)(z)(f_1*_hf_2)((h(B))^{-1}z)\om'_0(z)$,
$$h(b)(z)(f_1*_hf_2)((h(B))^{-1}z)=h(b)(z) \int_{F_l(c')}\lo^2(x') f_1(x') 
t_h(x',y') f_2(y'),$$ 
where $\got_h(x',y')=(x',(h(B))^{-1} z)$.\\
Let $x_1\in B$ be such
that $\el(x_1)=c$. Then $x'\mapsto x_1 x'=B x'$ is a diffeomorphism
$F_l(c')\lra F_l(c)$. Using this fact we can rewrite the last expression as: 
$h(b)(z)  \int_{F_l(c)}\lo^2(x) f_1(s(x_1) x) 
t_h(s(x_1)x,y') f_2(y')$, where $\got_h(s(x_1)x,y')=(s(x_1) x,(h(B))^{-1}
z)$. The situation is illustrated on the figure:
\begin{figure}[tbhp]\label{morbis}
\psfrag{hbr}{$h_b^R$}
\p{mh}{$m_h(x,y)$}
\p{B}{$B$}
\p{hB}{$h(B)$}
\p{a1}{$a_1$}
\p{c1}{$c_1$}
\centering
\fbox{\includegraphics[height=0.25\textheight,width=0.8\textwidth]
{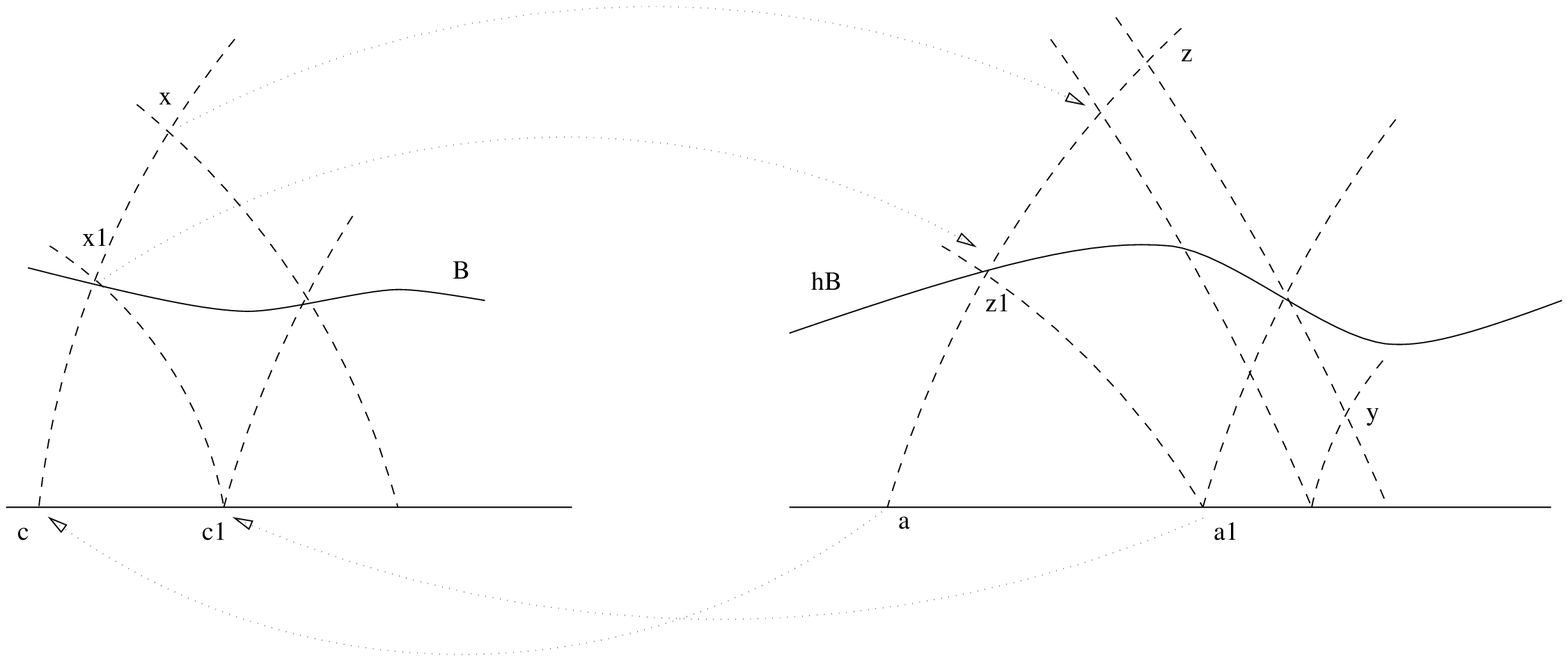}}
\caption{}
\end{figure}

\noindent
Now we use the following: 
\begin{lem} a) $\got_h(B\times id)=(B\times h(B))\got_h$\\
b) $b(B x) t_h(B x,y)=t_h(x,y) h(b)(h(B) m_h(x,y))$ for any
$(x,y)\in\G\timh\G'$.
\end{lem}
From the point a) we have $y'=y$ and from the point b): $b(x) t_h(x,y)=h(b)(z)
t_h(B^{-1}x,y)$. The statement a) is proven.\\
b) If we write $\Psi=\psi\Psi_0$ then the integrals appearing in the formula 
b) are the same as in a), so b) is proven.\vs\\}
Proof of the lemma: {\em a) Let $(x,y)\in\G\timh\G'\,,\,a:=\el'(y)$. 
$Bx=x_1 x$ for a
unique point $x_1\in B$. 
$$\got_h(B\times id)(x,y)=\got_h(x_1 x,y)=(x_1
x,h_a^R(x_1 x)y)= (x_1 x, h_b^R(x_1) h_a^R(x) y),$$ for $b:=\el'(h_a^R(x))$.
But this is equal to $$(x_1 x, h(B) h_a^R(x) y)=(B\times
h(B))(x,h_a^R(x)y)=(B\times h(B))\got_h(x,y).$$
b) Let $(x,y)\in \G\timh\G'\,,\,c:=\er'(y)$. $\ro(x)\mt\ro'(y)$ defines half
density on $T_{(x,y)}(\G\timh \G_c')$. 
$$\got_h(B\times
id)(\ro(x)\mt\ro'(y))=b(Bx)\got_h(\ro(Bx)\mt\ro'(y))=
b(Bx)t_h(Bx,y)(\lo(Bx)\mt\ro'(z')),$$ 
where $\got_h(Bx,y)=:(Bx,z')$, i.e. $z'=h(B) m_h(x,y)$. \\
From the other side: 
$$(B\times h(B))\got_h(\ro(x)\mt\ro'(y))=(B\times
h(B))t_h(x,y)(\lo(x)\mt\ro'(z))=t_h(x,y) h(b)(z')(\lo(Bx)\mt\ro'(z')).$$}
\dow
\end{prop}
\begin{col} $B\in MC^*(\G)$ and $C^*(h)B=h(B)$. {\em Indeed, from statement
b) of the previous proposition: $||\pi_h(B\om)||\leq ||\pi_h(\om)||$ so also
$||B \om||_{C^*} \leq ||\om||_{C^*}$, so $B$ can be extended to a bounded
linear mapping of $C^*(\G)$. Moreover from Lem. \ref{bis-gw} 
$B\in MC^*(\G)$ and $B^*=s(B)$. From the point a) of the 
previous proposition $C^*(h)(B)=h(B)$.}
\end{col}

\noindent {\bf Functions on $\G^0$ as affiliated elements.}\vs

Let $g$ be a smooth function on $\G^0$. Define the action of $g$ on 
bidensities by the formula $(g\om)(x):=g(\el(x))\om(x)$. This is clearly
a linear mapping. The following lemma is easy to prove:
\begin{lem} a) $\om_1^*(g \om_2)=(g^*\om_1)^*\om_2\,,\,\om_1,\om_2\in 
\sA(\G)$ and $g^*$ is a complex conjugation of $g$.\\
b) For $h:\G\rel\G'$ - a morphism of differential groupoids let $h(g):E'\lra C$
be defined by: $h(g)(a):=g(f_h(a))$, then\\ 
$\hat{h}(g\om_1)\om_2=h(g)\hat{h}(\om_1)\om_2\,,\,\,\om_1\in \sA(\G),\om_2\in
\sA(\G')$\\
$\pi_{h}(g\om_1)\psi=h(g)\pi_h(\om_1)\psi\,,\,\om_1\in\sA(\G),\psi\in L^2(\G')$
- smooth with compact support 
and $h(g)$ is viewed as operator on $L^2(\G')$ by multiplication.\\
If $g$ is bounded then $||h(g)||\leq \,sup|g|$.\\
\dowl
\end{lem}
From the lemma it follows that if $g$ is smooth and bounded then $g$ defines
multiplier of $C^*(\G)$. Moreover, using morphism $l:\G\lra \G\times \G$ one 
can see that $||g||_{C^*}=sup|g|$, so we get isometric *-homomorphism from
algebra of continous, bounded function on $\G^0$ to a multiplier algebra
$MC^*(\G)$.

Now we prove that continous functions on $\G^0$ are affiliated to $C^*(\G)$.
Start with the following:
\begin{lem}
Let $g$ be a continous, bounded function on $\G^0$ and assume that $f(a)\neq 0$
for $a\in \G^0$. Then $g C^*(\G)$ is dense in $C^*(\G)$.\vs\\
Proof: {\em Take any $\om\in \sA(\G)$ and let $g_0$ be a smooth function on
$\G^0$ with compact support such that $g_0|_{\el(supp\,\om)}=1$. 
Then $\frac{g_0}{g}$ is bounded and continous, so it defines a multiplier
of $C^*(\G)$. We have $g(\frac{g_0}{g}\om)=(g\frac{g_0}{g})\om=g_0\om=\om$.
Since $\sA(\G)$ is dense in $C^*(\G)$ this proves the assertion.}\\
\dowl
\end{lem}
Now let $g$ be a continous function on $\G^0$. Define 
$z_g:=\frac{g}{\sqrt{1+|g|^2}}$. Then $z_g$ is a multiplier of $C^*(\G)$, 
$||z_g||\leq 1$ and $(1-z_g^*z_g)^{1/2}C^*(\G)$ is dense in $C^*(\G)$. 
So defining $T_g:C^*(\G)\supset D(T_g)\lra C^*(\G)$ by: 
$$(x\in D(T_g)\,\,{\rm and\,}\,y=T_gx)\iff 
(\exists\,a\in C^*(\G)\,:\,x=(1-z_g^*z_g)^{1/2}a\,\,{\rm and\,}\,y=z_g a)$$ 
we get an  element affiliated with $C^*(\G)$. (see appendix E). 
If $g$ is smooth and $\om\in\sA(\G)$ then $T_g \om=g \om$.\vs\\
{\bf Final remarks.}\\
The category of differential groupoids introduced by S. Zakrzewski led us
to functorial construction of $C^*$-algebra of differential groupoid. It seems
to be natural and generalize several well known examples. Will it be useful ?
It is almost certain that any double Lie group ( so for example any Iwasawa
decomposition) leads to a quantum group on the $C^*$ algebra level. 
On the other hand, Iwasawa decompositions (or very similar to them) defines
Poisson Lie structures on semidirect products, among them is Poincare group, so
by the groupoid approach we can get corresponding quantum groups (if they 
exist) directely on the $C^*$ level. However Poisson Lie structures coming from
double Lie groups don't exhaust all possible Poisson structures and it is 
necessary  to investigate more general situations (see Appendix A).

We end this work with the list of open questions:\vs\\
1. In  investigations of foliations we frequently meet groupoids
which are not Hausdorff manifolds. Is it possible to generalize effectively
our approach to such situations ?\vs\\
2. Is there any connections between modular functions ( or more precisely its
cohomology class) and equality $\cred=C^*(\G)$ ?\vs\\
3. Are sections of Lie algebroid of a differential groupoid affiliated
with $\cred$ or $C^*(\G)$ ?\vs\\
4. If we agree that morphisms of differential groupoids are relations, this
leads us to a new notion of morphisms between Lie algebroids. Is this notion
more natural or more useful in investigations of Lie algebroids ?\vs\\

\section{Appendixes.}

{\bf Appendix A: Lie groupoids and C$^*$-algebras}

\noindent
{\em This is a part of a preliminary version of introduction left by
S. Zakrzewski. }\vs

In this work we construct a functor from the category of smooth
groupoids (with suitably defined morphisms) to
the category of C$^*$-algebras (with morphisms defined as in the
context of {\em locally compact noncommutative spaces}, cf for
instance \cite{slw}). The existence of such a functor was
expected, once it was constructed for finite groupoids 
and since it has become clear what is the C$^*$-algebra of a
smooth groupoid.

We recall the definition of a smooth groupoid (called also {\em
Lie groupoid} or {\em differential groupoid}) below.
Apart from the standard defintion, we recall our
definition of groupoids as `algebras in the category of binary
relations', given previously in \cite{SZ1,SZ2}, since it is
crucial for the correct understanding of morphisms of groupoids
in our sense.

In the rest of this introduction, let us explain the role
of our construction in establishing relations between `classical'
and `quantum' theories. Recall that symplectic manifolds
correspond to (play a similar role as) Hilbert spaces
(possibly projective) and symplectic diffeomorphisms correspond
to unitaires.

In order to have a procedure which relates some concrete
symplectic manifolds to some concrete Hilbert spaces
and some concrete symplectic diffeomorphisms to unitaries
one has to consider more special situation. Suppose we are
given a manifold $Q$ (playing the role of `configurations'). We
have then immediately the corresponding phase space $S=T^*Q$ and
also the Hilbert space $H=L^2(Q)$ (of square-integrable complex
half-densities on $Q$). To any diffeomorphism $\phi$ of $Q$
there corresponds a symplectomorphism $u:=\phi _*$ (the
push-forward of covectors) and also a unitary operator $U:=\phi
_*$ (the push-forward of half-densities). 
It is clear that in these cicumstances we have a 1--1
correspondence between such $u$'s and $U$'s,
illustrated by the following diagrams:\vs

\begin{figure}[tbhp]\label{diag1}
\p{s=t*q}{$S=T^*Q$}
\p{h=l2(q)}{$H=L^2(Q)$}
\p{q}{$Q$}
\p{u=f*}{$u=\phi_*$}
\p{fdq}{$\phi\in Diff(Q)$}
\p{U=f*}{$U=\phi_*$}
\centering
\includegraphics[height=0.15\textheight,width=0.9\textwidth]
{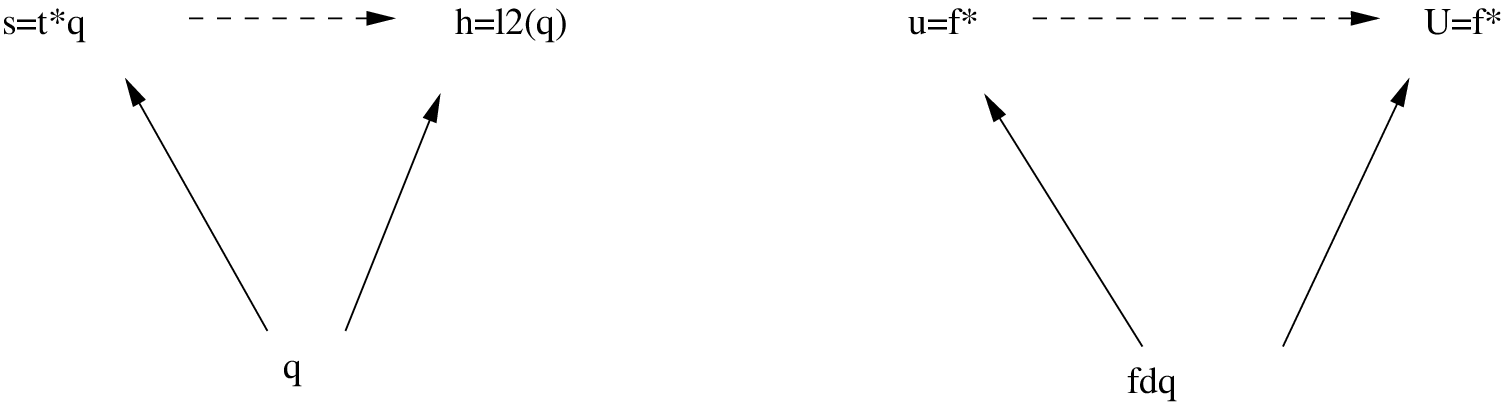}
\caption{}
\end{figure}
We see that the transition from the classical level to the
quantum level is possible in this case due to the common
`configuration level':

\begin{figure}[tbhp]\label{diag2}
\p{classical}{$\begin{array}{c}\mbox{classical}\\\mbox{level}\end{array}$}
\p{quantum}{$\begin{array}{c}\mbox{quantum}\\\mbox{level}\end{array}$}
\p{configuration}{$\begin{array}{c}\\ \mbox{configuration}\\\mbox{level}
\end{array}$}
\p{u=f*}{$u=\phi_*$}
\p{fdq}{$\phi\in Diff(Q)$}
\p{U=f*}{$U=\phi_*$}
\centering
\includegraphics[height=0.15\textheight,width=0.3\textwidth]
{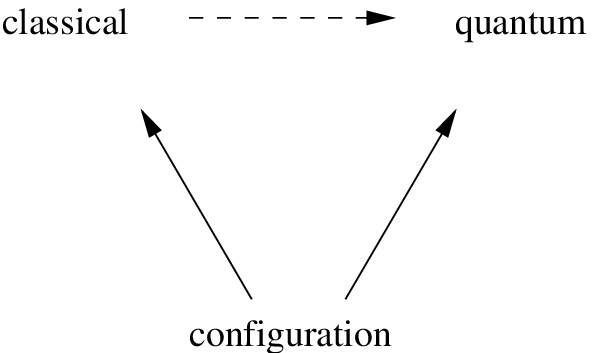}
\caption{}
\end{figure}
The symplectic diffeomorphisms of $T^*Q$
being just the natural lifts of diffeomorphisms
of $Q$ are said to be {\em point transformations}. It turns out
that not only these can be `quantized'. Namely, as a second step
consider {\em phase shifts} of $T^*Q$, that is symplectic
diffeomorphisms $v$ of $T^*Q$ of the form
\[ 
T^*Q\ni\xi\mapsto \xi + df (\pi (\xi ))\in T^*Q,
\]
where $f$ is a smooth function on $Q$ and $\pi\colon T^*Q\to Q$
is the cotangent bundle projection. It is natural to associate
with $f$ also the unitary operator $V$ in $L^2(Q)$ of
multiplication by $e^{if}$. Symbolically, we have

\begin{figure}[tbhp]\label{diag3}
\p{v=+df}{$v=+df$}
\p{V=.eif}{$V=\cdot e^{i f}$}
\p{f:qstrzr}{$f:Q\lra R$}
\p{u=f*}{$u=\phi_*$}
\p{configuration}{{\small configuration}}
\p{U=f*}{$U=\phi_*$}
\centering
\includegraphics[height=0.15\textheight,width=0.3\textwidth]
{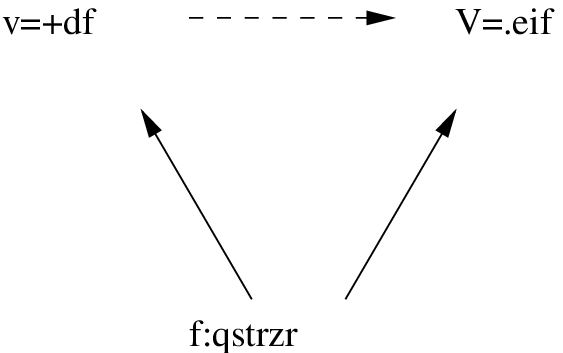}
\caption{}
\end{figure}
It means that we can associate (projective) unitaries with phase shifts. 
Moreover, symplectic diffeomorphisms of the form $vu$ form a group, which 
can be then naturally mapped into the unitary opearators by the rule
\[
vu\mapsto VU
\]
(it works modulo the phase factor).
What we here obtain is (essentially) the quantization of
symplectic diffeomorphisms which preserve the natural
polarization of the cotangent bundle (i.e. map fibers onto
fibers). In fact, to construct (projective) $H$ from $S$, the
polarization is sufficient (the change of the lagrangian
section playig the role of the `zero section' is then
implemented by the corresponding unitary transformations of type
$V$).

We may summarize the above discussion as follows. A symplectic
manifold $S$ may serve to construct a quantum-mechanical
Hilbert space $H$ if $S$ comes from configurations, $S=T^*Q$, or if
at least $S$ is equipped with projection on $Q$ with lagrangian
fibers (which essentially means that $S$ is equipped with a
polarization). Then a symplectic diffeomorphism of $S$ may serve
to construct a unitary operator in $H$ if it is a point
transformation, or, at least if it preserves the polarization.

Now the point is that such important for quantum mechanics
structures as operator algebras (in particular, C$^*$-algebras)
have also classical counterparts, namely symplectic groupoids.
In this context, above diagrams have the following form:

\begin{figure}[tbhp]\label{diag4}
\p{classical}{\hspace{-2em}$\begin{array}{c}\mbox{symplectic}\\
\mbox{groupoids} \end{array}$}
\p{quantum}{\hspace{-1.5em}$\begin{array}{c}\mbox{$C^*$-algebras}
\end{array}$}
\p{configuration}{$\begin{array}{c} \mbox{Lie groupoids}
\end{array}$}
\p{u=f*}{$u=\phi_*$}
\p{fdq}{$\phi\in Diff(Q)$}
\p{U=f*}{$U=\phi_*$}
\centering
\includegraphics[height=0.15\textheight,width=0.3\textwidth]
{figure7.eps}
\caption{}
\end{figure}
\FloatBarrier
\noindent
whose concrete realization is

\begin{figure}[tbhp]\label{diag5}
\p{v=+df}{\hspace{1em}$T^*\G$}
\p{V=.eif}{$C^*(\G)$}
\p{f:qstrzr}{$\hspace{2.5em}\G$}
\p{u=f*}{$u=\phi_*$}
\p{configuration}{{\small configuration}}
\p{U=f*}{$U=\phi_*$}
\centering
\includegraphics[height=0.15\textheight,width=0.3\textwidth]
{figure8.eps}
\caption{}
\end{figure}

\noindent
Here $\G$ is a Lie groupoid, $T^*\G$ is its cotangent (symplectic)
groupoid and $C^*(\G)$ is the C$^*$-algebra of $\G$. Similarly, for
morphisms, we shall have (as a result of the present paper)

\begin{figure}[tbhp]\label{diag6}
\p{classical}{{\small $\hspace{-9em}Ph(h)\in Mor(T^*(\G),T^*(\G'))$}}
\p{quantum}{{\small $C^*(h)\in Mor(C^*(\G),C^*(\G'))$}}
\p{configuration}{{\small $h\in Mor(\G,\G')$}}
\p{u=f*}{$u=\phi_*$}
\p{fdq}{$\phi\in Diff(Q)$}
\p{U=f*}{$U=\phi_*$}
\centering
\includegraphics[height=0.15\textheight,width=0.3\textwidth]
{figure7.eps}
\caption{}
\end{figure}

This corresponds to the `point case'. There is also the second
step, admitting also `phase shifts'. It consists in considering
symplectic groupoids, which are projectable on Lie groupoids (in
the sense that $\G$ is a cotangent bundle of some manifold $Q$
and the multiplication relation projects onto a Lie groupoid
multiplication relation on $Q$). It turns out that the
symplectic groupoid structure of $\G$ is the cotangent lift of
the Lie groupoid structure, shifted by a `2-cocycle' and the
previous diagram is generalized to 

\begin{figure}[tbhp]\label{diag7}
\p{classical}{\hspace{-10em}$\begin{array}{r}
\mbox{symplectic groupoids }\\ \mbox{projectable on Lie groupoids}
\end{array}$}
\p{quantum}{$C^*$-algebras}
\p{configuration}{\hspace{-1em}$\begin{array}{r} \\ \mbox{Lie groupoids}\\
\mbox{with 2-cocycle}\end{array}$}
\p{u=f*}{$u=\phi_*$}
\p{fdq}{$\phi\in Diff(Q)$}
\p{U=f*}{$U=\phi_*$}
\centering
\includegraphics[height=0.15\textheight,width=0.3\textwidth]
{figure7.eps}
\caption{}
\end{figure}

\noindent
The situation with 2-cocycles will be described in another
paper.\FloatBarrier

\vs\vs

\noindent
{\bf Appendix B: Cocycles and one parameters group on \cred}.
\begin{defi}
A (smooth) one cocycle on $\G$ {\em is a smooth function $\sigma :\G\lra C$ 
which satisfies condition:  $(x,y)\in\G^2\Rightarrow f(x y)=f(x) f(y)$.
}\end{defi}
For $x\in \G^0$ we have $\sigma(x)=\sigma(x x)=\sigma(x)^2$ so 
$\sigma |_{\G^0}$ is either 1 or 0. We assume that $\sigma |_{\G^0}=1$, then
for any $x\in G\,,\, \sigma(x)\neq 0$. From now on the word {\em cocycle}
will mean smooth, nonvanishing one cocycle.
\begin{ex} {\em a) If $f$ is a non vanishing, smooth function on $E$ then 
$\sigma(x):=\frac{f(\el(x))}{f(\er(x))}$ is a one cocycle on $\G$.\\
b) $\D$ given in Example \ref{red-mod}b) is a one cocycle on $\G$.\\
c) Let $(G;A,B)$ be a double Lie group. The decomposition $G=AB$ defines
decomposition $\g:=\ab\ms\bb$ of a Lie algebra of $G$. Then operators of
the adjoint representation of $G$ can be written as 
$Ad_g:=\left(\begin{array}{cc} \alpha_1 &\alpha_2\\\alpha_3 &\alpha_4
\end{array}\right)$ for $\alpha_1:\ab\lra\bb\,,\,\alpha_2:\bb\lra\ab$, etc. 
Then the function $Q(g):=|\frac{\det Ad_g}{\det \alpha_1 \det \alpha_4}|$ is 
a one cocycle on $G_A$ and $G_B$. This function plays an importent role in 
a definition of the quantum group structure on $C^*(G_A)$.}
\end{ex}  
For a cocycle $\sigma$ we define the mapping $\sigma:\sA\ni\om\mapsto
\sigma(\om):=\sigma \om\in \sA$. Note the following:
\begin{lem}
{\em 1.  $\sigma$ is a homomorphism of $\sA$\\
2.  If $|\sigma|=1$ then $\sigma$ is a *-homomorphism of $\sA$.}\vs\\
Proof: {\em Using the formulas given in Col. \ref{algebra} one performs
simple computations.}\\
\dowl
\end{lem}
In this way for any cocycle $\sigma:\G\lra ]0,\infty[$ 
the mapping  $R\ni t\mapsto \sigma^{i t}$ defines one parameter group of
*-automorphisms of $\sA$.
 Also it defines strongly continous one parameter group
$U_{\sigma}(t)$ of unitary operators on $L^2(\G):
(U_{\sigma}(t)\psi)(x):=\sigma^{i t}(x)\psi(x)$. Moreover
$\pi_{id}(\sigma_t\om)=U_{\sigma}(t)\pi_{id}(\om)U_{\sigma}(-t)\,$ for
$\om\in\sA(\G)$. 
\begin{prop}
For any smooth cocycle $\sigma: \G\lra ]0,\infty[\,$, $\sigma_t$ is strongly
continous one parameter group  on \cred. (see appendix E)\vs\\
Proof: {\em First note that if $U(t)$ is a one parameter group of unitaries
on a Hilbert space $H$ then the set $C:=\{a\in B(H)\,:\,\lim_{t\rightarrow
0}||U(t)aU(-t)-a||=0\}$ is closed * - subalgebra of $B(H)$. We will show that
$\lim_{t\rightarrow 0}||U_{\sigma}(t) \pi_{id}(\om)
U_{\sigma}(-t)-\pi_{id}(\om)||=0$ for $\om\in \sA(\G)$.  \\
Now take some $\om_0$ and let $||.||_l,\,||.||_r,\,||.||_0$ denote the
associated norms on $\sA(\G)$. From prop. \ref{rep} we have
$||\pi_{id}(\om)||\leq||\om||_0$ and $$||U_{\sigma}(t) \pi_{id}\om)
U_{\sigma}(-t)-\pi_{id}(\om)||=||\pi_{id}(\sigma_t(\om)-\om)||\leq
||\sigma_t(\om)-\om||_0.$$ 
$||\sigma_t(\om)-\om||_l=sup_{a\in\G^0}\int_{F_l(a)}\lo(x)^2 |\sigma(x)^{i t}
f(x)-f(x)|$ where $\om=f\om_0$. Since supp$\,\om$ is compact, one can find
$M$ such that $|\log \sigma(x)|\leq M$ for $x\in$ supp$\,\om$. So for $x\in$
 supp$\om$ :  $|\sigma(x)^{i t} f(x)-f(x)|\leq |t| M|f(x)|$
and $||\sigma_t(\om)-\om||_l\leq  |t| M||\om||_l$. In the same way we
have $||\sigma_t(\om)-\om||_r\leq |t| M||\om||_r$ and as a consequence:
$\lim_{t\rightarrow 0}||U_{\sigma}(t) \pi_{id}(\om)
U_{\sigma}(-t)-\pi_{id}(\om)||=0$ for $\om\in \sA(\G)$. Now the assertion
follows from the remark on the beginning of the proof.}\\
\dow
\end{prop}
Let $\sigma_i$ be an analytic generator of $\sigma_t$ (see appendix E).
\begin{prop}\label{core} 
$\sA(\G)$ is a core for $\sigma_i$ and for
$\om\in\sA(\G)\,\,\sigma_i(\om)=\sigma^{-1}\om $.\vs\\
Proof: {\em Start with the following:
\begin{lem}
Let $\sigma:\G\lra ]0,\infty[$ be a smooth cocycle and $\om\in\sA(\G)$. The
function $C\ni z\mapsto \sigma^{i z} \om\in C^*(\G)$ (or \cred) is an
entire analytic. \vs\\
Proof: {\em Choose some $\om_0$, then straightforward computation shows that
$\lim_{h\rightarrow 0}
||\frac{\sigma^{i(z+h)}\om-\sigma^{i z}\om}{h}-i \log(\sigma)
\sigma^{i z}\om||_0=0.$ 
Since $||\om||_{C^*_{red}}\leq||\om||_{C^*}\leq||\om||_0$ we have the desired
result.}\\
\dowl
\end{lem}
From the lemma we have that $\sA(\G)\subset D(\sigma_i)$ and
$\sigma_i(\om)=\sigma^{-1}\om$. For $a\in B(L^2(\G))$ we set:
$$R_n(a):=\frac{n}{\sqrt{\pi}} \int_R d t e^{-n^2 t^2} U_{\sigma}(t) a
U_{\sigma}(-t)$$ 
(since $U_{\sigma}(t)$ is strongly continous $U_{\sigma}(t) a
U_{\sigma}(-t)$ is $\sigma$ - weakly continous, and the integral is defined
by the property that for any normal functional $\varphi$ on $B(L^2(\G))$:
$\varphi(R_n(a))=\frac{n}{\sqrt{\pi}} \int_R d t e^{-n^2 t^2}
\varphi(U_{\sigma}(t) a U_{\sigma}(-t))$). We know \cite{slw} that the set
$\{R_n(a)\,:\,n\in N,a=\pi_l(\om),\om\in\sA(\G)\}$ is a core for $\sigma_i$.
So it is enaugh to show that $R_n(\pi_l(\sA))\subset \pi_l(\sA)$. Let
$\sigma_n(x):=\frac{n}{\sqrt{\pi}} \int_R d t e^{-n^2 t^2} \sigma(x)^{i t}$.
It is clear that $\sigma_n \sA\subset \sA$. We claim that
$R_n(\pi_l(\om))=\pi_l(\sigma_n\om)$. Indeed, take any $\psi$ - smooth, half
density with compact support on $\G$, then the equality
$(\psi,R_n(\pi_l(\om))\psi)=(\psi,\pi_l(\sigma_n\om)\psi)$ follows from
Fubini theorem.}\\
\dow
\end{prop}

\vs\vs

\noindent
{\bf Appendix C: Weights on \cred}\vs

Choose some $\ro$ and real, smooth half density $\nu$ on $\G^0$.
Such  data define linear functional
$\varphi$ on $\sA$ as follows: write $\om=f \om_0$ and put 
$\varphi(\om):=\int_{\G^0} \nu^2\, f.$
Also they define linear mapping  $\hat{\varphi}:\sA(\G)\ni\om\mapsto 
f \ro\mt \nu \in L^2(\G)$. 
 Note the following:
\begin{lem} a) $\varphi(\om^*\om)=(\hat{\varphi}(\om),\hat{\varphi}(\om))$. 
 (so $\varphi$ is a positive, linear functional on $\sA(\G)\,$)\\
b) If $\pi_{id}$ is the left regular representation of $\sA(\G)$ on $L^2(\G)$ 
then: $\hat{\varphi}(\om_1\, \om_2)=\pi_{id}(\om_1)\hat{\varphi}(\om_2)$\vs\\
Proof: {\em Straightforward computation.}\\
\dowl
\end{lem}
Now we assume that $\nu$ is non vanishing. 
Due to the above proposition we can identify the Hilbert space given by 
the GNS construction for $\varphi$ and associated representation of $\sA$ 
with $L^2(\G)$ and left regular representation. Next we show that $\varphi$ c
an be extended to a weight on $C^*(\G)$. We recall the following 
theorem \cite{W-M-N}:
\begin{tw}\label{waga} Let $\pi$ be a nondegenerate representation 
of a separable
$C^*$ algebra $A$ on a separable Hilbert space $H$. Morover let 
$\hat{\varphi} :A\supset D(\hat{\varphi})\lra H$ be a closed, densely defined
linear mapping with the dense range such that $D(\hat{\varphi})$ is a left 
ideal of $A$  and $\hat{\varphi}(a b)=\pi(a)\hat{\varphi}(b)\,,\,a\in A\,,\,
b\in D(\hat{\varphi})$. Then the formula $\varphi(a^* a):=
\left\{ \begin{array}{ll} (\hat{\varphi}(a),\hat{\varphi}(a)) & if\,\,a\in
D(\hat{\varphi})\\ \infty & otherwise \end{array} \right.$
defines a locally finite, lower semicontinous weight on $A$.\\
\dow
\end{tw}
First we are going  to prove the following:
\begin{prop} $\hat{\varphi}: C^*_{red}(\G)\lra L^2(\G)$ {\em is closable.}\vs\\
Proof: {\em We have to show, that if
$\lim_{n\rightarrow\infty}\pi_l(\om_n)=0$ - in $BL^2(\G)$ and
$\lim_{n\rightarrow\infty}\hat{\varphi}(\om_n)=\psi$ - in $L^2(\G)$ then
$\psi=0$.\\  
The groupoid inverse $s$ is a diffeomorphism of $\G$ so it defines the unitary
operator $S:\,L^2(\G)\lra L^2(\G)$, also it defines linear,
antimultiplicative bijection, which we denote also by $s: \sA\lra\sA$. Let
$\D$ be a modular function for $(\ro,\nu)$, see Ex. \ref{red-mod} b) i.e.
 $\D$ is given by: $\D(x)(\ro\mt\nu)(x)=(\lo\mt\nu)(x)$.
\begin{lem} $S\, \hat{\varphi}(\om)
=\hat{\varphi}\, \D\, s\,(\om)$ $\om\in\sA$. \vs\\
Proof: {\em First let us note that $S(\ro\mt\nu)=\lo\mt\nu$. Indeed, let
$\{X^1,...,X^m\}$ be a basis in $T_x\G$ such that $\{X_1,...,X^k\}$ is a
basis in $T_x^l\G$ (then of course $\{s(X^1),..., s(X^k)\}$ is a basis in
$T_{s(x)}^r\G$ ). Compute: 
$$(S(\ro\mt\nu))(x)(X^1\wedge ...\wedge
X^m):=(\ro\mt\nu)(s(x))(s(X^1)\wedge ...\wedge s(X^m))=$$
$$=\ro(s(x))(s(X^1)\wedge
...\wedge s(X^k))\,\nu(\er(s(x)))(\er s(X^{k+1})\wedge ...\wedge \er s(X^m))=$$
$$=\lo(x)(X^1\wedge ...\wedge X^k)\,\nu(\el(x))(\el X^{k+1}\wedge ...\wedge \el
X^m)=(\lo\mt\nu)(x)(X^1\wedge ...\wedge X^m),$$ since $\ro s=\lo$ and $\er
s=\el$.\\
Now let $\om_0$  $\om=f\om_0$. Then $(S\hat{\varphi}(\om))(x)=
(S(f\,\ro\mt\nu))(x)=f(s(x))(\lo\mt\nu)(x)$ and 
$\D\, s\,(f\om_0)=\tilde{f}\om_0$ for $\tilde{f}(x):=\D(x) f(s(x))$ so
$$(\hat{\varphi}\D\,s\,(\om))(x)=(\hat{\varphi}(\tilde{f}\om_0))(x)=
\tilde{f}(x)(\ro\mt\nu)(x)=\D(x)
f(s(x))(\ro\mt\nu)(x)=f(s(x))(\lo\mt\nu)(x).$$} 
\dowl
\end{lem}
For any $\om\in\sA$:
$$S\,\pi_{id}(\om)\,S\,\psi=S\,\pi_{id}(\om)\,S\,(\lim_{n\rightarrow\infty}
\hat{\varphi}(\om_n)))= S(\lim_{n\rightarrow\infty}\pi_{id}(\om)
S\hat{\varphi}(\om_n))=$$
$$=S(\lim_{n\rightarrow\infty}\pi_{id}(\om)\hat{\varphi}(\D s(\om_n)))=  
= S(\lim_{n\rightarrow\infty}\hat{\varphi}(\om \D  s(\om_n)))=$$
$$=S(\lim_{n\rightarrow\infty}\hat{\varphi}(\D s (\om_n s \D^{-1}(\om))))=
S(\lim_{n\rightarrow\infty}S\hat{\varphi}(\om_n s \D^{-1}(\om)))=
\lim_{n\rightarrow\infty}\pi_{id}(\om_n)\hat{\varphi}(s \D^{-1}(\om))=0.$$
So $\pi_{id}(\om) \,S\,\psi=0$, for any $\om\in\sA(\G)$. 
Since left regular representation is
nondegenerate $S(\psi)=0$ and $\psi=0$.}\\
\dow
\end{prop}
We denote the closure of $\hat{\varphi}$ by the same symbol. So
$\hat{\varphi}: C^*_{red}(\G)\supset D(\hat{\varphi})\lra L^2(\G)$ where
$D(\hat{\varphi}:=\{a\in C^*_{red}(\G)\,:\,a=\lim
a_n\,,\,a_n\in\sA(\G)\,,\,$and there exists $ x\in L^2(\G)\,,\,x=\lim
\hat{\varphi}(a_n)\}$. It is also clear that the range of $\hfi$ is dense in
$L^2(\G)$.    
\begin{lem}
$D(\hfi)$ is a left ideal in \cred and $\hfi(a b)=\pi_{id}(a)\hfi(b)$ for $a\in
\cred\,,\, b\in D(\hfi)$.\vs\\
Proof: {\em Let $a\in\cred\,,\,b\in D(\hfi)$, $b=\lim b_n$ and $a=\lim
a_m$ for some $a_m,b_n \in \sA(\G)$. Then $a b=\lim a_n b_n$, 
$\pi_{id}(a)=\lim \pi_{id}(a_n)$ and $x:=\lim \hfi(b_n)$. 
$\hfi(a_n b_n)=\pi_{id}(a_n) \hfi(b_n)$
and this sequence tends to $\pi_{id}(a) x$. 
It follows that $a b\in D(\hfi)$ and $\hfi(a b)=\pi_{id}(a)\hfi(b).$}\\
\dowl
\end{lem}
Using theorem \ref{waga} $\varphi$ extends to a locally finite, lower 
semicontinous weight on \cred.

\newcommand{\sip}{\sigma_{\frac{i}{2}}}

Now we prove that $\varphi$ is a KMS-weight with $\sigma_t:=|\D|^{-2 i t}$ 
as a modular group. Since $\D$ is non vanishing $|\D|$ is a smooth cocycle
and $\sigma_t$ is strongly continous one parameter group on \cred. Using
the same arguments as in Prop. \ref{core} one can show that $\sA(\G)$ is a core
for $\sip$ and for $\om\in \sA(\G)$: 
$\sip(\om)= \frac{1}{|\D|}\om$. Consider the following:
\begin{lem} For $a\in D(\sip)$: $a\in D(\hfi)\iff(\sip(a))^*\in D(\hfi)$.\vs\\
Proof: {\em For $\om=f \om_0\in\sA(\G)$ we have: 
$(\sip(\om))^*=|\D| f^*\om_0$ and $\hfi((\sip(\om))^*)=|\D|f^*\ro\mt\nu$.
$$|||\D|f^*\ro\mt\nu||^2=\int_\G (\ro^2(x)\mt\nu^2(x)) |\D(x)|^2 |f(s(x))|^2=
\int_\G(\lo^2(x)\mt\nu^2(x))|\D(s(x))|^2 |f(x)|^2=$$
$$=\int_\G (\ro^2(x)\mt\nu^2(x))|f(x)|^2=||\hfi(\om)||^2,$$ 
where in the second 
equality we use $S$ and in the third the definition of $\D$. So for 
$\om\in\sA(\G)$ we have: $||\hfi((\sip(\om))^*)||=||\hfi(\om)||$.\\
Now let $a\in D(\sip)\,,\, a=\lim \om_n\,,\,\om_n=f_n\om_0\in\sA(\G)$ and 
$\sip(a)=\lim\sip(\om_n)$. So $a\in D(\hfi)$ is equivalent to a convergence
of a sequence $\hfi(\om_n)$ and due to the above equality this is equivalent
to a convergence of a sequence $\hfi((\sip(\om_n)^*)$.}\\
\dowl
\end{lem}
From the above lemma we have: $\varphi(a^*a)=(\hfi(a),\hfi(a))=
(\hfi((\sip(a))^*),\hfi((\sip(a))^*))=\varphi(\sip(a) \sip(a)^*)\,,\,a\in 
D(\sip)$. So to prove that $\varphi$ is a KMS weight it remains to
show that $\varphi \sigma_t=\varphi$ and this is straightforward
computation. In this way we prove:
\begin{prop} Let $\ro$ be smooth, real, non vanishing, rightinvariant
half density on $\G$ and $\nu$ be real, smooth, non vanishing half density
on $\G^0$. Let $\D$ be a modular function for this pair. Then the formula:
$\varphi(f\om_0):=\int_{\G^0} f \nu^2$ defines locally finite, lower 
semicontinous weigth on \cred which is KMS weight with a modular group: 
$\sigma_t:=|\D|^{-2 i t}.$\\
\dow
\end{prop}

\noindent
{\bf Appendix D: Subgroupoids and homogenous spaces.}
\newcommand{\tE}{\tilde{E}}
\newcommand{\tG}{\mbox{$\tilde{\G}$}}
\newcommand{\tm}{\mbox{$\tilde{m}$}}
\newcommand{\ts}{\mbox{$\tilde{s}$}}

In this appendix we discuss the notion of subgroupoid and present
some constructions related to it. Loosely speaking subgroupoid is subset
closed with respect to multiplication and involution. We keep in mind
the basic examples: subgroup of a group, subset of a set, 
equivalence relation  and cartesian product $A\times A\,,\,A\subset X$ for
pair groupoid $X\times X$. On the level of  suitable algebras of 
functions in these examples we can relate to them the following constructions.
Having a subgroup $H\subset G$ we can form the homogoneus space $G/H$ 
together with an action of $G$. This action defines a morphism from $G$ to a 
pair groupoid $G/H\times G/H$ which in turns defines a representation of group
algebra of $G$. On the level of sets, if $A\subset X$ then restriction of 
a function on $X$ to $A$ defines a morphism from algebra of functions on $X$
to an algebra of functions on $A$. For an equivalence relation 
$R\subset X\times X$ we can perform construction simmilar to the subgroup case
and get a morphism from pair groupoid $X\times X$ to another pair groupoid.
For the last example $A\times A\subset X\times X$ we get just subalgebra of
 functions on $X\times X$ which is loosely related to the whole algebra.

After these preliminary remarks we define objects and present constructions.
According to the general line of this work we start from pure algebraic
situation and later on add differential structure.
\begin{defi} {\em Let $\tG \subset \G$. Define $\tm:\tG\times \tG\rel\tG$,
 $\te: \{1\}\rel \tG$ and $\ts$  by: 
$Gr(\tm):=Gr(m)\cap(\tG\times \tG\times \tG)$,
$Gr(\te):=Gr(e)\cap(\tG\times\{1\})$ and $\ts:=s|_{\tG}$. \tG} is 
 a subgroupoid of $\G$ 
{\em iff  $(\tG,\tm,\ts,\te)$ is a groupoid.}
\end{defi}

We denote:  $\tE:=E\cap \tG$. For any subset $A\subset E$ the set 
$\tG:=\el^{-1}(A)\cap\er^{-1}(A)$ is a subgroupoid of $\G$. In the following
we restrict our attention to subgroupoids of special kind which directly 
correspond to subgroups, subsets and equivalence relations. 
A subgroupoid $\tG\subset \G$ is {\em horizontal} iff $\tE=E$; it is 
{\em vertical } iff for any $x\in\tG$ the  fiber $F_l(x)$ is contained
in $\tG$ (then of course also the right fiber $F_r(x)$ is contained in $\tG$.)
\begin{lem} Let $\tG\subset\G$ be a subgroupoid and let $i:\tG\lra\G$ be the 
inclusion map. Then:\\
a) $\tG$ is horizontal iff $i:\tG\rel \G$ is a morphism.\\
b) $\tG$ is vertical iff $i^T:\G\rel\tG$ is a morphism.\\
Proof: {\em The proof is straightforward application of the definitions.}
\dowl
\end{lem}
Horizontal subgroupoids are also called {\em wide}.
\begin{ex}{\em a) If $\G$ is a group then every subgroup is horizontal
subgroupoid and the only one vertical subgroupoid is $\G$ itself.\\
b) If $\G$ is a set then every subset is vertical subgroupoid and the only
one horizontal subgroupoid is $\G$ itself.\\
c) If $\G:=X\times X$ is a pair groupoid, then every equivalence relation on 
$X$ is a horizontal subgroupoid and the only one vertical subgroupoid is $\G$
itself. (since $\G$ has only one orbit.)\\
d) If $\G$ is a bundle of groups over $E$, then each subbundle is a horizontal 
subgroupoid and restriction of $\G$ to a subset of $E$ is a vertical 
subgroupoid.}
\end{ex}
Let $\tG$ be a horizontal subgroupoid of $\G$. Consider the relation $R$ on 
$\G$: $(x, y)\in R\iff s(x) y\in \tG$.
It is easy to see that this is an equivalence relation. Let $\G/\tG=:Y$ denote 
the set of equivalence classes.
\begin{lem}
a) The mapping $R\ni (x,y)\mapsto x\in\G$ is surjective.\\
b) The mapping $f: Y\ni [x]\mapsto \el(x)\in\G^0$ is well defined and 
surjective.\\
c) Let $\G\times_f Y:=\{(x,y)\in \G\times Y\,:\,\er(x)=f(y)\}$. The mapping
$g: \G\times_f Y\ni(x,[x'])\mapsto [x x']\in Y$ is well defined and 
surjective.\\
d) The relation $h:\G\rel Y \times Y$ given by: 
$Gr(h):=\{(g(x,y),y;x)\,:\,(x,y)\in\G\times_f Y\}$ is a morphism from $\G$ to a
pair  groupoid $Y\times Y$.\\
Proof: {\em a) Let $x\in \G$. Since $\tG$ is horizontal, $\tE=E$ and there 
exists $z\in\tG$ with $\er(x)=\el(z).$ Then $(x,xz)\in R$.\\
b) From the definition of $R$ it is clear that if $(x ,y)\in R$ then 
$\el(x)=\el(y)$, so $f$ is well defined. Moreover for any $e\in E=\tE$ we 
have: $e=f([e]).$\\
c) It is clear that if $(x_1,x_2)\in R$ and $\er(x)=f([x_1])$ 
than $(x x_1,x x_2)\in R$ so the definition of $g$ is correct. Moreover
$g(f(y),y)=y$ for any $y\in Y$ so $g$ is surjection.\\
d) Let us check, that $h$ is a morphism of groupoids.\\
i) $hE=diag(Y\times Y)$. \\
If $e\in E$ and $f([x])=e$ than $[ex]=[x]$, 
so $hE\subset diag(Y\times Y)$. Also for any $[x]\in Y$ we have:
$([x],[x];f([x]))\in h$.\\
ii) $hs=s'h$. \\
$(y_1,y_2,x)\in hs\iff (y_1,y_2;s(x))\in h\iff \el(x)=f(y_2)
\,and \, y_1=s(x) y_2\iff \el(x)=f(y_2)\,and\,\er(x)=f(y_2)\,and\,
y_2=x y_1\iff (y_2,y_1;x)\in h\iff (y_1,y_2;x)\in s'h$.\\
iii) $hm=m'(h\times h)$. \\
$(y_1,y_2;x_1,x_2)\in hm\iff \er(x_1)=\el(x_2)\,and\,
(y_1,y_2;x_1 x_2)\in h\iff \er(x_1)=\el(x_2)\,and\,\er(x_2)=f(y_2)\,and\,
y_1=x_1 x_2 z_2\Rightarrow \er(x_1)=\el(x_2)\,and\,(x_2 y_2,y_2;x_2)\in h\,
and\,(y_1,x_2 y_2;x_1)\in h\Rightarrow (y_1,y_2;x_1,x_2)\in m'(h\times h)$.\\ 
Conversely, for $(y_1,y_2;x_1,x_2)\in m'(h\times h)\iff (y_1,y_3,x_1)\in h\,
and\,(y_3,y_2;x_2)\in h$ for some $y_3\in Y$, so $\er(x_1)=f(y_3)$,
$\er(x_2)=f(y_2)$ and $y_3=x_2 y_2$. Now $f(y_3)=\el(x_2)$ so $x_1,x_2$
are composable and $(y_1,y_2;x_1 x_2)\in h$.\\}
\dowl
\end{lem}
In the differential setting we define:
\begin{defi} {\em Let $\G$ be a differential groupoid and $\tG\subset\G$ be
a submanifold. $\tG$ is } a differential subgroupoid of $\G$ {\em iff 
$(\tG,\tm,\ts,\te)$ is a differential groupoid.}
\end{defi}
Let $\tG\subset\G$ be a horizontal subgroupoid of $\G$ and let the relation
$R$, the mappings $f$ and $g$ be as above.
We have the following:
\begin{lem} a) $R$ is a submanifold of $\G\times\G$\\
b) The mapping: $R\ni (x,y)\mapsto x\in \G$ is a surjective submersion.\\
c) If $\tG$ is closed then from a) and b) we have that $Y$ is a manifold 
 and then  the mapping $f$ is a surjective submersion.\\
d) $\G\times_f Y$ is a submanifold of $\G\times Y$ and $g$ is a surjective 
submersion.\\
e) The relation $h$ defined in the previous lemma is a morphism of 
differential groupoids.\vs\\
Proof: {\em a) Consider $\{(x,y)\in \Gd\,;\,m(x,y)\in \tG\}$. Since $m$ 
restricted  to $\Gd$ is surjective submersion and $\tG$ is a submanifold, this
set is a submanifold in $\Gd$ so also  in $\G \times \G$. $R$ is the image 
of this submanifold by the diffeomorphism: $(x,y)\mapsto (s(x),y)$.\\
b) It is clear that this mapping is smooth and surjective. Let 
$(x_0,y_0)\in R$ and $s(x_0) y_0=z_0\in \tG$. Let $x(t)$ be a curve through 
$x_0$, then $\er(x(t))$ is a curve in $\G^0$ through $\er(x_0)=\el(z_0)$. 
Since $\te_L$ is a surjective submersion one
can lift it to a curve $z(t)\in\tG$ through $z_0$. Then $(x(t),x(t) z(t))$ 
is a curve in $R$ through $(x_0,y_0)$.\\
c) This is clear, since $\el=f \pi$ where $\pi: \G\lra Y$ is a 
surjective submersion.\\
d) Since $\G\times_f Y=(\er\times f)^{-1}(diag(\G^0\times\G^0)$ and 
$\er\times f$ is a surjective submersion it is clear that $\G\times_f Y$ is a 
submanifold in $\G\times Y$. On $\Gd$ consider the relation 
$\tilde{R}:=id\times R$ i.e. 
$(x_1,x_2;x_3,x_4)\in \tilde{R}\iff x_1=x_3\,and\,(x_2,x_4)\in R$.
 Then this is a regular equivalence relation and $\Gd/\tilde{R}=\G\times_f Y$.
Let $\tilde{pi}:\Gd\lra\Gd/\tilde{R}$ be the cannonical projection. 
Then $g$ is 
determined by: $\underline{m} \pi=g \tilde{pi}$ where 
$\underline{m}:\Gd\lra\G$ is a restriction of $m$ to $\Gd$. Since this 
restriction is  surjective submersion $g$ is a surjective submersion.\\
e) From the previous lemma we know that $h$ is a morphism of groupoids. We 
have to check that $Gr(h)$ is a submanifold and two transversality condition 
$m'\tran (h\times h)$ and $h\tran e$ hold. Let 
$\tilde{Y}:=\{(y,x)\in Y\times X\,:\,(x,y)\in\G\times_f Y\}$. Then $\tilde{Y}$
is a submanifold in $Y\times\G$ and $\tilde{g}:\tilde{Y}\ni(y,x)\mapsto
g(x,y)\in Y$ is a surjective submersion. So $Y\times\tilde{Y}$ is a 
submanifold in $Y\times Y\times \G$ and $Gr(h)=(\pi_1\times \tilde{g})^{-1}
(diag(Y\times Y)$ where $\pi_1:Y\times\tilde{Y}\ni(y,\tilde{y})\mapsto y\in Y$.
Since $\pi_1\times \tilde{g}$ is a surjective submersion $Gr(h)$ is a 
submanifold in $Y\times\tilde{Y}$ and so in $Y\times Y\times \G$.

Now we check transversality condition. First we take $m'\tran (h\times h)$.
Since $m'(h\times h)=h m$ and $m$ is a differentiable reduction we know
\cite{SZ2} that $Gr(hm)=Gr(m'(h\times h))$ is a submanifold so it is enaugh to 
show that $Tm'$ and $T(h\times h)$ and $Pm'$ and $P(h\times h)$ have simple 
composition.
Let $(v;u,w)\in Tm'T(h\times h)$ be vectors tangent to $Gr(m'(h\times h)$ at 
points $(y_1,y_2;x_1,x_2)$. We can write $v=(v_1,v_2)$ for 
$v_1\in T_{y_1}Y\,,\,
v_2\in T_{y_2}Y$. Then we have $(v_1,v_3;u)\in Th$ and 
$(v_3,v_2;w)\in Th$ for some $v_3\in T_{y_3}Y$ where $y_3$ is a unique point
in $Y$ such that $(y_1,y_3;x_1)\in h$ and $(y_3,y_2;x_2)\in h$. But then
$v_3=g(w,v_2)$ and since $g$ is a mapping $v_3$ is uniquely determined.
So $Tm'$ and $T(h\times h)$ have a simple composition.

Now let $(\varphi_1,\varphi_2;\psi_1,\psi_2)\in Pm' P(h\times h)$ where
$\varphi_1\in T_{y_1}^*Y,\varphi_2\in T_{y_2}^*Y,\psi_1\in T_{x_1}^*\G,
\psi_2\in T_{x_2}^*\G$ and $(y_1,y_2;x_1,x_2)\in m'(h\times h)$. It is easy
to see that in the case of pair groupoids 
$$Gr(Pm):=\{(\varphi_1,\varphi_2;\varphi_1,\varphi_3,-\varphi_3,\varphi_2)\,:\,
\varphi_1\in T_{y_1}^*Y,\varphi_2\in T_{y_2}^*Y,\varphi_3\in T_{y_3}^*Y\}.$$
So our condition on $\varphi_1,\varphi_2,\psi_1,\psi_2$ is equivalent to
existence of some $\varphi_3\in T_{y_3}^*Y$ for unique $y_3$ - as above which
satisfies $(\varphi_1,\varphi_3;\psi_1)\in Ph$ and 
$(-\varphi_3,\varphi_2;\psi_2)\in Ph$. From the second condition: 
$-<\varphi_3,g(u_2,w_2)>=-<\varphi_2,w_2> +<\psi_2, u_2>$ for any
 $(u_2,w_2)\in T_{x_2,y_2}(\G\times_f Y)$. Since $g$ is a submersion this 
equality determines $\varphi_3$. In this way $Pm'$ and $P(h\times h)$ have 
simple composition. 

Now we check that $h\tran E$. The submanifold property is clear so we have to
verify only simple compositions conditions. If $(v_1,v_2;w)\in T(he)$ then
$v_1=v_2$ and $(w,v_2)\in T_{(x,y)}(\G\times_f Y)$ for some $w\in T_x\G^0$. 
But it means that $w=\er(w)=f(v_2)$ so $w$ is uniquely determined. If 
$(\varphi_1,\varphi_2;\psi)\in P(he)$ then 
$\psi\in T_x^*E\subset T_x^*\G$, $x=f(y)$,   
$<\varphi_1,g(w,v_2)>+<\varphi_2,v_2>=<\psi,w>$ for any $w=f(v_2)$ and 
$v_2\in T_{y_2}Y$. But $g(f(v_2),v_2)=v_2$ so the value of $\psi$ on the 
vectors $f(v_2)$ is determined and since $f$ is a submersion onto $E$ this
determines $\psi$ uniquely.
\dowl}
\end{lem}

\noindent
{\bf Appendix E: $C^*$-algebras.}\vs\\
{\em Multiplier algebra}. Let $A$ be a $C^*$-algebra. By $B(A)$ we denote the
algebra of bounded linear mappings acting on $A$. Let $a,b\in B(A)$, we say
that {\em $b$ is a hermitian adjoint of $a$} and write $b=a^*$ if:\\
$y^*(a x)=(b y)^*x$ for all $x,y\in A$. It follows that the set of those
$a\in B(A)$ which have a hermitian conjugate is a $C^*$-algebra. This is a
{\em multiplier algebra of $A$} and will be denoted by $M(A)$. $A$ can be
embedded into $M(A)$ via the left multiplication and the image is ideal in
$M(A)$. $M(A)=A$ iff $1\in A$.\vs\\
{\em $C^*$ - category.} Let $A,B$ be $C^*$ algebras. {\em
Morphism from $A$ to $B$} is 
*-homomorphism $\phi: A\lra M(B)$ such that the set $\phi(A) B$ is dense in
$B$.  That morphism can be composed follows from the fact that any such
$\phi$ extends uniquely to a $C^*$-homomorphism from $M(A)$ to $M(B)$, this
extension is defined by: $\hat{\phi}(m)(\phi(a) b):=\phi(m a) b$ for $m\in
M(A)\,,\,a\in A\,,\,b\in B$. If $\phi_1\in Mor(A,B)\,,\,\phi_2\in Mor(B,C)$
then composition is defined by $\hat{\phi_2}\phi_1: A\lra M(C)$. 
$C^*$ algebras with above defined morphisms form a {\em $C^*$ category}.\vs\\
{\em Affiliated elements.} Let $A$ be a $C^*$-algebra and 
$T: A\supset D(T)\lra A$ 
densely defined linear mapping. {\em $T$ is affiliated with $A$} iff there 
exists $z\in M(A)\,,\,||z||\leq 1$ such that 
$$(x\in D(T)\,\,{\rm and} \,\,y=Tx)\iff
(\exists\,a\in A : x=(I-z^* z)^{1/2} a\,\,{\rm and} \,\,y=z a).$$ 
If $1\in A$ then the set of affiliated elements is equal to $A$.\vs\\
{\em Weights on $C^*$ algebras.} Let $A$ be a $C^*$ algebra and $A_+$ be a
set of positive elements. {A weight $\varphi$ on $A$} is a mapping $\varphi:
A_+\lra [0,\infty]$ which satisfies:
$\varphi(a+b)=\varphi(a)+\varphi(b)\,,\,\varphi(\lambda a)=\lambda
\varphi(a)$ for any $\lambda\geq 0$ and $a,b\in A_+$. ($0\cdot\infty=0$).
The weight $\varphi$ is {\em densely defined} if the set $\{a\in
A_+\,:\,\varphi(a)<\infty\}$ is dense in $A_+$. $\varphi$ is {\em lower
semi-continous} if for any $t\in R_+$ the set $\{a\in A_+\,:\,\varphi(a)\leq
t\}$ is closed. \vs\\
{\em One parameter groups on $C^*$ algebras.} The homomorphism $\sigma:
R\lra Aut(A)$ - the group of * - automorphisms of $A$ such that
$||\sigma_t||\leq 1$ for any $t\in R$ is called {\em one
parameter group on $A$.} $\sigma$ is {\em strongly continous} iff for any
$a\in A$ the function: $R\ni t \mapsto \sigma_t(a)\in A$ is continous.\vs\\
{\em An analytic generator} of  a strongly continous one
parameter group $\sigma_t$ on $A$ is a linear mapping $\sigma_{i}:
D(\sigma_{i})\lra A$ defined in the following way: $a\in
D(\sigma_{i})$  iff there exists
continous function $f:\{z\in C\,:\,Im(z)\in[0,1]\}\lra A$, analytic
in the interior and such that $f(t)=\sigma_t(a)$ for all $t\in R$; then
$\sigma_{i}(a):=f(i)$. It can be proved that analytic
generators are densely defined, closed, multiplicative mappings.\vs \\
{\em KMS weights}. Let $\sigma_t$ be a strongly continous, one parameter
group on $A$, and $\varphi$ densely defined, lower semi-continous weight.
$\varphi$ is called {\em KMS-weight with modular group} $\sigma$ iff $\varphi
\sigma_t=\varphi$ and $\varphi(a^* a)=
\varphi(\sigma_{\frac{i}{2}}(a)(\sigma_{\frac{i}{2}}(a))^*)$ for any $a\in 
D(\sigma_{\frac{i}{2}})$.\vs\\
{\em $C^*$ dynamical systems and crossed products.} Let $G$ be a locally 
compact group, $A$ a $C^*$ algebra and  $Aut(A)$- group of *-automorphism of 
$A$. Let $G\ni g\mapsto \alpha_g\in Aut(A)$ be a strongly continous group 
homomorphism.  Then {\em a $C^*$dynamical system} is a triple $(G,A,\alpha)$.
A {\em covariant representation} of $(G,A,\alpha)$ is a pair $(\pi,U)$ where
$\pi$ is nondegenerate *-representation of $A$ and $U$ is strongly continous,
unitary representation of $G$ such that $\pi(\alpha_g(x))=U(g)\pi(x)U(g)^{-1}$.
Let $K(G,A)$ denote the linear space of compactly supported continous mapping 
 from $G$ to $A$. It follows that $K(G,A)$ is a normed $*$-algebra if we 
define:
$$(f_1\,f_2)(g):=\int_G \lambda(h)f_1(h)\alpha_h(f_2(h^{-1}g))\,,\,\,
(f*)(g):= \D(g)^{-1}(\alpha_g(f(g^{-1}))^*$$ and 
$||f||_1:=\int_G\lambda(g) ||f(g)||$,
where $\lambda$ is a left Haar measure and $\D$ corresponding modular 
function. Any covariant representation $(\pi,U)$ defines a representation 
$\rho$ of $K(G,A)$ by: $\rho(f):=\int_G \lambda(g) \pi(f(g)) U(g)$. One can 
show that $||\rho(f)||\leq||f||_1$.
It follows that formula $||f||:=sup_{\rho}||\rho(f)||$
where the supremum is taken over all covariant representations 
defines $C^*$ norm. The completion of $K(G,A)$ in this norm is 
{\em the crossed product} $C^*(G,A,\alpha)$.\vs\vs

\ed